\documentclass{amsart}

\usepackage{mathtools}

\usepackage{booktabs}
\usepackage{color}
\usepackage{amsmath}
\usepackage{amsfonts}
\usepackage{amssymb}
\usepackage{latexsym}
\usepackage{array}
\usepackage{amsthm}
\usepackage{graphics}
\usepackage{comment}
\usepackage{epsf,graphicx,subfigure}
\usepackage{epsfig}
\usepackage{epstopdf}
\usepackage{stmaryrd}
\usepackage{tikz}
\usepackage{hyperref}
\usepackage{xcolor}

\usepackage[noadjust]{cite}

\usetikzlibrary{shapes, arrows.meta, positioning}

\tikzstyle{arrow} = [thick, ->, >=stealth]

\tikzstyle{box} = [draw, rectangle, rounded corners, minimum width=1.4cm, minimum height=0.35cm, text centered]

\tikzstyle{box0} = [draw, rectangle, rounded corners, minimum width=1.4cm, minimum height=0.35cm, text centered,text width=3.8cm]

\tikzstyle{box0ex} = [draw, rectangle, rounded corners, minimum width=1.4cm, minimum height=0.35cm, text centered,text width=4.3cm]

\tikzstyle{box0ss} = [draw, rectangle, rounded corners, minimum width=1.4cm, minimum height=0.35cm, text centered,text width=3.2cm]

\tikzstyle{box0s} = [draw, rectangle, rounded corners, minimum width=1.4cm, minimum height=0.35cm, text centered,text width=3.1cm]

\tikzstyle{box1} = [draw, rectangle, rounded corners, minimum width=1.4cm, minimum height=0.35cm, text centered,text width=3.0cm]

\tikzstyle{box2} = [draw, rectangle, rounded corners, minimum width=1.4cm, minimum height=0.35cm, text centered,text width=2.0cm]

\tikzstyle{box3} = [draw, rectangle, rounded corners, minimum width=1.4cm, minimum height=0.35cm, text centered,text width=1.8cm]

\usetikzlibrary{mindmap}

\usetikzlibrary{arrows,positioning,shapes.geometric}

\usepackage{longtable}
\usepackage{rotating}
\usepackage{lscape}
\usepackage{adjustbox}

\usepackage{verbatim}
\usetikzlibrary{positioning,shadows,arrows}

\newtheorem{Proposition}{Proposition}
\newtheorem{Lemma}{Lemma}
\newtheorem{Corollary}{Corollary}
\newtheorem{Assumption}{Assumption}
\newtheorem{Remark}{Remark}
\newtheorem{Definition}{Definition}
\newtheorem{Example}{Example}

\begin{document}

\title[]{On Random Fields Associated with Analytic Wavelet Transform}
\author{Gi-Ren Liu}
\address{Department of Mathematics, National Cheng Kung University, Tainan, Taiwan}
\email{girenliu@gmail.com}

\author{Yuan-Chung Sheu}
\address{Department of Applied Mathematics, National Yang Ming Chiao Tung University, Hsinchu, Taiwan}
\email{sheu@math.nctu.edu.tw}

\author{Hau-Tieng Wu}
\address{Courant Institute of Mathematical Sciences, New York University, New York, NY, 10012 USA}
\email{hauwu@cims.nyu.edu}

\maketitle

\begin{abstract}
Despite the broad application of the analytic wavelet transform (AWT), a systematic statistical characterization of its magnitude and phase as inhomogeneous random fields on the time-frequency domain when the input is a random process remains underexplored.
In this work, we study the magnitude and phase of the AWT as random fields on the time-frequency domain when the observed signal is a deterministic function plus additive stationary Gaussian noise. We derive their marginal and joint distributions, establish concentration inequalities that depend on the signal-to-noise ratio (SNR), and analyze their covariance structures. Based on these results, we derive an upper bound on the probability of incorrectly identifying the time-scale ridge of the clean signal, explore the regularity of scalogram contours, and study the relationship between AWT magnitude and phase. Our findings lay the groundwork for developing rigorous AWT-based algorithms in noisy environments.
\end{abstract}

\section{Introduction}
The continuous wavelet transform (CWT) is a powerful time-frequency (TF) analysis tool widely used to study nonstationary time series across scientific and engineering disciplines \cite{daubechies1992ten,flandrin1998time,flandrin2018explorations,mallat1999wavelet}. To suppress spectral interference between positive and negative frequency components in real-valued signals, analytic mother wavelets are often used \cite{carmona1997characterization,delprat1992asymptotic,holighaus2019characterization,lilly2010analytic,meynard2018spectral}. The resulting transform, known as the AWT, maps a signal to a complex-valued function on the time-scale domain, yielding rich information about its dynamic structure.
Consider a signal modeled as a nonstationary random process $Y(t)=f(t)+\Phi(t)$, where $t\in \mathbb{R}$ is time, $f$ is a deterministic function representing the clean signal, and $\Phi$ is a mean zero real-valued random process modeling the stochastic part of the signal. The AWT of $Y$ is a complex valued function defined on the time-scale domain $\{(t,s);\,t\in \mathbb{R},\,s>0\}$, denoted as $W_{Y}:\{t\in \mathbb{R},\,s>0\}\to \mathbb{C}$, where $W_Y(t,s)$ is defined by convolving $Y$ with an analytic mother wavelet $\psi$ scaled at $s>0$ and centered at time $t\in \mathbb{R}$.
As a result, $W_{Y}$ is a complex-valued random field on the time-scale domain, and its magnitudes $|W_{Y}|$ and phases $\arg(W_{Y})$, where the magnitude and $\arg$ are applied to $W_Y$ pointwisely, also constitute random fields.
The AWT captures key features of signal dynamics, including instantaneous frequency and amplitude modulation \cite{delprat1992asymptotic}, transient phenomena like singularities and oscillatory bursts  \cite{daubechies1992ten}, and various scale-dependent structures \cite{jaffard1996wavelet}.
Numerous algorithms have been developed to exploit the dynamic information encoded in the AWT. For example, the magnitude of AWT, called scalogram, and its contour plot and ridges are usually applied to visualize the dynamics of the signal and handcrafted features \cite{mallat1999wavelet,mallat2002singularity,lilly2010analytic}.
The filter property of the AWT has broad applications, including, for example, the quantification of phase-amplitude coupling and phase-phase coupling in brain waves \cite{tort2010measuring,hulsemann2019quantification}, among other uses. The magnitude information in $|W_Y|$ underpins the scattering transform, a tool for convolutional neural network analysis \cite{mallat2012group,bruna2015intermittent}, while the phase  $\arg(W_{Y})$ is central to reassignment \cite{auger2002improving} and synchrosqueezing \cite{2011synchrosqueezed} techniques that sharpen time-frequency resolution.

Despite its broad applicability, interpreting the AWT requires care, particularly in the presence of noise. Noise distorts the time-scale representation and introduces uncertainty in extracted features. A rigorous statistical understanding of $W_Y$, especially the behavior of $|W_{Y}|$ and $\arg(W_{Y})$, is essential for reliable inference.
Yet, a systematic study of $|W_{Y}|$ and $\arg(W_{Y})$ as two random fields on the time-scale domain, and their interaction, particularly under the presence of a nonzero signal component (i.e., the non-null setting), remains largely underexplored.
Since $W_Y$ is trivially a complex Gaussian random field when $\Phi$ is a stationary Gaussian process, it may seem at first glance that characterizing the distributions of $|W_{Y}|$ and $\arg(W_{Y})$, as well as their interactions, is trivial. However, the affine structure of the AWT, the nonlinearities introduced by taking magnitude and phase, and the presence of a deterministic signal component in the non-null setting, render the problem substantially nontrivial.

In this work,
we study the statistical properties of $|W_{Y}|$ and $\arg(W_{Y})$ when $f$ satisfies appropriate regularity conditions and $\Phi$ is a stationary Gaussian process \cite{major1981lecture} with mild conditions.
Our contributions are multifaceted.
First, we derive the joint and marginal probability density functions of $|W_{Y}(t,s)|$ and $\arg(W_{Y}(t,s))$ for any $(t,s)\in \mathbb{R}\times (0,\infty)$.
These results yield concentration inequalities quantifying how rapidly the magnitudes and phases of the AWT of $Y$ converge to those of  $f$ as the SNR increases, and how the convergence depends on time and scale.
Second, we analyze the covariance structure of $|W_{Y}|$ and $\arg(W_{Y})$.
In the null case $f=0$, we derive explicit expressions for the covariances of $|W_{Y}|$ and $\arg(W_{Y})$. Using Kolmogorov's consistency theorem, we show that local covariance structures across discrete time-scale points extend to valid global processes. This provides a theoretical foundation for empirical covariance estimation and supports statistical procedures such as ridge detection and uncertainty quantification. Furthermore, we show the independence of $|W_{Y}|$ and $\arg(W_{Y})$ at different time-scale points.
Third, we derive an SNR-dependent upper bound on the probability of misidentifying the time-scale ridge of $f$ using $|W_Y|$, extending our previous ridge analysis in \cite{liu2024analyzing}, and study the regularity of scalogram contours, showing that almost all contours are differentiable.
Finally, to bridge continuous-time AWT theory with practical applications,
we introduce a discretized version of the AWT for signals observed through sampling, using a Riemann-sum
approximation. As the sampling rate increases, this discretized transform converges
 to the continuous-time AWT when the underlying signal has a continuous
 sample path. In practice, when the mother wavelet is of infinite support, the Riemann sum must be truncated in the numerical implementation. Intuitively, such truncation error tends to zero as the truncation threshold
 increases. We show that this intuition is correct provided that the wavelet is integrable.
We further establish a Gaussian approximation theorem for the discretized AWT of non-Gaussian causal non-stationary time series, showing that its distribution can be approximated by that of a suitable Gaussian sequence, with an explicit asymptotic upper bound on their L$\acute{\textup{e}}$vy-Prokhorov distance.
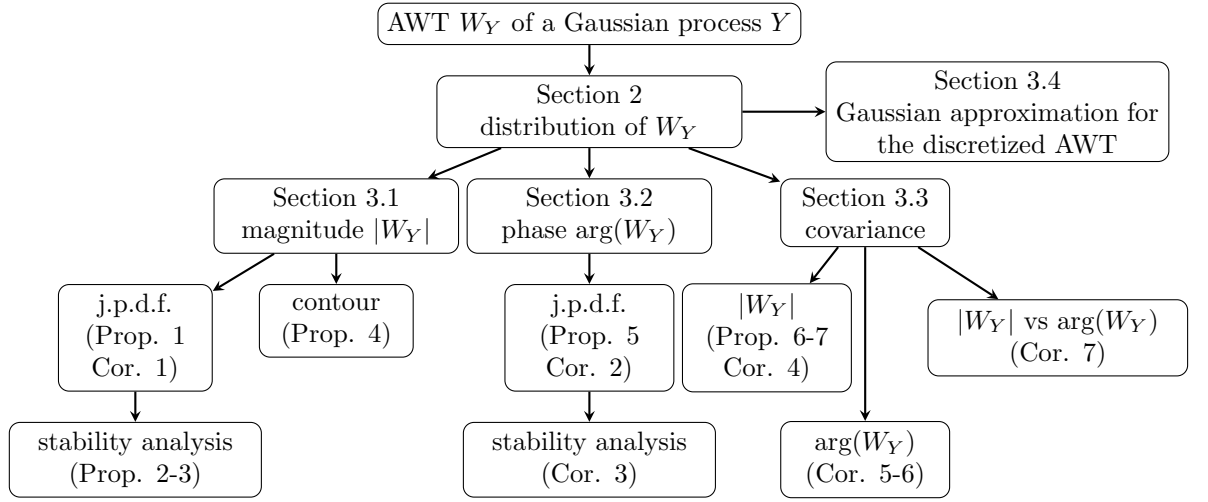
\begin{figure}[htb!]
\begin{center}
\begin{tikzpicture}[node distance=0.9cm and 1.2cm]

\node (root) [box] {AWT $W_Y$ of a Gaussian process $Y$};

\node (circular) [box0, below= 0.4cm of root] {Section \ref{sec:preliminary}\\ distribution of $W_{Y}$};

\node (gaussapprox) [box0ex, right=1.1cm of circular]
    {Section \ref{sec:Gaussian_approxi}\\Gaussian approximation for \\the discretized AWT};

\node (phase) [box1, below = 0.4cm of circular] {Section \ref{sec:phase}\\phase $\arg(W_Y)$};
\node (mag) [box1, left=0.1cm of phase] {Section \ref{sec:main:magnitude} magnitude $|W_Y|$ };

\node (covar) [box2, below right=0.4cm and 0.5cm of circular] {Section \ref{sec:cov}\\covariance};

\node (jpdf_mag) [box3, below left=0.4cm and 0.01cm of mag] {j.p.d.f.\\(Prop. \ref{prop:jpdf:Sn}\\Cor.  \ref{prop:jpdf:S2})};
\node (level_set) [box3, below= 0.4cm of mag] {contour\\(Prop. \ref{thm:contour:c1})};
\node (jpdf_phase) [box2, below=0.4 cm of phase] {j.p.d.f.\\(Prop. \ref{prop:jpdf:Thetan}\\Cor. \ref{corollary:phase:jpdf:null})};

\node (stab_mag) [box0s, below=0.4 cm of jpdf_mag] {stability analysis\\(Prop. \ref{lemma:convprob}-\ref{prop:ridge})};
\node (stab_phase) [box0s, below=0.4cm of jpdf_phase] {stability analysis\\(Cor. \ref{corollary:phase:expo})};


\node (magvsmag) [box2, right=0.1cm of jpdf_phase] {$|W_Y|$\\(Prop. \ref{cov:mag:generalpsi}-\ref{prop:Cov_general}\\Cor. \ref{prop:Cov_null})};

\node (phasevsphase) [box2, below=2.3cm  of covar] {$\arg(W_Y)$\\(Cor. \ref{prop:Cov2_phase_null}-\ref{prop:Cov_phase_null})};
\node (magvsphase) [box0ss, right=1.0cm of magvsmag] {$|W_Y|$ vs $\arg(W_Y)$\\(Cor. \ref{prop:Cov_phaseAM_null})};

\draw[arrow] (root) -- (circular);
\draw[arrow] (circular) -- (mag);
\draw[arrow] (circular) -- (phase);
\draw[arrow] (circular) -- (covar);

\draw[arrow] (mag) -- (jpdf_mag);
\draw[arrow] (mag) -- (level_set);
\draw[arrow] (jpdf_mag) -- (stab_mag);

\draw[arrow] (phase) -- (jpdf_phase);
\draw[arrow] (jpdf_phase) -- (stab_phase);

\draw[arrow] (covar) -- (magvsmag);
\draw[arrow] (covar) -- (phasevsphase);
\draw[arrow] (covar) -- (magvsphase);

\draw[arrow] (circular) -- (gaussapprox);

\end{tikzpicture}
\end{center}
\caption{Flowchart of the structure of this paper. j.p.d.f: joint probability density function.}
\label{fig:flowchart:structure}
\end{figure}

The rest of the paper is organized as follows.
In Section \ref{sec:preliminary}, we present some necessary background, including a review of the continuous wavelet transform and the spectral representation of stationary Gaussian processes.
As illustrated in Figure~\ref{fig:flowchart:structure}, our main findings are presented in Section \ref{sec:mainresult}, which is further divided into four subsections addressing: (1) the distributional properties of the magnitude of AWT, (2) the distributional properties of the phase of AWT, (3) the covariance structure of the AWT's magnitude and phase,
and (4)
the Gaussian approximation for the discretized AWT of non-Gaussian and non-stationary noise.
Section \ref{sec:conclusion} concludes the paper, and Section \ref{sec:proof} contains all proofs.

\section{Preliminary}\label{sec:preliminary}

Let $\psi: \mathbb{R}\rightarrow  \mathbb{C}$ be a function in $L^{1}(\mathbb{R})\cap L^{2}(\mathbb{R})$.
Denote the Fourier transform of $\psi$ by
$\widehat{\psi}(\omega) = \int_{\mathbb{R}}e^{-i t\omega}\psi(t)dt,\ \omega\in \mathbb{R}$.
\begin{Assumption}\label{assump:analytic}
The function $\psi$ is an analytic mother wavelet, which means that
$\widehat{\psi}(\omega)=0$ for $\omega\leq0$.
Furthermore, $|\widehat{\psi}|$ is assumed to be unimodal with a unique maximum at
\[
\omega_{\psi}=\frac{1}{2\pi}\operatorname*{arg\,max}_{\omega>0}|\widehat{\psi}(\omega)|,
\]
which we refer to as the center frequency of $\psi$.
\end{Assumption}
Consider a real-valued signal $f$, which we treat as a bounded function.
The CWT of $f$, denoted by $\{W_{f}(t,s)\mid t\in \mathbb{R}, s>0\}$, is defined as follows \cite{daubechies1992ten}:
\begin{align}\label{def:WT}
W_{f}(t,s)
=\int_{\mathbb{R}} f(\tau) \frac{1}{\sqrt{s}}\overline{\psi\left(\frac{\tau-t}{s}\right)} d\tau,
\end{align}
where $t\in \mathbb{R}$ represents time, $s>0$ represents scale.
Under Assumption \ref{assump:analytic}, $W_{f}$ is often referred to as the AWT of $f$.
In practice, signals are often contaminated by noise. Thus, we model the recorded signal as
\begin{align*}
Y(t) = f(t)+\Phi(t),\ t\in \mathbb{R},
\end{align*}
where $\{\Phi(t)\mid t\in \mathbb{R}\}$ represents the noise.
Similar to (\ref{def:WT}), we denote the AWT of the noise $\Phi$ and the noisy signal $Y$ by $W_{\Phi}$ and $W_{Y}$,
respectively.
Given the relative breadth of analytical tools for Gaussian processes \cite{breuer1983central,clausel2012large, dobrushin1979non,ivanov2013limit,taqqu1979convergence} and their ubiquity in applications \cite{makowiec2006long}, we adopt the following assumption.
\begin{Assumption}\label{assump:Gaussian}
The noise $\Phi$ is a stationary Gaussian process with a constant mean and
a covariance function $C_{\Phi}$:
\begin{equation*}
C_{\Phi}(t_{1}-t_{2}) = \mathbb{E}[\Phi(t_{1})\Phi(t_{2})]- \mathbb{E}[\Phi(t_{1})] \mathbb{E}[\Phi(t_{2})],\ t_{1},t_{2}\in \mathbb{R}.
\end{equation*}
\end{Assumption}
Under Assumption \ref{assump:Gaussian}, by the Bochner-Khinchin theorem \cite[Chapter 4]{krylov2002introduction},
there exists a unique nonnegative measure $F$ on $\mathbb{R}$ such that $F(\mathbb{R}) = C_{\Phi}(0)$
and the covariance function of $\Phi$ has the spectral representation
\begin{equation*}
C_{\Phi}(t) = \int_{\mathbb{R}}e^{i\lambda t}F(d\lambda),\ t\in \mathbb{R}.
\end{equation*}
The measure $F$ is called the spectral measure of the covariance function $C_{\Phi}$.
Furthermore, by the Karhunen Theorem, the process $\Phi$ has the representation
\begin{align}\label{spectral_rep:samplepath}
\Phi(t)=\mathbb{E}[\Phi(t)]+\int_{\mathbb{R}}
e^{i\lambda t}Z_{F}(d\lambda),\ t\in \mathbb{R},
\end{align}
where $Z_{F}$ is a complex-valued Gaussian random measure on $\mathbb{R}$
satisfying
\begin{align*}
Z_{F}(\Delta_{1})=\overline{Z_{F}(-\Delta_{1})}
\end{align*}
and
\begin{align}\label{ortho}
\mathbb{E}\left[Z_{F}(\Delta_{1})
\overline{Z_{F}(\Delta_{2})}\right]=F(\Delta_{1}\cap\Delta_{2})
\end{align}
for any measurable sets
 $\Delta_{1},\Delta_{2}$ in the $\sigma$-algebra of $\mathbb{R}$.
See, for example, \cite{major1981lecture} for the above facts.
\begin{Lemma}\label{prop:circular}
For any $n\in \mathbb{N}$, $t_{1},t_{2},...,t_{n}\in \mathbb{R}$ and $s_{1},s_{2},...,s_{n}>0$,
denote
\begin{align}\label{vector:WY}
\mathbf{W}_{Y} = \left[W_{Y}(t_{1},s_{1})\ W_{Y}(t_{2},s_{2})\ \cdots\ W_{Y}(t_{n},s_{n})\right]^{\top},
\end{align}
\begin{align}\label{def:Gamma}
\Gamma =& \mathbb{E}\left[\left(\mathbf{W}_{Y}-\mathbb{E}\left[\mathbf{W}_{Y}\right]\right)
\left(\mathbf{W}_{Y}-\mathbb{E}\left[\mathbf{W}_{Y}\right]\right)^{*}
\right],
\end{align}
and
\begin{align*}
\mathbf{C} =& \mathbb{E}[\left(\mathbf{W}_{Y}-\mathbb{E}\left[\mathbf{W}_{Y}\right]\right)
\left(\mathbf{W}_{Y}-\mathbb{E}\left[\mathbf{W}_{Y}\right]\right)^{\top}
],
\end{align*}
where $*$ represents the conjugate transpose, and $\top$ denotes the transpose.
Under Assumption \ref{assump:Gaussian}, for any $\ell,\ell'\in \{1,...,n\}$, the $(\ell,\ell')$-th entries of $\Gamma$
and $\mathbf{C}$ have the representation
\begin{align}\label{entry_Gamma}
\Gamma_{\ell\ell'}= \sqrt{s_{\ell}s_{\ell'}}\int_{\mathbb{R}}e^{i(t_{\ell}-t_{\ell'})\lambda}\overline{\widehat{\psi}(s_{\ell}\lambda)}
\widehat{\psi}(s_{\ell'}\lambda)F(d\lambda)
\end{align}
and
\begin{align}\label{entry_pseudoC}
C_{\ell\ell'}= \sqrt{s_{\ell}s_{\ell'}}\int_{\mathbb{R}}e^{i(t_{\ell}-t_{\ell'})\lambda}\overline{\widehat{\psi}(s_{\ell}\lambda)}
\overline{\widehat{\psi}(-s_{\ell'}\lambda)}F(d\lambda),
\end{align}
where $F$ is the spectral measure of $C_{\Phi}$.
Furthermore, if Assumption \ref{assump:analytic} holds,
$\mathbf{C}=\mathbf{0}_{n\times n}$.
\end{Lemma}
The proof of Lemma \ref{prop:circular} is provided in Section \ref{sec:proof:prop:circular}.
Because $\mathbb{E}[\mathbf{W}_{Y}]=\mathbb{E}[\mathbf{W}_{f}]$, where
\begin{align*}
\mathbf{W}_{f} = \left[W_{f}(t_{1},s_{1})\ W_{f}(t_{2},s_{2})\ \cdots\ W_{f}(t_{n},s_{n})\right]^{\top},
\end{align*}
Lemma \ref{prop:circular} indicates that $\mathbf{W}_{Y}-\mathbf{W}_{f}$ is a circularly symmetric Gaussian random vector under Assumptions \ref{assump:analytic} and \ref{assump:Gaussian}.
In this context, the probability density function of the complex Gaussian vector $\mathbf{W}_{Y}$
is completely characterized by its mean $\mathbf{W}_{f}$ and covariance matrix $\Gamma$ as follows
\cite{schreier2010statistical}:
\begin{align}\label{pdf_WY}
p\left(\mathbf{w};\mathbf{W}_{Y}\right)
= \frac{1}{\pi^{n}\textup{det}(\Gamma)}
\textup{exp}\left(-\left(\mathbf{w}-\mathbf{W}_{f}\right)^{*}\Gamma^{-1}\left(\mathbf{w}-\mathbf{W}_{f}\right)\right),
\end{align}
where $\mathbf{w}\in \mathbb{C}^{n}$.
Hereafter, we use the notation $p(\cdot;X)$ to represent the probability density function of a continuous random variable or vector $X$.
Note that (\ref{pdf_WY}) means that the real part $\mathfrak{R}(\mathbf{W}_{Y})$
and the imaginary part $\mathfrak{I}(\mathbf{W}_{Y})$ of $\mathbf{W}_{Y}$ have the joint probability density function
\begin{align}\notag
p\left(\left[\begin{array}{c}\mathbf{x}\\\mathbf{y}\end{array}\right];\left[\begin{array}{c}\mathfrak{R}(\mathbf{W}_{Y})\\
\mathfrak{I}(\mathbf{W}_{Y})\end{array}\right]\right)
=p\left(\mathbf{x}+i\mathbf{y};\mathbf{W}_{Y}\right),\ \mathbf{x},\mathbf{y}\in \mathbb{R}^{n}.
\end{align}
For any $t\in \mathbb{R}$ and $s>0$, we define the phase of $W_{Y}(t,s)$ as
$$\Theta_{Y}(t,s) = \arg(W_{Y}(t,s))\ \textup{mod}\ 2\pi.$$
That is,
\begin{align}\label{def:Theta_Y}
W_{Y}(t,s) = |W_{Y}(t,s)|e^{i\Theta_{Y}(t,s)}.
\end{align}
In this work, we use $[0,2\pi)$ as the representative interval for the quotient space $\mathbb{R}/2\pi \mathbb{Z}$,
so that $\Theta_{Y}(t,s)\in[0,2\pi).$

\begin{Remark}
From (\ref{pdf_WY}), the joint probability density function of $|W_{Y}(t,s)|$ and $\Theta_{Y}(t,s)$ is given by
 \begin{align}\label{jpdf_|W|Theta}
&p\left(r,\theta;|W_{Y}(t,s)|,\Theta_{Y}(t,s)\right)
\\\notag= & \frac{r}{\pi\mathbb{E}[|W_{\Phi}(t,s)|^2]}\textup{exp}\left(-\frac{|W_{f}(t,s)|^2}{\mathbb{E}[|W_{\Phi}(t,s)|^2]}\right)
\\\notag&\times\textup{exp}\left(-\frac{1}{\mathbb{E}[|W_{\Phi}(t,s)|^2]}\left(r^2-2r\Re(W_{f}(t,s)e^{-i\theta})\right)
\right)
\end{align}
for $r\geq0$ and $\theta\in [0,2\pi)$. Unlike the null case where $f=0$, the nonseparability of the joint probability density function in (\ref{jpdf_|W|Theta})
indicates that the magnitude $|W_{Y}(t,s)|$ and the phase $\Theta_{Y}(t,s)$ are dependent when $W_{f}(t,s)\neq0$.
\end{Remark}

\section{Main results}\label{sec:mainresult}
This section presents our main findings: the statistical properties of random fields $|W_{Y}|$ and $\arg(W_{Y})$.

\subsection{Marginal distributions of AWT magnitude}\label{sec:main:magnitude}
In the following, we denote the componentwise magnitude of the random vector $\mathbf{W}_{Y}$, as given in (\ref{vector:WY}), by
\begin{align*}
|\mathbf{W}_{Y}|:=\left[|W_{Y}(t_{1},s_{1})|\ \ |W_{Y}(t_{2},s_{2})|\ \ \cdots\ \ |W_{Y}(t_{n},s_{n})|\right]^{\top}.
\end{align*}

\begin{Definition}
Let $\{I_{\nu}\}_{\nu\in \mathbb{R}}$ denote the modified Bessel functions of the first kind,
defined by
\begin{equation}\label{def:;Bessel}
I_{\nu}(x) = \left(\frac{z}{2}\right)^{\nu} \overset{\infty}{\underset{k=0}{\sum}}\frac{(x^{2}/4)^{k}}{k!\Gamma(\nu+k+1)},\ x\in \mathbb{R},
\end{equation}
where $\Gamma(\cdot)$ is the Gamma function.
The probability density function of the Rice distribution \cite{miller1964multidimensional} with noncentrality parameter $m\geq0$ and spread parameter $\sigma^{2}>0$ is given by
\begin{align}\label{def:rice}
p_{\textup{Rice}}\left(r;m,\sigma^{2}\right) = \frac{r}{\sigma^{2}}
\textup{exp}\left(-\frac{r^2+m^{2}}{2\sigma^{2}}\right)
I_{0}\left(\frac{mr}{\sigma^{2}}\right)
\end{align}
for $r\geq 0$.
\end{Definition}

The Rice distribution arises naturally in signal processing and probability. For example, when we model the magnitude of a noisy complex-valued sinusoidal signal, the parameter $m$ reflects the strength of the underlying signal, $\sigma^2$ reflects the noise strength, and $r$ is the observed signal magnitude, which is inherently random.

\begin{Proposition}\label{prop:jpdf:Sn}
Assume that Assumptions \ref{assump:analytic} and \ref{assump:Gaussian} hold.
Additionally,
suppose that the covariance matrix $\Gamma$ of $\mathbf{W}_{Y}$
is invertible.

(a) The probability density function of the random vector $|\mathbf{W}_{Y}|$
can be expressed as
\begin{align}\label{jpdf:|WY|:n:nonnull}
&p\left(r_{1:n};|\mathbf{W}_{Y}|\right)
=\left(\overset{n}{\underset{\ell=1}{\prod}}r_{\ell}\right)\int_{0}^{2\pi}\hspace{-0.15cm}\int_{0}^{2\pi}\hspace{-0.2cm}\cdots \hspace{-0.15cm}\int_{0}^{2\pi}
\hspace{-0.1cm} p\left(\mathbf{w};\mathbf{W}_{Y}\right)d\theta_{1}d\theta_{2}\cdots d\theta_{n}
\end{align}
for $r_{1:n}=(r_{1},r_{2},\ldots,r_{n})\in \mathbb{R}_{\geq 0}^{n}$, where
$\mathbf{w}=[ r_{1}e^{i\theta_{1}} \  r_{2}e^{i\theta_{2}} \  \cdots \  r_{n}e^{i\theta_{n}}]^{\top}$
and $p\left(\mathbf{w};\mathbf{W}_{Y}\right)$ is defined in
(\ref{pdf_WY}).

(b) For any $t\in \mathbb{R}$ and $s>0$,
\begin{align*}
p\left(r;|W_{Y}(t,s)|\right)
=p_{\textup{Rice}}\left(r;|W_{f}(t,s)|,\frac{1}{2}\mathbb{E}[|W_{\Phi}(t,s)|^{2}]\right)
\end{align*}
for $r\geq0$, where $p_{\textup{Rice}}$ is the probability density function of the Rice distribution defined in (\ref{def:rice}).

(c) For the case $n=2$,  $p\left(r_{1:n};|\mathbf{W}_{Y}|\right)$ can be expressed as a single integral:
\begin{align}\notag
&p\left(r_{1},r_{2}; |W_{Y}(t_{1},s_{1})|,|W_{Y}(t_{2},s_{2})|\right)
\\\notag=&\frac{r_{1}}{\pi\Gamma_{11}}\int_{0}^{2\pi}
 p_{\textup{Rice}}\left(r_{2};m,
\frac{1}{2}\left(\Gamma_{22}-\frac{|\Gamma_{12}|^{2}}{\Gamma_{11}}\right)\right)
\\\label{jpdf:|WY|:n=2:nonnull}&\times\textup{exp}\left(-\frac{|r_{1}e^{i\theta}-W_{f}(t_{1},s_{1})|^{2}}{\Gamma_{11}}\right)d\theta,
\end{align}
where
$$
m=\left|\frac{\Gamma_{12}}{\Gamma_{11}}\left(r_{1}e^{i\theta}-W_{f}(t_{1},s_{1})\right)
+W_{f}(t_{2},s_{2})\right|.
$$
\end{Proposition}
The proof of Proposition \ref{prop:jpdf:Sn} is provided in Section \ref{sec:proof:prop:jpdf:Sn}.
When the magnitude of the AWT of the clean signal $f$ tends to zero, the above joint probability density functions reduce to those given in the following corollary.

\begin{Corollary}\label{prop:jpdf:S2}
Assume that Assumptions \ref{assump:analytic} and \ref{assump:Gaussian} hold.
Additionally,
suppose that the covariance matrix $\Gamma$ of $\mathbf{W}_{Y}$
is invertible.

(a) In the null case $f=0$, the joint probability density function $p\left(r_{1:n};|\mathbf{W}_{Y}|\right)$,
given in (\ref{jpdf:|WY|:n:nonnull}),
simplifies to the following form
\begin{align}\label{jpdf:S:n}
p\left(r_{1:n};|\mathbf{W}_{Y}|\right)
=
\frac{1}{\pi^{n}\textup{det}(\Gamma)}\left(\overset{n}{\underset{\ell=1}{\prod}}r_{\ell}\right)
J\left(\textup{diag}\left(r_{1:n}\right)\Gamma^{-1}\textup{diag}\left(r_{1:n}\right)\right),
\end{align}
where $J: \mathbb{R}^{n\times n} \mapsto [0,\infty)$ is a function defined for any $n\times n$ Hermitian matrix $\mathbf{M}$ by
\begin{align}\label{def:functionB}
J(\mathbf{M}) = \int_{0}^{2\pi}\hspace{-0.15cm}\int_{0}^{2\pi}\hspace{-0.17cm}\cdots \hspace{-0.12cm}\int_{0}^{2\pi}\hspace{-0.2cm} \exp
\left(-\mathbf{v}^{*}\mathbf{M}\mathbf{v}\right)d\theta_{1}d\theta_{2}\cdots d\theta_{n}
\end{align}
with $\mathbf{v}= [e^{i\theta_{1}} \ e^{i\theta_{2}} \ \cdots \ e^{i\theta_{n}}]^{\top}$.

(b) In the null case with $n=2$,
(\ref{jpdf:S:n})
can be further
simplified as follows
\begin{align}\label{jpdf:|WY|:n=2:null}
&p\left(r_{1:2};|\mathbf{W}_{\Phi}|\right)
=\frac{4r_{1}r_{2}}{\textup{det}(\Gamma)}
\exp\left(-\left(\frac{r_{1}^2}{\Gamma_{11}}+\frac{r_{2}^2}{\Gamma_{22}}\right)
\frac{\Gamma_{11}\Gamma_{22}}{\textup{det}(\Gamma)}\right)
 I_{0}\left(
\frac{2|\Gamma_{12}|r_{1}r_{2}} {\textup{det}(\Gamma)}\right),
\end{align}
where $\textup{det}(\Gamma)=\Gamma_{11}\Gamma_{22}-|\Gamma_{12}|^2$, and
$\Gamma_{11}$, $\Gamma_{22}$, and $\Gamma_{12}$ are defined in (\ref{entry_Gamma}).
\end{Corollary}
The proof of Corollary \ref{prop:jpdf:S2} is provided in Section \ref{sec:proof:prop:jpdf:S2}.

\begin{Remark}
By the formula \cite[eq. 8.976]{gradshteyn2014table}:
\begin{align*}
\sum_{k=0}^{\infty} z^k L_k(x) L_k(y)
=
\frac{1}{1 - z}
\exp\left( -\frac{z(x + y)}{1 - z} \right)
I_{0}\left( \frac{2\sqrt{zxy}}{1 - z} \right)
\end{align*}
for $z\in[0,1)$,
where $\{L_{k}\}_{k=0,1,2,...}$ are the Laguerre polynomials,
the probability density function (\ref{jpdf:|WY|:n=2:null})
can also be expressed in terms of a series as follows
\begin{align*}
p\left(r_{1},r_{2};|\mathbf{W}_{\Phi}|\right)  =&\frac{4r_{1}r_{2}}{\Gamma_{11}\Gamma_{22}}\exp\left(-\left(\frac{r_{1}^2}{\Gamma_{11}}
+\frac{r_{2}^2}{\Gamma_{22}}\right)\right)
\\&\times\overset{\infty}{\underset{k=0}{\sum}}\left(\frac{|\Gamma_{12}|^2}{\Gamma_{11}\Gamma_{22}}\right)^{k}
L_{k}\left(\frac{r_{1}^2}{\Gamma_{11}}\right)
L_{k}\left(\frac{r_{2}^2}{\Gamma_{22}}\right).
\end{align*}
\end{Remark}



We now return our attention to the non-null case $f\neq0$, focusing on the situation where the clean signal is stronger than the noise. From Proposition \ref{prop:jpdf:Sn}(b), we immediately have
\begin{align}\label{pdf:magnitude}
&p\left(r;\frac{|W_{Y}(t,s)|}{|W_{f}(t,s)|}\right)
=
2q(t,s)re^{-q(t,s)(1+r^2)}
I_{0}\left(2q(t,s)r\right)
\end{align}
for $r\geq0$, where  $q:\mathbb{R}\times  \mathbb{R}_{+} \rightarrow [0,\infty)$ is defined as
\begin{equation}\label{def:SNR}
q(t,s)=\frac{|W_{f}(t,s)|^2}{\mathbb{E}[|W_{\Phi}(t,s)|^2]}.
\end{equation}
Since $q(t,s)$ is the ratio of the energy of the clean signal $f$ to that of the noise $\Phi$ around time $t$
measured by a wavelet at scale $s$,  we also refer to it as the SNR, which is a measure of signal quality in the wavelet domain.
About the expression in (\ref{pdf:magnitude}), note that the modified Bessel function of the first kind has the asymptotic behavior
\begin{align}\label{approx:Bessel}
I_{0}(x) = \frac{e^{x}}{\sqrt{2\pi x}}\left(1+\mathcal{O}\left(\frac{1}{x}\right)\right)\ \textup{as}\ x\rightarrow\infty.
\end{align}
When $q(t,s)$ is sufficiently large, by using (\ref{approx:Bessel}), for each $r>0$ we have
\begin{align*}
p\left(r;\frac{|W_{Y}(t,s)|}{|W_{f}(t,s)|}\right)
=&
2q(t,s)r e^{-q(t,s)(1+r^2)}
 \frac{e^{2q(t,s)r}}{\sqrt{4\pi q(t,s)r}}\left(1+\mathcal{O}\left(\frac{1}{2q(t,s)r}\right)\right)
\\=&\sqrt{\frac{q(t,s)r}{\pi}}e^{-q(t,s)\left(1-r\right)^{2}}\left(1+\mathcal{O}\left(\frac{1}{2q(t,s)r}\right)\right)\,.
\end{align*}
This implies that the distribution of $|W_{Y}(t,s)|/|W_{f}(t,s)|$ will become increasingly concentrated around one as the SNR $q(t,s)$ grows larger.
The following proposition describes the concentration speed, whose proof is provided in Section \ref{sec:proof:lemma:convprob}.

\begin{Proposition}\label{lemma:convprob}
Assume that Assumptions \ref{assump:analytic} and \ref{assump:Gaussian} hold.
For any $t\in \mathbb{R}$, $s>0$, and $\varepsilon>0$,
\begin{align*}
\mathbb{P}\left(\left|\frac{|W_{Y}(t,s)|}{|W_{f}(t,s)|}-1\right|>\varepsilon\right)
\leq\textup{exp}\left(-\frac{\varepsilon^{2}}{16}q(t,s)\right)
+\textup{exp}\left(-\frac{\varepsilon}{2}q(t,s)\right).
\end{align*}
\end{Proposition}

\subsubsection{About ridge point}
When $f$ is an amplitude- and frequency-modulated signal, we define
\begin{equation*}
s_{f}(t) = \textup{arg}\underset{s>0}{\max}\ |W_{f}(t,s)|
\end{equation*}
which is the so-called {\em ridge point} in the time-scale representation of $f$ \cite{lilly2010analytic}, and
can be used to estimate the instantaneous frequency of $f$ \cite{delprat1992asymptotic}.
In the presence of noise $\Phi$, for a given time $t$, the event $\{|W_{Y}(t,s_{f}(t))|< |W_{Y}(t,s)|\}$ may occur for some $s\neq s_{f}(t)$,
In such cases, the Fourier frequency corresponding to scale $s$ may be mistakenly interpreted as the instantaneous frequency of $f$ at that time. The probability of this event is estimated as follows, which slightly extends the results given in \cite{liu2024analyzing}.

\begin{Proposition}\label{prop:ridge}
Assume that Assumptions \ref{assump:analytic} and \ref{assump:Gaussian} hold.
Let $f$ be an amplitude- and frequency-modulated signal such that $s_{f}(t)$ is a singleton for a given $t\in \mathbb{R}$.
Then, for any $s>0$ with $s\neq s_{f}(t)$ and any $\delta\in [0,1)$,
\begin{align*}
&\mathbb{P}\left(|W_{Y}(t,s_{f}(t))|<(1-\delta) |W_{Y}(t,s)|\right)
\\
\leq&
\textup{exp}\left(-\frac{\varepsilon^{2}}{16}q(t,s_{f}(t))\right)
+\textup{exp}\left(-\frac{\varepsilon}{2}q(t,s_{f}(t))\right)
\\&+\textup{exp}\left(-\frac{\varepsilon^{2}}{16}q(t,s)\right)
+\textup{exp}\left(-\frac{\varepsilon}{2}q(t,s)\right),
\end{align*}
where
\begin{align}\label{threshold}
\varepsilon=\frac{|W_{f}(t,s_{f}(t))|-(1-\delta)|W_{f}(t,s)|}{|W_{f}(t,s_{f}(t))|+(1-\delta)|W_{f}(t,s)|}.
\end{align}
\end{Proposition}
The proof of Proposition \ref{prop:ridge} can be found in Section \ref{sec:proof:prop:ridge}.
This proposition describes how the ridge point is impacted by SNR.
The chosen mother wavelet $\psi$ typically has a Fourier transform $\widehat{\psi}(\zeta)$
that is concentrated around a central frequency and decays rapidly as
$\zeta$ moves away from it.
As a result of this spectral localization, $|W_f(t, s)|$ decreases as $s$ deviates from $s_{f}(t)$,
causing the deviation measure $\varepsilon$ to increase.
When $s_f(t)$ is far from $s$ so that $W_f(t,s)$ is small, we have $\varepsilon\approx 1$.
In this case, assuming $q(t,s)$ is not too small,
the probability of the event $\{|W_{Y}(t,s_{f}(t))|< |W_{Y}(t,s)|\}$ is small.
On the other hand, when $s$ is close to $s_f(t)$, we have $|W_f(t,s)|\approx |W_f(t,s_f(t))|$, implying $\varepsilon\approx 0$. In this case, the event $\{|W_{Y}(t,s_{f}(t))|< |W_{Y}(t,s)|\}$ remains unlikely provided that both $q(t,s)$ and $q(t,s_f(t))$ are sufficiently large. In other words, accurate ridge identification is more likely when the signal is strong across neighboring scales.

\subsubsection{About level contours}
\begin{figure*}[htb!]
\centering
\includegraphics[trim=70 10 70 20,clip,width=\textwidth]{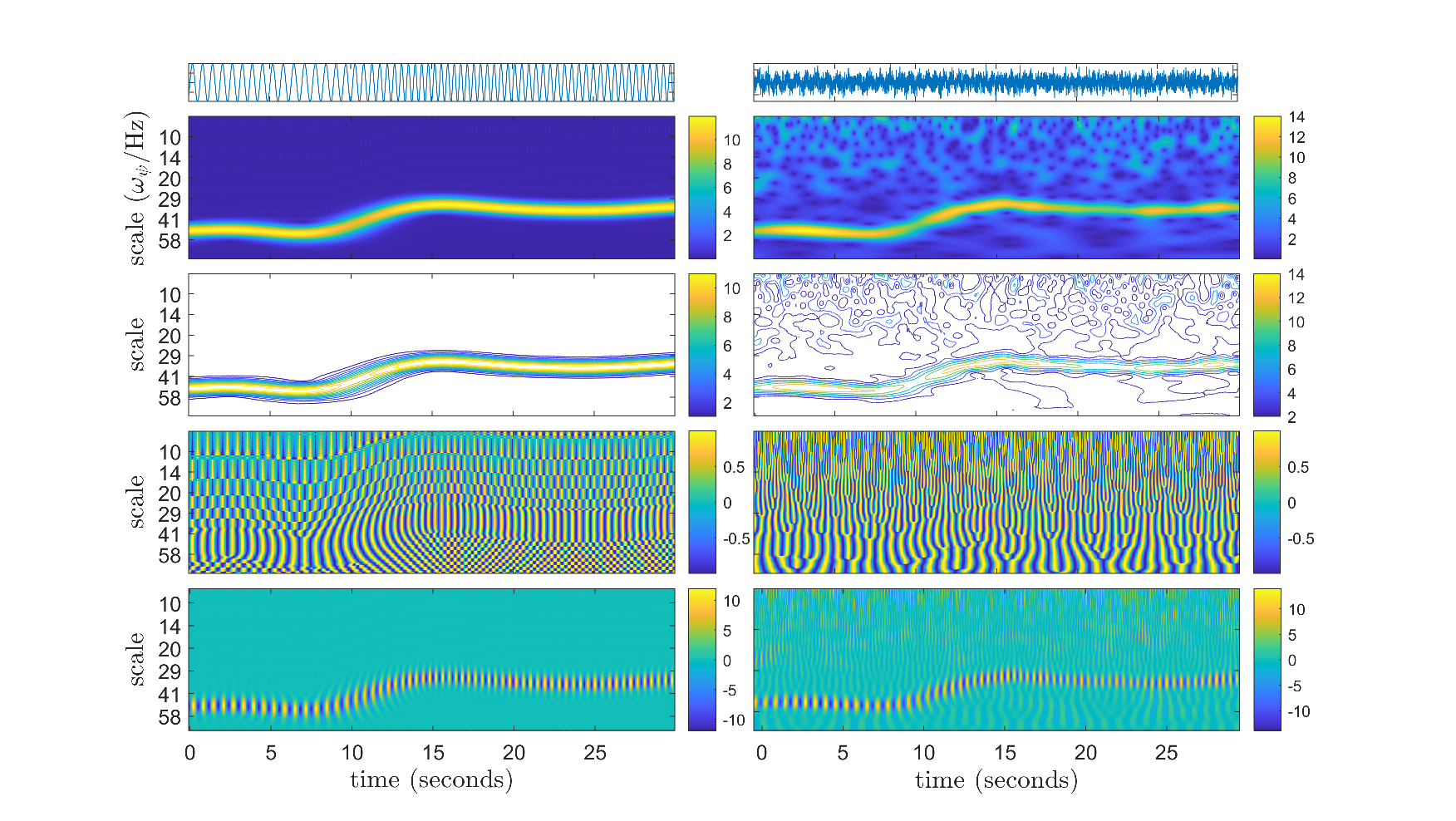}
 \caption{Time-scale representations of the clean signal $f$ (left column) and the noisy signal $Y$ (right column).
 The clean signal $f$, shown in the top-left panel, is a frequency-modulated cosine function.
The noisy signal $Y$, shown in the top-right panel, is generated by adding a sample path of Gaussian noise to $f$.
The second row displays the magnitudes of the AWT of $f$ and $Y$, i.e.,  the scalograms $|W_{f}|$ and $|W_{Y}|$.
Here, the signals $f$ and $Y$ are sampled at 200 Hz, and the mother wavelet $\psi$ has center frequency $\omega_{\psi}= 80$ Hz, as defined in Assumption \ref{assump:analytic}.
The third row shows the level curves of the corresponding scalograms.
The fourth row presents the fields $\cos(\arg(W_{f}))$ and $\cos(\arg(W_{Y}))$, where $\arg(W_{f})$ and $\arg(W_{Y})$
denote the phases of the AWT of $f$ and $Y$, respectively.
The fifth row displays  $|W_{f}|\cos(\arg(W_{f}))$ and $|W_{Y}|\cos(\arg(W_{Y}))$, i.e., the multiplication
of the amplitude and the cosine of the phase fields.
 }
 \label{fig:contour}
\end{figure*}

We conclude this subsection by presenting a result on the local differentiability of level curves in the graph of $|W_{Y}|$.
The third row of Figure~\ref{fig:contour} illustrates the scalogram contours of a clean signal $f$ and its noisy counterpart $Y$.



\begin{Proposition}\label{thm:contour:c1}
Let $f$ be a square-integrable function and let $\Phi$ be a Gaussian process satisfying Assumption \ref{assump:Gaussian}.
Assume further that the Fourier transform $\widehat{\psi}$ of the analytic mother wavelet $\psi$ is
three times differentiable on $(0,\infty)$, and the following conditions hold:
\begin{align*}
&\textup{($\mathrm{D}^{0}_{1}$-$\mathrm{D}^{0}_2$):}\  \underset{\lambda>1}{\sup}|\lambda^{p}\widehat{\psi}(\lambda)|<\infty,\ p=1,2;\ \ \
\\&\textup{($\mathrm{D}^{1}_{0}$-$\mathrm{D}^{1}_2$):}\  \underset{\lambda>0}{\sup}|\lambda^{p}D\widehat{\psi}(\lambda)|<\infty,\ p=0,1,2;
\\
&\textup{($\mathrm{D}^{2}_2$-$\mathrm{D}^{2}_3$):}\  \underset{\lambda>0}{\sup}|\lambda^{p}D^{2}\widehat{\psi}(\lambda)|<\infty,\ p=2,3;\ \ \
\\&\textup{($\mathrm{D}^{3}_3$):}\  \underset{\lambda>0}{\sup}|\lambda^{3}D^{3}\widehat{\psi}(\lambda)|<\infty,
\end{align*}
where $D^{k}$ represents the $k$th derivative operator.
Then, with probability one, there exists a subset $I\subset \mathbb{R}_{+}$ such that
$\mathbb{R}_{+}\setminus I$ has Lebesgue measure zero, and for any $c\in I$ and any point $(t_{0},s_{0})\in \mathbb{R}\times \mathbb{R}_{+}$ satisfying
$$
|W_{Y}(t_0,s_0)| = c,
$$
one of the following results holds:
\begin{enumerate}
\item There exists a differentiable function $\varphi_{1}$, defined in a neighborhood of $t_0$ with $\varphi_{1}(t_{0})=s_{0}$, and some $\varepsilon>0$ such that
$$
|W_{Y}(t,\varphi_{1}(t))| = c
$$
for all $t\in (t_{0}-\varepsilon,t_{0}+\varepsilon)$.

\item There exists a differentiable function $\varphi_{2}$, defined in a neighborhood of $s_0$ with $\varphi_{2}(s_{0})=t_{0}$, and some $\varepsilon>0$ such that
$$
|W_{Y}(\varphi_{2}(s),s)| = c
$$
for all $s\in (s_{0}-\varepsilon,s_{0}+\varepsilon)$.
\end{enumerate}
Furthermore, for any $c\in I$, the level set
$L_c = \left\{(t, s) \in \mathbb{R}\times \mathbb{R}_{+}  \mid |W_Y(t, s)| = c\right\}$
is an embedded submanifold of $\mathbb{R}\times \mathbb{R}_{+}$ whose dimension is equal to one.
\end{Proposition}

This proposition essentially states that the contours, or level sets, of $|W_Y|$ are locally $C^1$. The level sets of homogeneous Gaussian random fields, defined by a constant threshold $c$, are well studied in the context of excursion theory \cite{adler2007random}. In our setting, $|W_\Phi|$ is neither Gaussian
nor homogeneous. Hence, the existing results from excursion theory are not directly applicable.
The proof of Proposition \ref{thm:contour:c1} is provided in Section \ref{sec:proof:thm:contour:c1}.
A deeper analysis of the geometry and topology of $L_c$ would require extending
existing results to the inhomogeneous setting, which lies beyond the scope of this paper and will be pursued in future work.

\begin{Example}\label{example:wavelet}
The set of wavelets that satisfy the assumptions in Proposition \ref{thm:contour:c1}
includes some of the generalized Morse wavelets discussed in \cite[(14)]{lilly2010analytic}
and the analytic wavelets derived in \cite[(6)]{holighaus2019characterization}.
For instance, the generalized Morse wavelets, whose Fourier transform takes the form
$\widehat{\psi}(\lambda)= a_{\beta_{1},\beta_{2}}\lambda^{\beta_{1}}e^{-\lambda^{\beta_{2}}}1_{[0,\infty)}(\lambda),
\ \beta_{1}\geq1,\ \beta_{2}>0$,
where $a_{\beta_{1},\beta_{2}}$ is a normalizing constant.
Another example is the Klauder wavelet, whose Fourier transform is given by
$\widehat{\psi}(\lambda) = \lambda^{\alpha}e^{-\gamma \lambda}e^{i\beta \log \lambda} 1_{[0,\infty)}(\lambda)$,
where $\alpha\geq1$, $\beta\in \mathbb{R}$, and $\gamma\in \mathbb{C}$ with $\textup{Re}(\gamma)>0$.
\end{Example}

\subsection{Marginal distributions of AWT phase}\label{sec:phase}

\begin{figure*}[htb!]
\centering
\subfigure[][]
{\includegraphics[scale = 0.45]{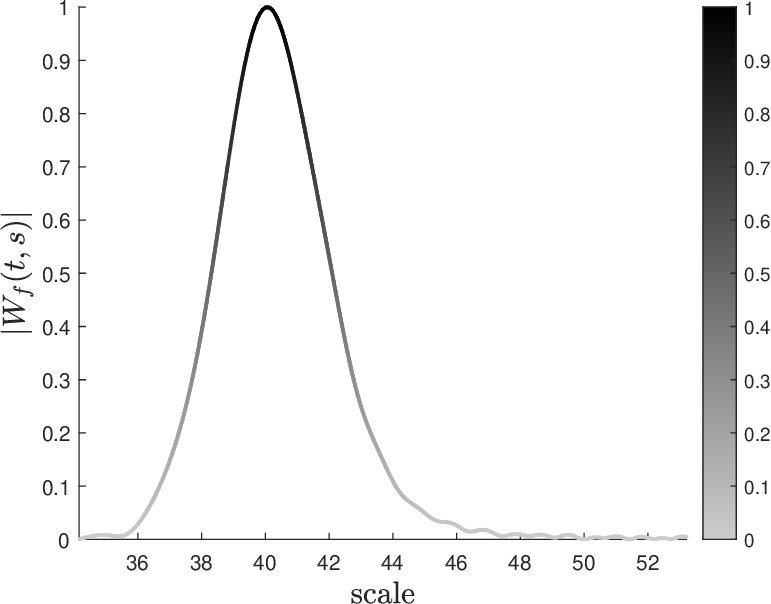}}
\subfigure[][]
{\includegraphics[scale = 0.45]{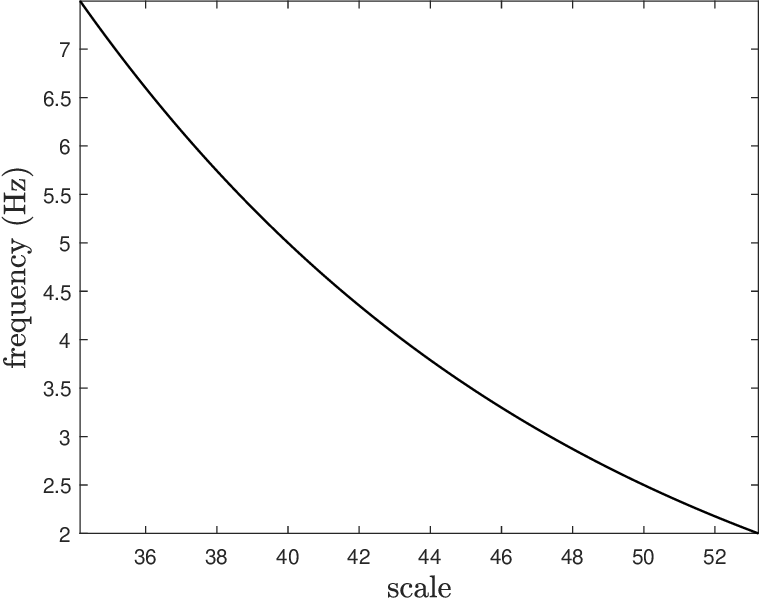}}\\
\subfigure[][]
{\includegraphics[scale = 0.42]{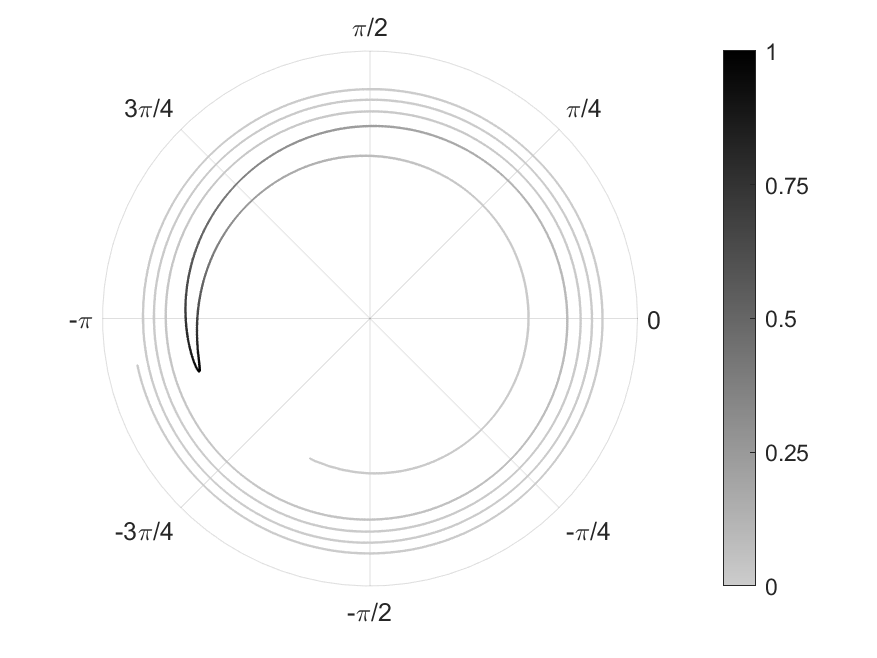}}
\subfigure[][]
{\includegraphics[scale = 0.42]{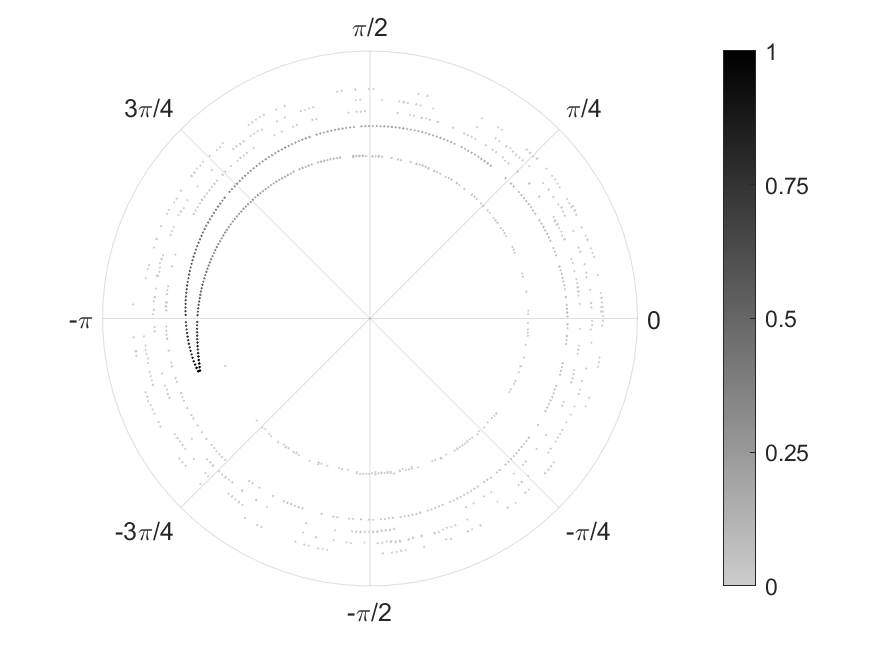}}
 \caption{(a) Magnitude of
 $W_{f}(t,s)$ as a function of the scale variable, where $f(t) = \cos(2\pi \phi(t))$ with $\phi(t)=\frac{1}{2}t^{2}$ (i.e., $\phi'(t)=t$). The signal $f$ is sampled at 200 Hz and the AWT of $f$ is evaluated at 5 seconds.
 Here, we use a mother wavelet $\psi$ with center frequency $\omega_{\psi}=80$ Hz.
 (b) Relationship between the scale variable and frequency. The peak in (a) occurs approximately  at scale $s=40$,
 corresponding to a frequency of $5$ Hz.
 (c) Polar plot of the phase of $W_{f}(t,s)$ as a function of the scale variable, i.e., a scatter plot of $(s\cos(\Theta_{Y}(5,s)),s\sin(\Theta_{Y}(5,s)))$ for various $s$.
  The color gradient indicates the variation in the magnitude $|W_{f}(5,s)|$, normalized to range from zero to one.
 (d) Polar plot of the peak locations of the histogram of the random variable $\Theta_{Y}(5,s)$.
 Here, the color is determined by the normalized magnitude $|W_{f}(5,s)|$.
 The histogram of $\Theta_{Y}(5,s)$ is generated using $10^4$ sample paths of the random process $\Phi$.}
 \label{fig:phase_polarplot_square}
\end{figure*}

In this section, we analyze the random field arising from the phase of the AWT of the random process $Y$.
Visualizations of this phase-related field can be found in the fourth and fifth rows of Figure~\ref{fig:contour}.
For any $n\in \mathbb{N}$, times $t_{1}$,$t_{2}$,...,$t_{n}\in \mathbb{R}$, and scales $s_{1}$,$s_{2}$,...,$s_{n}>0$,
we consider the $n$-dimensional random vector
\begin{align*}
\Theta_{Y}=\left[\Theta_{Y}(t_{1},s_{1})\ \Theta_{Y}(t_{2},s_{2})\ \cdots\ \Theta_{Y}(t_{n},s_{n})\right]^{\top}\in [0,2\pi)^{n}.
\end{align*}
\begin{Proposition}\label{prop:jpdf:Thetan}
Assume that Assumptions \ref{assump:analytic} and \ref{assump:Gaussian} hold.
Additionally,
suppose that the covariance matrix $\Gamma$ of $\mathbf{W}_{Y}$
is invertible.

(a) The probability density function of the random vector $\Theta_{Y}$
can be expressed as
\begin{align}\label{jpdf:ThetaY:n:nonnull}
&p\left(\theta_{1:n};\Theta_{Y}\right)
=\int_{0}^{\infty}\hspace{-0.15cm}\int_{0}^{\infty}\hspace{-0.2cm}\cdots \hspace{-0.15cm}\int_{0}^{\infty}\left(\overset{n}{\underset{\ell=1}{\prod}}r_{\ell}\right)
\hspace{-0.1cm} p\left(\mathbf{w};\mathbf{W}_{Y}\right)dr_{1}dr_{2}\cdots dr_{n}
\end{align}
for $\theta_{1:n}=(\theta_{1},\theta_{2},\ldots,\theta_{n})\in [0,2\pi)^{n}$,
where $\mathbf{w}=[r_{1}e^{i\theta_{1}} \  r_{2}e^{i\theta_{2}} \ \cdots \ r_{n}e^{i\theta_{n}}]^{\top}$.

(b) In the case $n=1$,
for any $t\in \mathbb{R}$ and $s>0$, the probability density function of $\Theta_{Y}(t,s)$ has the representation
\begin{align}\label{margianl_pdf_phase}
&p(\theta;\Theta_{Y}(t,s))
=\frac{1}{2\pi }e^{-q}
+
\frac{\sqrt{q}}{2\sqrt{\pi}}
\cos(\theta-\theta_{f}) \textup{exp}\left\{\hspace{-0.05cm}-q\left[1-\cos^{2}(\theta-\theta_{f})\right]\hspace{-0.05cm}\right\}
\hspace{-0.1cm}B_{q}(\theta-\theta_{f}),
\end{align}
for $\theta\in [0,2\pi)$, where
$\theta_{f}=\arg(W_{f}(t,s))$,
 $q=q(t,s)$ is a measure of signal quality defined in (\ref{def:SNR}),
\begin{equation}\label{def:B_q}
B_{q}(u)=1+\textup{erf}\left(\sqrt{q}\cos(u)\right),\ u\in \mathbb{R},
\end{equation}
and
\begin{equation*}
\textup{erf}(z) = \frac{2}{\sqrt{\pi}}\int_{0}^{z}e^{-v^{2}}dv,\ z\in \mathbb{R}.
\end{equation*}

\end{Proposition}
The proof of Proposition \ref{prop:jpdf:Thetan} is provided in Section \ref{sec:proof:prop:jpdf:Thetan}.
In the null case, the probability density function of the random vector $\Theta_{Y}$, given in (\ref{jpdf:ThetaY:n:nonnull}), simplifies to the form presented in Corollary \ref{corollary:phase:jpdf:null} below.
On the other hand, (\ref{margianl_pdf_phase}) shows that the probability density function $p(\theta;\Theta_{Y}(t,s))$ attains its maximum at $\arg(W_{f}(t,s))$.
This relationship is illustrated in Figure \ref{fig:phase_polarplot_square}.

\begin{Corollary}\label{corollary:phase:jpdf:null}
Assume that Assumptions \ref{assump:analytic} and \ref{assump:Gaussian} hold.
Additionally,
suppose that the covariance matrix $\Gamma$ of $\mathbf{W}_{Y}$
is invertible.

(a) In the null case, that is, $Y=\Phi$,
\begin{align*}
p\left(\theta_{1:n};\Theta_{\Phi}\right)
=
\frac{(n-1)!}{2\pi^{n}\textup{det}(\Gamma) }
\int_{S_{+}^{n-1}}
\left(\overset{n}{\underset{\ell=1}{\prod}}\omega_{\ell}\right)
\left(\omega^{\top}\mathbf{M}\omega\right)^{-n}d\omega
\end{align*}
for $\theta_{1:n}\in [0,2\pi)^n$, where
$S_{+}^{n-1}=
\{\omega=[\omega_{1},\ldots,\omega_{n}]^{\top}\in\mathbb{R}^{n}\mid \|\omega\|=1\ \textup{and}\ \omega_{\ell}\geq 0\ \textup{for}\ \ell=1,2,\ldots,n\}$,
$d\omega$ is the surface measure on $S_{+}^{n-1}$, and
\begin{align*}
\mathbf{M} = \textup{diag}\left(e^{-i\theta_1}, \dots, e^{-i\theta_n} \right)\Gamma^{-1} \textup{diag}\left(e^{i\theta_1}, \dots, e^{i\theta_n}\right).
\end{align*}

(b) In the null case with $n=1$, $\Theta_{\Phi}$ is uniformly distributed over $[0,2\pi)$. When $n=2$,
 $p\left(\theta_{1:2};\Theta_{\Phi}\right)$
depends only on the difference $\theta_{2}-\theta_{1}$ as follows
\begin{align}\label{jpdf:Theta2:null}
&p\left(\theta_{1:2};\Theta_{\Phi}\right)
=\frac{1 }{4\pi^{2}}\frac{\textup{det}(\Gamma)}{D(\theta_{2}-\theta_{1})}
+\frac{1 }{4\pi^{2}}\frac{\textup{det}(\Gamma)C(\theta_{2}-\theta_{1})}
{D(\theta_{2}-\theta_{1})^{3/2}}
\Big[\frac{\pi}{2}+
\tan^{-1}\Big(\frac{C(\theta_{2}-\theta_{1})}
{\sqrt{D(\theta_{2}-\theta_{1})}}
\Big)\Big],
\end{align}
where
\begin{align*}
C(\varphi)=\Re(\Gamma_{12}e^{i\varphi}),\ \varphi\in \mathbb{R},
\end{align*}
and
\begin{align*}
D(\varphi)= \Gamma_{11}\Gamma_{22}-|C(\varphi)|^2.
\end{align*}
\end{Corollary}
For any $t_{1},t_{2}\in \mathbb{R}$ and $s_{1},s_{2}>0$, the joint probability density function of $\Theta_{\Phi}(t_{1},s_{1})$ and $\Theta_{\Phi}(t_{2},s_{2})$ given by (\ref{jpdf:Theta2:null})
attains its maximum when
\begin{align}\notag
\theta_{2} - \theta_{1} =& \arg(\overline{\Gamma_{12}}) \\\label{observation_phase}=&\arg\left(\mathbb{E}\left[\overline{W_{\Phi}(t_{1},s_{1})}W_{\Phi}(t_{2},s_{2})\right]\right).
\end{align}
The relationship (\ref{observation_phase}) shows
that the most likely instantaneous phase difference $\theta_{2} - \theta_{1}$
between the wavelet coefficients
$W_{\Phi}(t_{1},s_{1})$ and $W_{\Phi}(t_{2},s_{2})$ aligns with the phase of the covariance
$\mathbb{E}[\overline{W_{\Phi}(t_{1},s_{1})}W_{\Phi}(t_{2},s_{2})]$.
The proof of  Corollary \ref{corollary:phase:jpdf:null} is provided in Section \ref{sec:proof:corollary:phase:jpdf:null}.

We conclude this section by presenting the convergence rate of the random variable $\Theta_{Y}(t,s)$ to $\arg(W_{f}(t,s)).$
For any $\varepsilon\in (0,\pi/2)$,
consider an arc on the quotient torus $\mathbb{R}/2\pi \mathbb{Z}$ as follows
\begin{align*}
A_{\varepsilon} = \left\{\theta\in [0,2\pi) \mid \theta\notin (\theta _{f}-\varepsilon,\theta_{f}+\varepsilon)+2n\pi\ \textup{for any}\ n\in \mathbb{Z}\right\}.
\end{align*}
Because $0\leq\cos(\theta-\theta_{f})\leq \cos(\varepsilon)$,
$1-\cos^{2}(\theta-\theta_{f})>1-\cos^{2}(\varepsilon)>0$,
and
$
B_{q}(\theta-\theta_{f})\leq 2
$
for any $\theta\in A_{\varepsilon}$ and $q\geq 0$,
Proposition \ref{prop:jpdf:Thetan} leads to the following corollary.

\begin{Corollary}\label{corollary:phase:expo}
Assume that all the conditions in Proposition \ref{prop:jpdf:Thetan} hold and adopt its notation. Then, for any $\varepsilon\in (0,\pi/2)$,
\begin{align}\notag
\mathbb{P}\left(d\left(\Theta_{Y}(t,s),\arg(W_{f}(t,s))\right)>\varepsilon\right)
\leq \frac{\varepsilon}{\pi }e^{-q}
+\frac{2\varepsilon \sqrt{q}}{\sqrt{\pi}}\cos(\varepsilon)\textup{exp}\left\{-q\left[1-\cos^{2}(\varepsilon)\right]\right\},
\end{align}
where $d(\cdot,\cdot)$ denotes the cosine distance, defined as
\begin{align*}
d\left(\Theta_{Y}(t,s),\arg(W_{f}(t,s))\right)
= 1-\left|\cos\big(\Theta_{Y}(t,s)-\arg(W_{f}(t,s))\big)\right|.
\end{align*}
\end{Corollary}

\begin{figure}[hbt!]
\centering
\subfigure[][]
{\includegraphics[scale = 0.27]{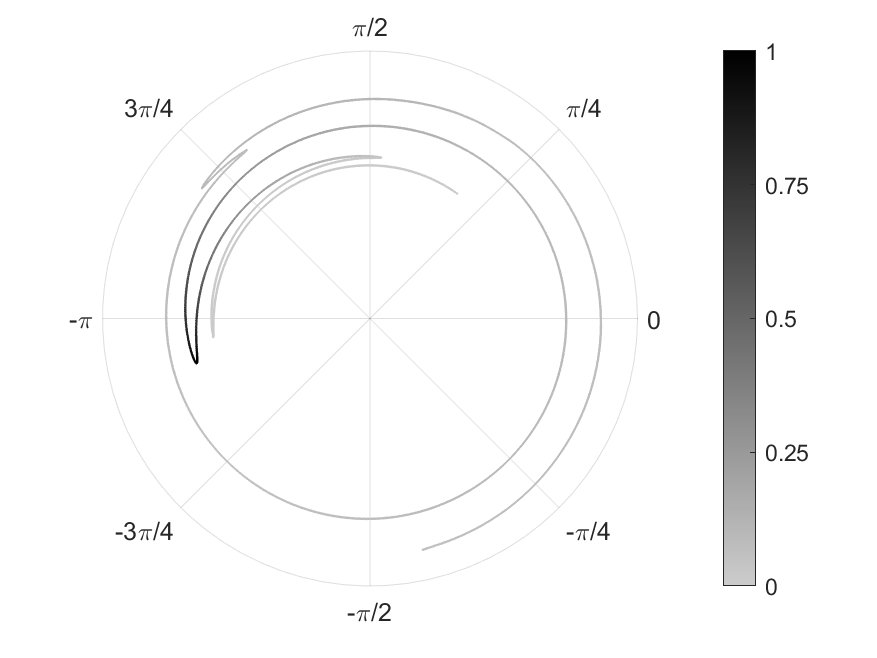}}
\subfigure[][]
{\includegraphics[scale = 0.27]{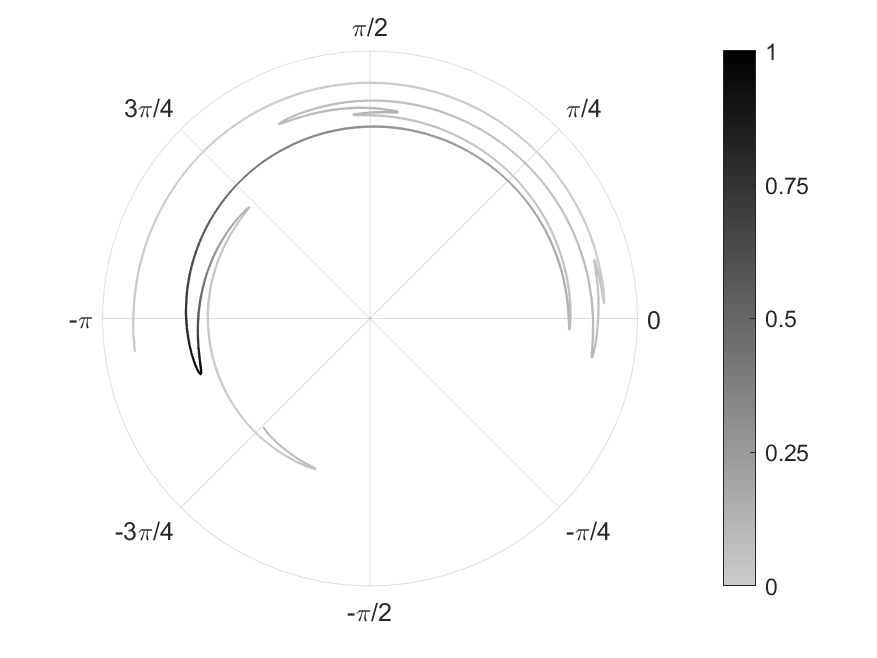}}
\subfigure[][]
{\includegraphics[scale = 0.27]{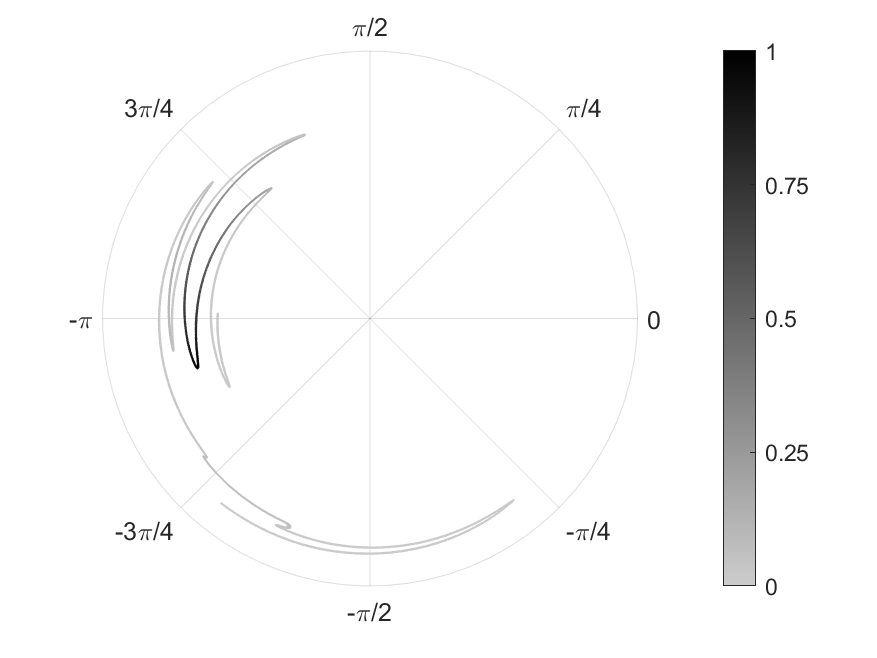}}\\
\subfigure[][]
{\includegraphics[scale = 0.27]{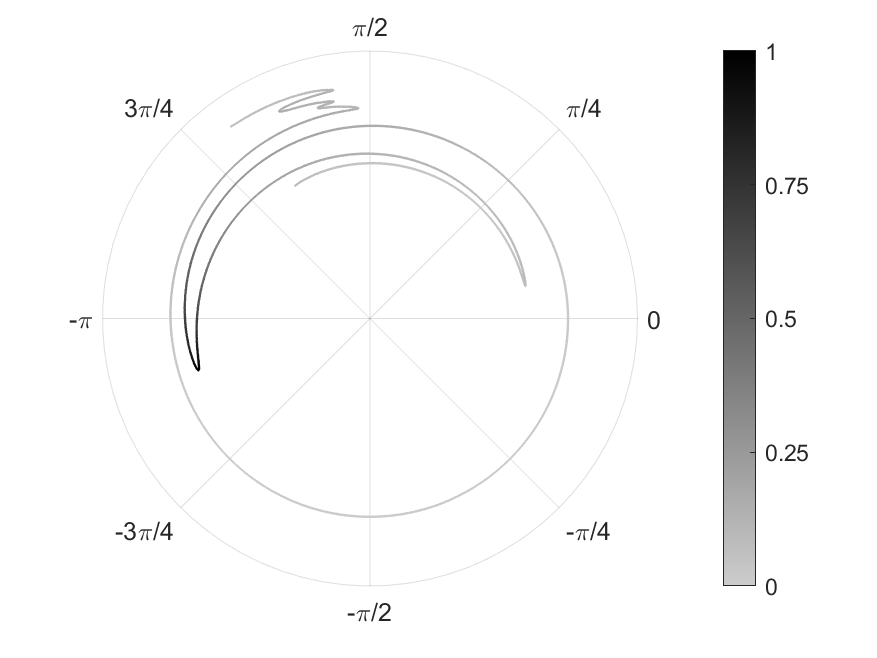}}
\subfigure[][]
{\includegraphics[scale = 0.27]{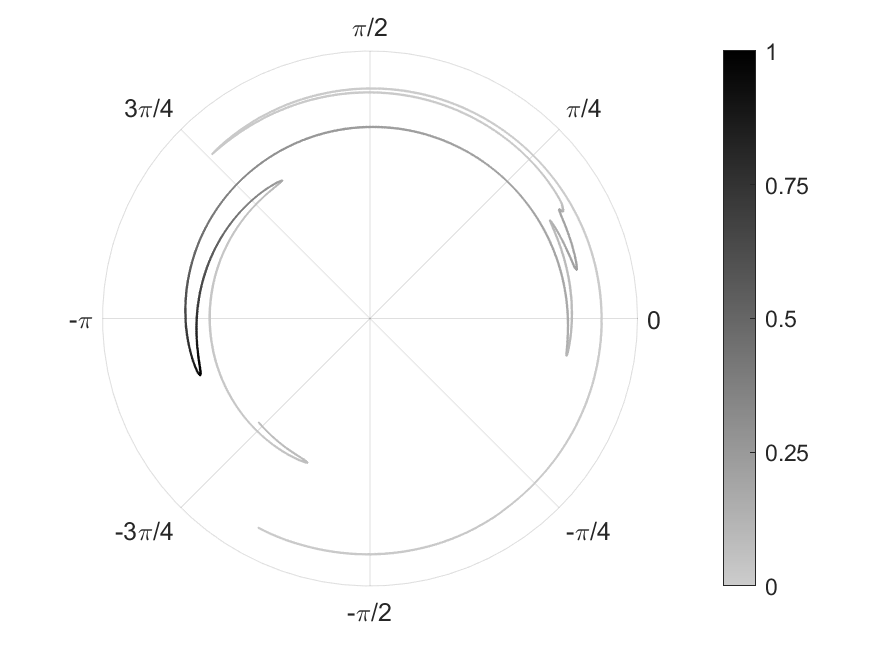}}
\subfigure[][]
{\includegraphics[scale = 0.27]{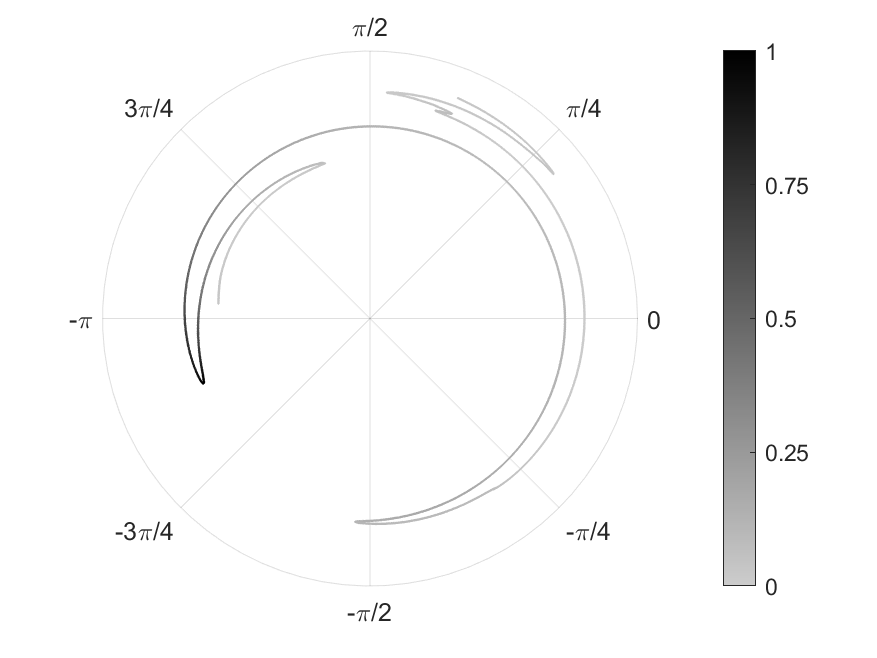}}
 \caption{Sample paths of $(s\cos(\Theta_{Y}(t,s)),s\sin(\Theta_{Y}(t,s)))$ for $t=5$, where $Y=f+\Phi$ and $f(t) = \cos(2\pi \phi(t))$ with $\phi(t)=\frac{1}{2}t^{2}$. The color gradient in each subfigure represents the variation in the corresponding sample path of $|W_{Y}(5,s)|$, normalized to range from zero to one.}
 \label{fig:sample_phase_polarplot_square}
\end{figure}

Corollary \ref{corollary:phase:expo} shows that the phase of the AWT of the noisy signal $Y$ converges to that of the clean signal $f$ in probability
when the signal-to-noise ratio $q(t,s)\rightarrow \infty$.
Furthermore,
the upper bound of $\mathbb{P}\left(d\left(\Theta_{Y}(t,s),\arg(W_{f}(t,s))\right)>\varepsilon\right)$ decays exponentially to zero as the noise strength decreases to zero (i.e., $q(t,s)\rightarrow\infty$).
Figure \ref{fig:sample_phase_polarplot_square} presents numerical simulations
showing that individual samples of $\Theta_{Y}(t,s)$ closely align with
$\theta_{f}(t,s)$ at scale $s$ where the corresponding samples of $|W_{Y}(t,s)|$ are large.

\subsection{Covariance of AWT's magnitudes and phases}\label{sec:cov}
In this section, we investigate the following quantities:
$$\textup{Cov}\left(|W_{Y}(t_{1},s_{1})|^2, |W_{Y}(t_{2},s_{2})|^2\right),$$
$$\textup{Cov}\left(|W_{Y}(t_{1},s_{1})|, |W_{Y}(t_{2},s_{2})|\right),$$
and
$$\mathbb{E}\left[e^{i(\Theta_{\Phi}(t_{1},s_{1})-\Theta_{\Phi}(t_{2},s_{2}))}\right],$$
where $t_{1},t_{2}\in \mathbb{R}$ and $s_{1},s_{2}>0$ are arbitrary.
Our goal is to characterize the relationship between these covariances and the covariance matrix $\Gamma$ of
the AWT of the noise $\Phi$, which is defined in (\ref{def:Gamma}) and assumed to be invertible.

\begin{Proposition}\label{cov:mag:generalpsi}
Assume that Assumption \ref{assump:Gaussian} holds.
For any $t_{1},t_{2}\in \mathbb{R}$ and $s_{1},s_{2}>0$,
\begin{align}\notag
&\textup{Cov}\left(|W_{Y}(t_{1},s_{1})|^2, |W_{Y}(t_{2},s_{2})|^2\right)
\\\notag=&4\underset{\Box,\triangle \in \{\Re,\Im\}}{\sum} \Box(W_{f}(t_{1},s_{1}))\triangle(W_{f}(t_{2},s_{2}))
\sqrt{s_{1}s_{2}}\int_{\mathbb{R}}e^{i(t_{1}-t_{2})\lambda} \overline{\widehat{\Box(\psi)}}(s_{1}\lambda)\widehat{\triangle(\psi)}(s_{2}\lambda)F(d\lambda)
\\\notag&+2s_{1}s_{2}\left[\int_{\mathbb{R}}e^{i(t_{1}-t_{2})\lambda} \overline{\widehat{\Box(\psi)}}(s_{1}u)
\widehat{\triangle(\psi)}(s_{2}u)F(du)\right]^{2},
\end{align}
where $F$ is the spectral measure of $C_{\Phi}$.
\end{Proposition}
We note that  Proposition \ref{cov:mag:generalpsi} does not require the wavelet $\psi$ to be analytic.
Its proof, provided in Section \ref{sec:proof:cov:mag:generalpsi}, is based on It$\hat{\textup{o}}$'s formula.
In some applications \cite{anden2014deep,mallat2012group}, the magnitude of $W_{Y}$ (without squaring) is preferred for feature extraction in signal segment classification tasks. This preference arises because the absolute value operator is non-expansive, a desirable property in both signal processing and machine learning.
However, deriving
$\textup{Cov}\left(|W_{Y}(t_{1},s_{1})|, |W_{Y}(t_{2},s_{2})|\right)$ is more challenging since computing the absolute value of a complex number involves the square root. The Malliavin calculus method used in the proof of Proposition \ref{cov:mag:generalpsi} is not applicable in this case.
The computation of the covariance of the componentwise magnitudes of $\mathbf{W}_{Y}$ relies on the probability density function (\ref{pdf_WY}) obtained from Lemma \ref{prop:circular}.

\begin{Proposition}\label{prop:Cov_general}
Assume that Assumptions \ref{assump:analytic} and \ref{assump:Gaussian} hold.
For any $t_{1},t_{2}\in \mathbb{R}$ and $s_{1},s_{2}>0$,
the covariance of the magnitudes $|W_{Y}(t_{1},s_{1})|$ and $|W_{Y}(t_{2},s_{2})|$ has the representation
\begin{align}\label{format:cov:mag}
&\textup{Cov}\left(|W_{Y}(t_{1},s_{1})|,|W_{Y}(t_{2},s_{2})|\right)
\\\notag=&\int_{\mathbb{C}}|w|p\left(w;W_{Y}(t_{1},s_{1})W_{Y}(t_{2},s_{2})\right)dw
-\overset{2}{\underset{\ell=1}{\prod}}\int_{\mathbb{C}}|w|p\left(w;W_{Y}(t_{\ell},s_{\ell}) \right)dw,
\end{align}
where
\begin{align}\notag
&p\left(w;W_{Y}(t_{1},s_{1})W_{Y}(t_{2},s_{2})\right)
\\\notag=&
\int_{0}^{\infty}\int_{0}^{2\pi} \frac{1}{\pi^{2}\textup{det}(\Gamma)r}
\textup{exp}\left(-\left(
\left[\begin{array}{cc}re^{i\theta}\\ wr^{-1}e^{-i\theta}\end{array}\right]-\left[\begin{array}{cc}W_{f}(t_{1},s_{1})\\ W_{f}(t_{2},s_{2})\end{array}\right]\right)^{*}\right.
\\\label{prop:pdf:product}&\left.\Gamma^{-1}
\left(\left[\begin{array}{cc}re^{i\theta}\\ wr^{-1}e^{-i\theta}\end{array}\right]-\left[\begin{array}{cc}W_{f}(t_{1},s_{1})\\ W_{f}(t_{2},s_{2})\end{array}\right]\right)\right)d\theta dr
\end{align}
and
\begin{align*}
p\left(w;W_{Y}(t_{\ell},s_{\ell}) \right) = \frac{1}{\pi\Gamma_{\ell\ell}}
\textup{exp}\left(-\frac{1}{\Gamma_{\ell\ell}}\left|w-W_{f}(t_{\ell},s_{\ell})\right|^{2}\right)
\end{align*}
for $\ell=1,2.$
\end{Proposition}
The proof of Proposition \ref{prop:Cov_general} is provided in Section \ref{sec:proof:prop:Cov_general}.
Denote the hypergeometric functions by
\begin{align*}
{}_2F_1\left(a,b;c;x\right)=
\frac{\Gamma(c)}{\Gamma(a)\Gamma(b)}\overset{\infty}{\underset{k=0}{\sum}}
\frac{\Gamma\left(k+a\right)\Gamma\left(k+b\right)}
{\Gamma(k+c)}\frac{x^{k}}{k!},
\end{align*}
where $x\in[0,1)$ and $a,b,c>0$.

\begin{figure}[hbt!]
\centering
\subfigure[][${}_2F_1\left(\frac{3}{2},\frac{3}{2};1;x\right)$]
{\includegraphics[scale = 0.5]{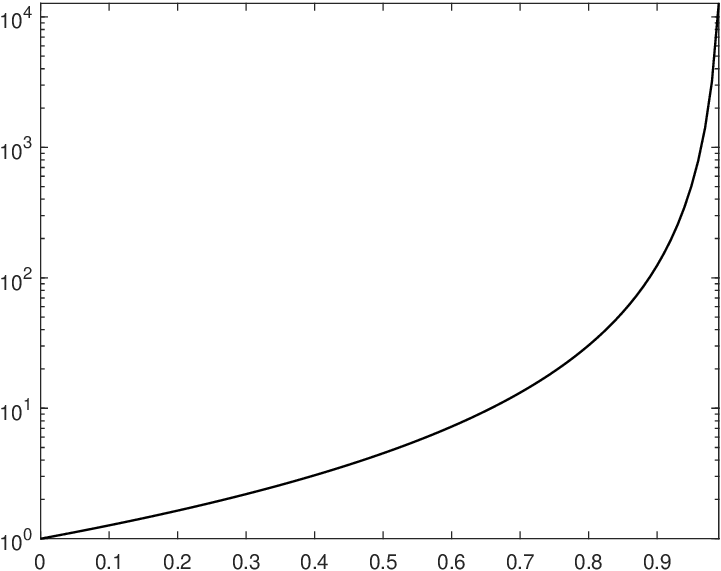}}
\subfigure[][$(1-x)^{2}{}_2F_1\left(\frac{3}{2},\frac{3}{2};1;x\right)$]
{\includegraphics[scale = 0.5]{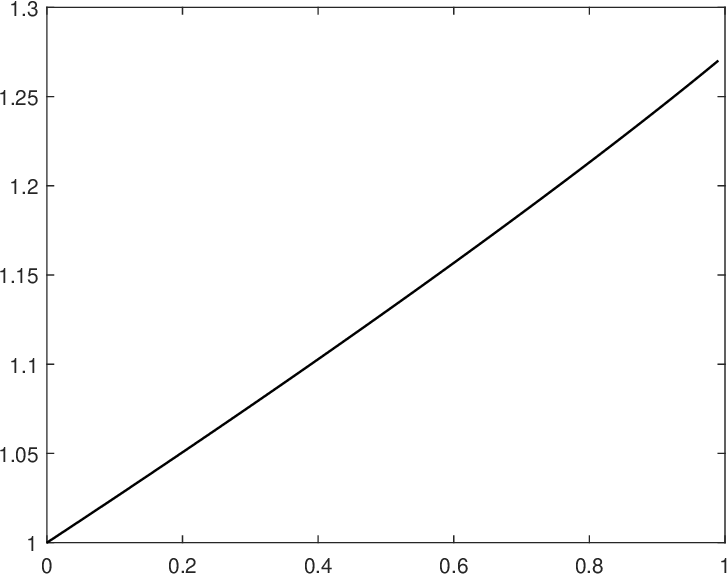}}
 \caption{Graphs of the  hypergeometric-related functions in (\ref{eq:correlation_mag})}
 \label{fig:F1}
\end{figure}
\begin{Corollary}\label{prop:Cov_null}
Assume that Assumptions \ref{assump:analytic} and \ref{assump:Gaussian} hold.
For the null case, i.e.,  $Y=\Phi$,
\begin{align}\notag
&\textup{Cov}\left(\left|W_{\Phi}(t_{1},s_{1})\right|,\left|W_{\Phi}(t_{2},s_{2})\right|\right)
\\\notag=&\frac{\pi}{4}\sqrt{\Gamma_{11}\Gamma_{22}}
\left[\hspace{-0.1cm}\left(\hspace{-0.06cm}1-\frac{|\Gamma_{12}|^{2}}{\Gamma_{11}\Gamma_{22}}\hspace{-0.06cm}\right)^2
{}_2F_1\left(\hspace{-0.06cm}\frac{3}{2},\frac{3}{2};1;\frac{|\Gamma_{12}|^{2}}{\Gamma_{11}\Gamma_{22}}\hspace{-0.06cm}\right)-1\hspace{-0.05cm}\right],
\end{align}
where $\Gamma_{11},$ $\Gamma_{22}$, and $\Gamma_{12}$ are given in (\ref{entry_Gamma}).
The correlation coefficient between
$|W_{\Phi}(t_{1},s_{1})|$ and $|W_{\Phi}(t_{2},s_{2})|$ is given by
\begin{align}\label{eq:correlation_mag}
&\textup{Corr}\left(\left|W_{\Phi}(t_{1},s_{1})\right|,\left|W_{\Phi}(t_{2},s_{2})\right|\right)
\\\notag=&\left(\hspace{-0.06cm}\frac{4}{\pi}-1\hspace{-0.06cm}\right)^{-1}
\hspace{-0.06cm}\left[\hspace{-0.06cm}\left(\hspace{-0.06cm}1-\frac{|\Gamma_{12}|^{2}}{\Gamma_{11}\Gamma_{22}}\hspace{-0.06cm}\right)^2
\hspace{-0.06cm}{}_2F_1\left(\frac{3}{2},\frac{3}{2};1;\frac{|\Gamma_{12}|^{2}}{\Gamma_{11}\Gamma_{22}}\right)-1\hspace{-0.06cm}\right].
\end{align}
When $|\Gamma_{12}|^{2}/(\Gamma_{11}\Gamma_{22})\rightarrow 0,$
\begin{align}\label{implication:corollary:mag}
\textup{Corr}\left(\left|W_{\Phi}(t_{1},s_{1})\right|,\left|W_{\Phi}(t_{2},s_{2})\right|\right)
\sim \left(\frac{16}{\pi}-4\right)^{-1}\frac{|\Gamma_{12}|^{2}}{\Gamma_{11}\Gamma_{22}}.
\end{align}
\end{Corollary}
Here, for two functions $g_{1}$ and $g_{2}$, the notation $g_{1}(\varepsilon)\sim g_{2}(\varepsilon)$ as $\varepsilon\rightarrow 0$ means that $\underset{\varepsilon\rightarrow 0}{\lim}g_{1}(\varepsilon)/g_{2}(\varepsilon)= 1$. The observation (\ref{implication:corollary:mag}) provides the relationship between the correlation decay of
the magnitudes $|W_{\Phi}(t_{1},s_{1})|$ and $|W_{\Phi}(t_{2},s_{2})|$
and that of $W_{\Phi}(t_{1},s_{1})$ and $W_{\Phi}(t_{2},s_{2})$.
The proof of Corollary \ref{prop:Cov_null} is provided in Section \ref{sec:proof:prop:Cov_null}.

\begin{Remark}
If $\Phi$ is a Gaussian white noise, meaning that $F(d\lambda)=d\lambda$, and
$\psi$ is a Morse wavelet with
$$
\widehat{\psi}(\lambda) = \lambda^{\alpha}e^{-\lambda}1_{[0,\infty)}(\lambda),
$$
then
$ \Gamma_{11}=\Gamma_{22}= \Gamma(2\alpha+1)2^{-(2\alpha+1)}$ and
\begin{align}\notag
\Gamma_{12} =& \left(s_{1}s_{2}\right)^{\alpha+\frac{1}{2}}\int_{0}^{\infty}e^{i(t_{1}-t_{2})\lambda}\lambda^{2\alpha}e^{-(s_{1}+s_{2})\lambda}
d\lambda
\\\notag=&\Gamma(2\alpha+1)
\frac{\left(s_{1}s_{2}\right)^{\alpha+\frac{1}{2}}}{\left[(s_{1}+s_{2})-i(t_{1}-t_{2})\right]^{2\alpha+1}},
\end{align}
which implies that
\begin{align*}
\frac{\Gamma_{12}}{\sqrt{\Gamma_{11}\Gamma_{22}}} = 2^{2\alpha+1}
\frac{\left(s_{1}s_{2}\right)^{\alpha+\frac{1}{2}}}{\left[(s_{1}+s_{2})-i(t_{1}-t_{2})\right]^{2\alpha+1}}.
\end{align*}
\begin{itemize}
\item When $t_{1}=t_{2}$, $\log_{2}(s_{1})=j_{1}$,  $\log_{2}(s_{2})=j_{2}$, and $\Delta = |j_{1}-j_{2}|$,
\begin{align}\notag
\frac{|\Gamma_{12}|}{\sqrt{\Gamma_{11}\Gamma_{22}}}=&2^{2\alpha+1}
\frac{\left(s_{1}s_{2}\right)^{\alpha+\frac{1}{2}}}{(s_{1}+s_{2})^{2\alpha+1}}
=\left(\frac{2^{\Delta/2}+2^{-\Delta/2}}{2}\right)^{-(2\alpha+1)}.
\end{align}
This observation implies that the correlation between $\left|W_{\Phi}(t_{1},s_{1})\right|$ and $\left|W_{\Phi}(t_{2},s_{2})\right|$ exhibits polynomial decay with respect to the affine scale difference $|s_{1}-s_{2}|$,
or equivalently, exponential decay with respect to the logarithmic scale difference $|j_{1}-j_{2}|$.

\item When $s_{1}=s_{2}=s$,
\begin{align}\notag
\frac{|\Gamma_{12}|}{\sqrt{\Gamma_{11}\Gamma_{22}}}=&2^{2\alpha+1}
\frac{s^{2\alpha+1}}{\left|2s-i(t_{1}-t_{2})\right|^{2\alpha+1}}
=
\left[1+\left(\frac{t_{1}-t_{2}}{2s}\right)^{2}\right]^{-(\alpha+1/2)}.
\end{align}
This result indicates that
the correlation between $\left|W_{\Phi}(t_{1},s_{1})\right|$ and $\left|W_{\Phi}(t_{2},s_{2})\right|$ exhibits exponential decay with respect to the time difference $|t_{1}-t_{2}|$, with a faster decay rate at smaller scales.
\end{itemize}
\end{Remark}

In addition to analyzing the covariance of the magnitude of $W_{Y}$, we also investigate the covariance of its phase $\Theta_{Y}$,
as defined in (\ref{def:Theta_Y}).

\begin{figure}[hbt!]
  \centering
  \includegraphics[scale=0.7]{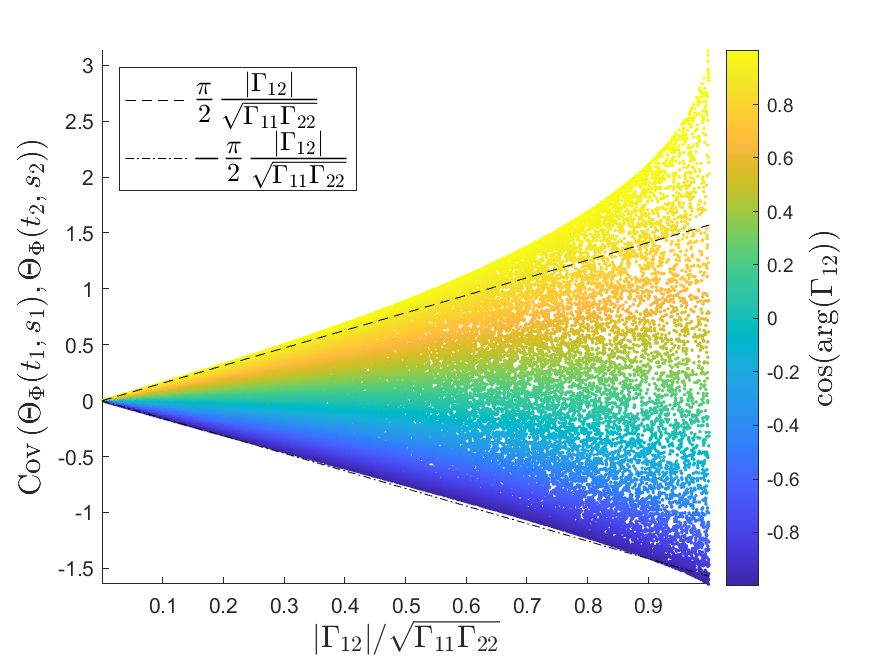}
  \caption{Effect of $|\Gamma_{12}|/\sqrt{\Gamma_{11}\Gamma_{22}}$ and $\cos(\arg(\Gamma_{12}))$
  on $\textup{Cov}\left(\Theta_{\Phi}(t_{1},s_{1}),\Theta_{\Phi}(t_{2},s_{2})\right)$.
  }
\label{fig:cov_phase}
\end{figure}

\begin{Corollary}\label{prop:Cov2_phase_null}
Assume that Assumptions \ref{assump:analytic} and \ref{assump:Gaussian} hold.
For any $t_{1},t_{2}\in \mathbb{R}$ and $s_{1},s_{2}>0$, the covariance of the phases $\Theta_{\Phi}(t_{1},s_{1})$ and $\Theta_{\Phi}(t_{2},s_{2})$ is a function of $|\Gamma_{12}|/\sqrt{\Gamma_{11}\Gamma_{22}}$ and $\arg(\Gamma_{12})$. Moreover, when $|\Gamma_{12}|/\sqrt{\Gamma_{11}\Gamma_{22}}$ tends to zero,
\begin{align}\notag
\textup{Cov}\left(\Theta_{\Phi}(t_{1},s_{1}),\Theta_{\Phi}(t_{2},s_{2})\right)\sim \frac{\pi}{2}\cos\big(\arg(\Gamma_{12})\big) \frac{|\Gamma_{12}|}{\sqrt{\Gamma_{11}\Gamma_{22}}},
\end{align}
where the notation $\sim$ means that
\begin{align}\notag
\left(\frac{\pi}{2}\cos\big(\arg(\Gamma_{12})\big) \frac{|\Gamma_{12}|}{\sqrt{\Gamma_{11}\Gamma_{22}}}\right)^{-1}\textup{Cov}\left(\Theta_{\Phi}(t_{1},s_{1}),\Theta_{\Phi}(t_{2},s_{2})\right)
\end{align}
converges to one when $|\Gamma_{12}|/\sqrt{\Gamma_{11}\Gamma_{22}}\rightarrow 0.$
\end{Corollary}
Based on the joint probability density function of the phases $\Theta_{\Phi}(t_{1},s_{1})$ and $\Theta_{\Phi}(t_{2},s_{2})$ given in Corollary \ref{corollary:phase:jpdf:null},
we plot $\textup{Cov}\left(\Theta_{\Phi}(t_{1},s_{1}),\Theta_{\Phi}(t_{2},s_{2})\right)$
in Figure \ref{fig:cov_phase} for various randomly generated $2\times 2$ positive definite matrices $\Gamma$.
This figure shows that when $\cos(\arg(\Gamma_{12}))$ is positive, the phases $\Theta_{\Phi}(t_{1},s_{1})$ and $\Theta_{\Phi}(t_{2},s_{2})$ tend to be positively correlated; conversely, when $\cos(\arg(\Gamma_{12}))$ is negative, the phases tend to be negatively correlated.
This observation aligns with the result stated in Corollary~\ref{prop:Cov2_phase_null}.
The proof of Corollary \ref{prop:Cov2_phase_null} is provided in Section \ref{sec:proof:prop:Cov2_phase_null}.

Because $\Theta_{\Phi}$ is a random field taking values in the quotient space $\mathbb{R}/2\pi \mathbb{Z}$,
we also consider the circular covariance between $\Theta_{\Phi}(t_{1},s_{1})$ and $\Theta_{\Phi}(t_{2},s_{2})$,
defined as
\begin{align}\notag
\mathbb{E}\left[e^{i(\Theta_{\Phi}(t_{1},s_{1})-\Theta_{\Phi}(t_{2},s_{2}))}\right]
= \int_{\mathbb{C}}\frac{w}{|w|}p\left(w;W_{\Phi}(t_{1},s_{1})\overline{W_{\Phi}(t_{2},s_{2})}\right)dw.
\end{align}

\begin{Corollary}\label{prop:Cov_phase_null}
Assume that Assumptions \ref{assump:analytic} and \ref{assump:Gaussian} hold.
For any $t_{1},t_{2}\in \mathbb{R}$ and $s_{1},s_{2}>0$, we have
\begin{align}\notag
\mathbb{E}\left[e^{i(\Theta_{\Phi}(t_{1},s_{1})-\Theta_{\Phi}(t_{2},s_{2}))}\right]
=&\frac{\pi}{4}e^{i\arg\left(\Gamma_{12}\right)}
\left(1-\frac{|\Gamma_{12}|^{2}}{\Gamma_{11}\Gamma_{22}}\right)
\\\label{phase_cov_null_result}&\times\frac{|\Gamma_{12}|}{\sqrt{\Gamma_{11}\Gamma_{22}}}
{}_2F_1\left(\frac{3}{2},\frac{3}{2};2;\frac{|\Gamma_{12}|^{2}}{\Gamma_{11}\Gamma_{22}}\right)
\end{align}
and $\mathbb{E}\left[e^{i(\Theta_{\Phi}(t_{1},s_{1})+\Theta_{\Phi}(t_{2},s_{2}))}\right]=0$.
\end{Corollary}

\begin{figure}[hbt!]
\centering
\subfigure[][${}_2F_1\left(\frac{3}{2},\frac{3}{2};2;x\right)$]
{\includegraphics[scale = 0.5]{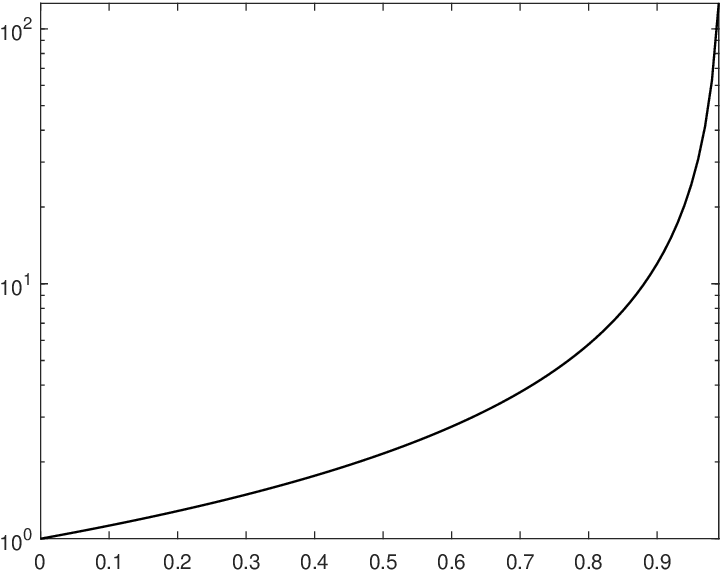}}
\subfigure[][$(1-x){}_2F_1\left(\frac{3}{2},\frac{3}{2};2;x\right)$]
{\includegraphics[scale = 0.5]{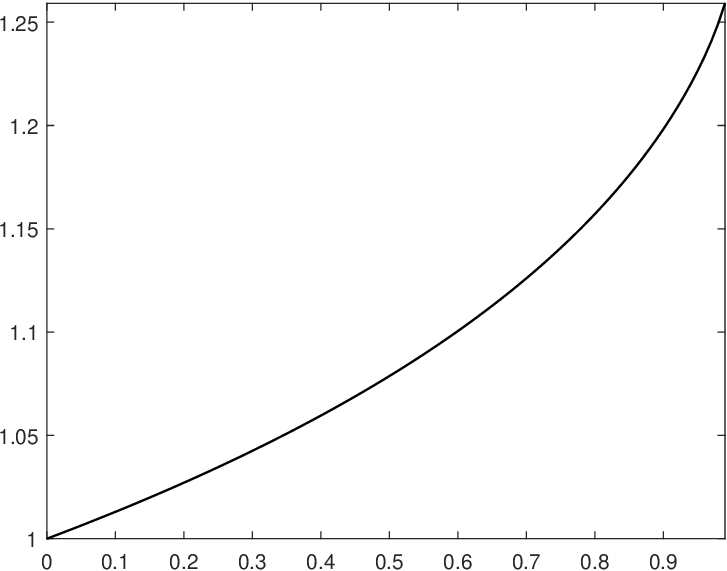}}
 \caption{Graphs of the  hypergeometric-related functions in (\ref{phase_cov_null_result})}
\end{figure}
The proof of Corollary \ref{prop:Cov_phase_null} is provided in Section \ref{sec:proof:prop:Cov_phase_null}.
Corollary \ref{prop:Cov_phase_null} shows that the circular covariance of the phases depends on the phase of $\Gamma_{12}$ and the magnitude of the correlation coefficient of the complex Gaussian random vector $\mathbf{W}_{\Phi}$.
Because $\underset{x\rightarrow 0+}{\lim}{}_2F_1\left(\frac{3}{2},\frac{3}{2};2;x\right)=1$, when $|\Gamma_{12}|^{2}/(\Gamma_{11}\Gamma_{22})$ is sufficiently small, we have the approximation
\begin{align}\notag
&\mathbb{E}\left[e^{i(\Theta_{\Phi}(t_{1},s_{1})-\Theta_{\Phi}(t_{2},s_{2}))}\right]
\sim
\frac{\pi}{4}e^{i\arg(\Gamma_{12})}
\frac{|\Gamma_{12}|}{\sqrt{\Gamma_{11}\Gamma_{22}}}.
\end{align}
On the other hand, by the fact $\underset{x\rightarrow 1-}{\lim} (1-x){}_2F_1\left(\frac{3}{2},\frac{3}{2};2;x\right)=\frac{4}{\pi}$,
when $|\Gamma_{12}|^2/(\Gamma_{11}\Gamma_{22})\rightarrow 1^{-}$, we have
\begin{align}\notag
&\mathbb{E}\left[e^{i(\Theta_{\Phi}(t_{1},s_{1})-\Theta_{\Phi}(t_{2},s_{2}))}\right]
\rightarrow
e^{i\arg(\Gamma_{12})}.
\end{align}

Finally, while it is well known that for a circularly symmetric Gaussian random vector the magnitude and phase of a single entry are independent, the dependence between the magnitude of one entry and the phase of another is less well understood. We address this gap in the following corollary, with the proof provided in Section~\ref{sec:proof:prop:Cov_phaseAM_null}.

\begin{Corollary}\label{prop:Cov_phaseAM_null}
Assume that Assumptions \ref{assump:analytic} and \ref{assump:Gaussian} hold.
For any $t_{1},t_{2}\in \mathbb{R}$ and $s_{1},s_{2}>0$ with $(t_{1},s_{1})\neq(t_{2},s_{2})$,
$|W_{\Phi}(t_{1},s_{1})|$ and $\Theta_{\Phi}(t_{2},s_{2})$ are independent.
\end{Corollary}

 For any $m,n\in \mathbb{N}$, let $t_{1},t_{2},...,t_{m+n}\in \mathbb{R}$ and $s_{1},s_{2},...,s_{m+n}\in \mathbb{R}_{+}$ with $(t_{i},s_{i})\neq(t_{j},s_{j})$ for $i\neq j$. Denote
\begin{align*}
\mathbf{A} = \left[\arg(W_{\Phi}(t_{1},s_{1}))\ \cdots\ \arg(W_{\Phi}(t_{m},s_{m}))\right]\in [0,2\pi)^{m}
\end{align*}
and
\begin{align*}
\mathbf{B} = \left[|W_{\Phi}(t_{m+1},s_{m+1})|\ \cdots\ |W_{\Phi}(t_{m+n},s_{m+n})|\right]\in \mathbb{R}^{n}_{\geq 0}.
\end{align*}
Corollary \ref{prop:Cov_phaseAM_null} shows that the elements of $\mathbf{A}$ and $\mathbf{B}$ are pairwise independent. It also follows that
the random vectors $\mathbf{A}$ and $\mathbf{B}$ are uncorrelated; that is, their $m\times n$ covariance matrix is the
zero matrix. However,  this does not imply that $\mathbf{A}$ and $\mathbf{B}$ are jointly independent because they are not Gaussian distributed.
This finding can be regarded as an initial step toward understanding the quantifications of phase-amplitude coupling and phase-phase coupling \cite{tort2010measuring,hulsemann2019quantification}. Choose $m=n$, times $t_j=t_{m+j}=j\Delta$ for all $j=1,\ldots,m$, where $\Delta$ is a discretization period, and scales so that $s_1=\ldots=s_m$ much larger than $s_{m+1}=\ldots=s_{2m}$. Under this configuration, $\mathbf{A}$ and $\mathbf{B}$ can be interpreted as the phase and amplitude, respectively, of bandpass-filtered signals. The uncorrelation between $\mathbf{A}$ and $\mathbf{B}$ when the input is pure noise highlights the need for caution when interpreting indices designed to quantify phase-amplitude and phase-phase coupling \cite{tort2010measuring,hulsemann2019quantification}, particularly when higher order moments are involved. A more thorough investigation of this issue will be pursued in future work.

\subsection{Gaussian approximation for the discretized AWT of non-Gaussian noise}\label{sec:Gaussian_approxi}

In practical applications, signals are usually observed in discrete time through a sampling process, and the noise is often non-Gaussian, nonstationary, and causal. The AWT is usually implemented via a Riemann-sum approximation,
which we refer to as the {\em discretized AWT}.
It is natural to ask whether the time-frequency representation generated by the discretized AWT is Gaussian, or asymptotically Gaussian in some sense, so that the developed random field results can be applied.
Clearly, if the input random process is Gaussian, then its time-frequency representation determined by
the discretized AWT is Gaussian. Therefore, it is less interesting to consider the Gaussian process discussed above and its discretization. The situation becomes interesting when the input random process is non-Gaussian. Since the AWT is defined through a kernel integration, it is intuitive to guess that the resulting time-frequency representation is approximately Gaussian, even if the noise is non-Gaussian but with sufficient moments control, and nonstationary but with sufficiently short-range dependence.
In this section, we show that this intuition is correct under mild assumptions.

We now introduce a discretized version of the AWT together with a control of the truncation error arising from approximating the integral by a finite Riemann sum.
Assume that a random process $\{\epsilon(t)\mid t\in \mathbb{R}\}$ is observed only at the discrete time points $\{t_{j}=j/\sqrt{n}\mid j\in \mathbb{Z}\}$, where $\sqrt{n}$ represents the sampling rate on the time axis.
If the sample paths of $\{\epsilon(t)\mid t\in \mathbb{R}\}$ are continuous, the AWT of $\epsilon$ can be approximated by the following Riemann sum:
\begin{align}\notag
W_{\epsilon}(t_{i},s) \approx&
\frac{1}{\sqrt{n}}\underset{j\in \mathbb{Z}}{\sum}\epsilon(t_{j})\frac{1}{\sqrt{s}}\overline{\psi\left(\frac{t_{j}-t_{i}}{s}\right)}
\\\notag=&\frac{1}{\sqrt{n}}\underset{j\in \mathbb{Z}}{\sum}\epsilon_{j}\frac{1}{\sqrt{s}}\overline{\psi\left(\frac{j-i}{\sqrt{n}s}\right)}
\\\label{def:discreteCWT1}=&
\frac{1}{\sqrt{n}}\underset{j\in \mathbb{Z}}{\sum}\epsilon_{i+j}\frac{1}{\sqrt{s}}\overline{\psi\left(\frac{j}{\sqrt{n}s}\right)}
\end{align}
for $i\in \mathbb{Z}$ and $s>0$, where $\epsilon_{j}=\epsilon(t_{j})$.
Let the scales of interest be $0<s_{1}<s_{2}<\cdots<s_{d}$, where $d\in \mathbb{N}$.
Let $[-A,A]$ be an interval outside of which $|\psi|$ is negligibly small, where $A>0$.
The right-hand side of (\ref{def:discreteCWT1}) can be rewritten as follows
\begin{align}\label{def:discreteCWT1s}
W_{\epsilon}(t_{i},s_{\ell}) \approx&
\frac{1}{\sqrt{n}}\underset{j\in D^{(\ell)}}{\sum}\epsilon_{i+j}\frac{1}{\sqrt{s_{\ell}}}\overline{\psi\left(\frac{j}{\sqrt{n}s_{\ell}}\right)}
+R_{A},
\end{align}
where
\begin{align*}
D^{(\ell)} = \left[-\lceil \sqrt{n}As_{\ell}\rceil,\lceil \sqrt{n}As_{\ell}\rceil\right]\cap \mathbb{Z},
\end{align*}
$\lceil \sqrt{n}As_{\ell} \rceil = \min\left\{k\in \mathbb{Z} \mid k\geq \sqrt{n}As_{\ell}\right\}$, and
\begin{align*}
R_{A}=\frac{1}{\sqrt{n}}\underset{j\in  \mathbb{Z} \setminus D^{(\ell)}}{\sum}\epsilon_{i+j}\frac{1}{\sqrt{s_{\ell}}}\overline{\psi\left(\frac{j}{\sqrt{n}s_{\ell}}\right)}.
\end{align*}
Note that
\begin{align}\notag
\mathbb{E}[|R_{A}|]\leq&
\sqrt{s_{\ell}}\left(\underset{k\in \mathbb{Z}}{\max}\ \mathbb{E}\left[|\epsilon_{k}|\right]\right)
\left(\frac{1}{\sqrt{n} s_{\ell}}\underset{\mathbb{Z} \setminus D^{(\ell)}}{\sum}\left|\psi\left(\frac{j}{\sqrt{n}s_{\ell}}\right)\right|\right).
\end{align}
Because $|\psi|\in L^{1}$,
\begin{align*}
\frac{1}{\sqrt{n}s_{\ell}}\underset{j\in \mathbb{Z} \setminus D^{(\ell)}}{\sum}\left|\psi\left(\frac{j}{\sqrt{n}s_{\ell}}\right)\right|
\rightarrow  \int_{\mathbb{R}\setminus[-A,A]}\left|\psi\left(x\right)\right| dx
\end{align*}
when the sampling rate $\sqrt{n}$ tends to infinity.
It follows that $\mathbb{E}[|R_{A}|]$ converges to zero as $A\rightarrow\infty$ if
$\max_{k\in \mathbb{Z}}\ \mathbb{E}\left[|\epsilon_{k}|\right]<\infty$.
Motivated by (\ref{def:discreteCWT1s}) and the above estimate of the remainder term $R_{A}$, we consider a discretized AWT defined by
\begin{align}\label{def:discreteCWT2}
\widetilde{W}_{\epsilon}(t_{i},s_{\ell}) =\frac{1}{n^{1/4}}\underset{j\in D^{(\ell)}}{\sum}\epsilon_{i+j}\overline{\psi_{j}^{(\ell)}},
\end{align}
where
\begin{align}\label{def:discrete_psi}
\psi_{j}^{(\ell)}
=\frac{1}{n^{1/4}s_{\ell}^{1/2}}\psi\left(\frac{j}{\sqrt{n}s_{\ell}}\right).
\end{align}
\begin{Assumption}\label{assumption:causal}
The sequence $\{\epsilon_{j}\}_{j\in \mathbb{Z}}$
admits the causal representation \cite{wu2005nonlinear}
\begin{align}\label{causal_rep}
\epsilon_{j} = H_{j}(\mathcal{F}_{j}),
\end{align}
where $\mathcal{F}_{j}=(\ldots,e_{j-2},e_{j-1},e_{j})$, $\{e_{j}\}_{j\in \mathbb{Z}}$  is a sequence of i.i.d. random variables, and
$H_{j}: \mathbb{R}^{\infty}\mapsto \mathbb{R}$ is a measurable function
describing the causal filtration mechanism.
\end{Assumption}
The function $H_j$ models the noise generation process. When $H_j=H$ for a fixed function $H$,
the noise generation process is fixed, and the random sequence $\{\epsilon_j\}_{j\in \mathbb{Z}}$
is stationary.
In practice, this noise generation process is meaningless without further assumptions.
We focus on two controls: moment conditions and temporal dependence.
We follow the dependence measure proposed in \cite{wu2005nonlinear} to control the dependence structure of $\{\epsilon_j\}_{j\in \mathbb{Z}}$.
Let $\{\widehat{e}_{j}\}_{j\in \mathbb{Z}}$ be an independent copy of $\{e_{j}\}_{j\in \mathbb{Z}}$.
For $k\in \mathbb{N}\cup\{0\}$, denote $$\mathcal{F}_{j,j-k}=(\ldots,e_{j-k-1},\widehat{e}_{j-k},e_{j-k+1},\ldots,e_{j-2},e_{j-1},e_{j}).$$

\begin{figure}[hbt!]
  \centering
  \includegraphics[scale=0.7]{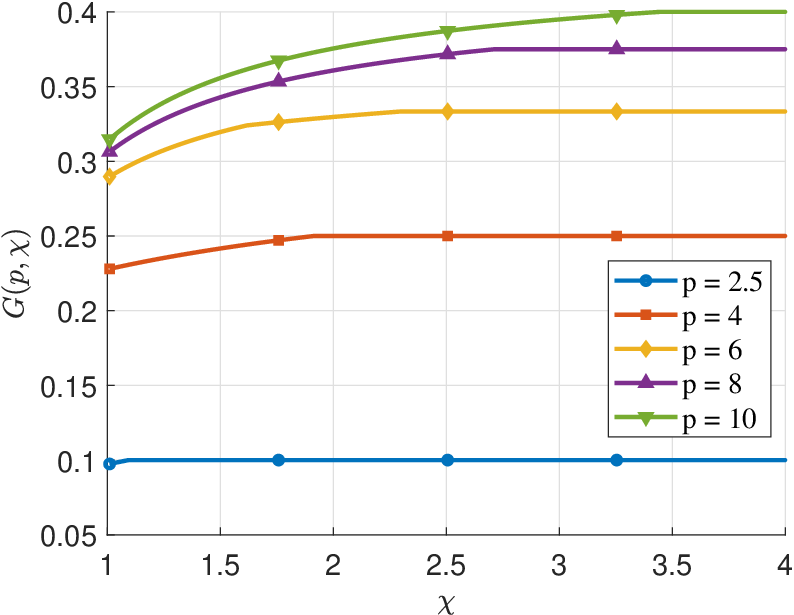}
  \caption{Graph of function $G(p,\cdot)$ for several values of $p$.}\label{fig:graph:G}
\end{figure}

\begin{Assumption}\label{assumption:weak}
Assume that $\sup_{j\in \mathbb{Z}}\ \mathbb{E}[|\epsilon_{j}|^{p}]<\infty$ for some $p>2$, and that
\begin{align}\label{def dependence measure}
&\sup_{j\in \mathbb{Z}}\ \left(\mathbb{E}\left[\left|H_{j}(\mathcal{F}_{j})-H_{j}(\mathcal{F}_{j,j-k})\right|^{p}\right]\right)^{1/p}
=\mathcal{O}\left( (k+1)^{-(\chi+1)}(\log(k+1))^{-a}\right)\nonumber
\end{align}
for some constants $\chi>1$ and $a>\sqrt{\chi}+1$; that is, the influence of the $k$-step historical input $e_{j-k}$ on the output $\epsilon_{j}$ decays polynomially.
\end{Assumption}
Define a function $G:(2,\infty)\times (1,\infty)\mapsto (0,\infty)$ by \cite{wu2025uncertainty}
$$G(p,\chi) = \frac{1}{1-\frac{1}{\kappa}-\frac{1}{p}} \left(\frac{1}{2}-\frac{1}{p}-\frac{1}{\kappa}+\frac{2}{p\kappa}\right)
$$
where
\begin{equation*}
\begin{aligned}
\frac{1}{\kappa}
= \max\Biggl\{
    \frac{1}{p},\;
    \frac{1}{\sqrt{\chi+1}\,p}
    &+
    \left(\frac12 - \frac{1}{\sqrt{\chi+1}\,p}\right)
    \max\left\{
        \frac{2}{\sqrt{\chi+1}\,p},\;
        \frac{1}{\chi}\!\left(\frac12 - \frac{1}{p}\right)
    \right\}
\Biggr\}.
\end{aligned}
\end{equation*}
Note that $G(p,\chi)\in(0,1/2-1/p)$ for $p>2$ and $\chi>1$.
The graph of $G(p,\cdot)$ for several values of $p$ is shown in Figure \ref{fig:graph:G}.

Recall that
the L$\acute{\textup{e}}$vy-Prokhorov distance between two probability measures
$\nu_{1}$ and $\nu_{2}$ on $\mathbb{C}^{n\times d}$ is defined by
\begin{align*}
\pi\!\left(\nu_{1},\nu_{2}\right)
= \inf\left\{\varepsilon>0:\;
\begin{aligned}
&\nu_{1}(B)
   \le \nu_{2}(B^{\varepsilon})+\varepsilon,\\
&\nu_{2}(B)
   \le \nu_{1}(B^{\varepsilon})+\varepsilon,\\
&\forall\, B\in\mathcal{B}(\mathbb{C}^{n\times d})
\end{aligned}
\right\},
\end{align*}
where $\mathcal{B}(\mathbb{C}^{n\times d})$  denotes the Borel $\sigma$-algebra on $\mathbb{C}^{n\times d}$, and $B^{\varepsilon}=\{x\in \mathbb{C}^{n\times d}\mid \underset{x'\in B}{\inf} |x-x'|_{\infty}<\varepsilon\}$
is the $\varepsilon$-enlargement of $B$ with respect to the supremum-norm metric. We have the following Gaussian approximation result for the discretized AWT.

\begin{Proposition}\label{prop:approximation}
Suppose that the mother wavelet $\psi$ is absolutely continuous with $\psi'\in L^{1}$.
Fix $0<s_1<\ldots<s_d$.
Under Assumptions  \ref{assumption:causal} and \ref{assumption:weak},
there exist mean-zero independent Gaussian random
variables $\{z_{j}\}_{j\in \in \mathbb{Z}}$ satisfying $\textup{Var}(z_{j})=\textup{Var}(\epsilon_{j})$
for all $j\in \mathbb{Z}$ such that
the L$\acute{\textup{e}}$vy-Prokhorov distance $\pi(\cdot,\cdot)$ between the probability measures of $\{\sqrt{s_{\ell}}\widetilde{W}_{\epsilon}(t_{i},s_{\ell})\mid 1\leq i\leq n,\ 1\leq \ell \leq d\}$
and $\{\sqrt{s_{\ell}}\widetilde{W}_{z}(t_{i},s_{\ell})\mid 1\leq i\leq n,\ 1\leq \ell \leq d\}$,
denoted respectively by $\nu_{\epsilon}^{(n)}$
and $\nu_{z}^{(n)}$, satisfies
\begin{align}\label{prop:LP_bound}
\pi(\nu_{\epsilon}^{(n)},\nu_{z}^{(n)})=o\left(n^{-G(p,\chi)/3}(\log n)^{4/3+\delta}\right)
\end{align}
as $n$ tends to infinity,
where $\delta>0$ is arbitrary, and $\widetilde{W}_{z}$ is defined in the same way as in (\ref{def:discreteCWT2})
but with the sequence $\{\epsilon_{j}\}_{j\in \mathbb{Z}}$ replaced by the Gaussian sequence
$\{z_{j}\}_{j\in \in \mathbb{Z}}$.
\end{Proposition}

The proof is provided in Section \ref{sec:proof:prop:approximation}.
Proposition \ref{prop:approximation} shows that, under the moment and dependence control assumptions, the distribution of the discretized AWT of a causal non-Gaussian time series, which can be nonstationary, can be approximated by that of a suitably constructed Gaussian sequence asymptotically as $n\to \infty$, with an explicit convergence rate as the sample size increases.
The appearance of the factor $\sqrt{s_{\ell}}$ in front of $\widetilde{W}_{\epsilon}(t_{i},s_{\ell})$ and $\widetilde{W}_{z}(t_{i},s_{\ell})$ shows that the Gaussian approximation of the distribution of $\widetilde{W}_{\epsilon}(t_{i},s_{\ell})$ by that of $\widetilde{W}_{z}(t_{i},s_{\ell})$ becomes increasingly accurate as the scale $s_{\ell}$ grows. A scale-dependent improvement in Gaussian approximation has also been observed in \cite{fageot2015wavelet}.

To take a closer look at the bound (\ref{prop:LP_bound}), for $i=1,\ldots,n$, define the $d$-dimensional vector
$$\widetilde{\mathbf{W}}_{\widehat{\epsilon}}(t_{i})= \left[\widetilde{W}_{\widehat{\epsilon}}(t_{i},s_{1})\  \widetilde{W}_{\widehat{\epsilon}}(t_{i},s_{2})\ \cdots\
\widetilde{W}_{\widehat{\epsilon}}(t_{i},s_{d})\right]\,,
$$
where $\{\widehat{\epsilon}_{j}\}_{j\in \mathbb{Z}}$ is a random sequence with the same distribution as $\{\epsilon_{j}\}_{j\in \mathbb{Z}}$.
The vector $\widetilde{\mathbf{W}}_{z}(t_{i})$ is defined analogously.
In parallel with \cite[Theorem~3.2]{wu2025uncertainty}, if
the scale values are bounded, that is, $\{s_{\ell}\}_{1\leq \ell\leq d}\subset [s_{*},s^{*}]$ for some constants  $s_{*}$ and $s^{*}\in(0,\infty)$,
then the Euclidean norm in $\mathbb{R}^{d}$ of $\widetilde{\mathbf{W}}_{\widehat{\epsilon}}(t_{i})-\widetilde{\mathbf{W}}_{z}(t_{i})$
admits the following estimate
\begin{align}\notag
\mathbb{E}\left[\underset{1\leq i\leq n}{\max}\left|\widetilde{\mathbf{W}}_{\widehat{\epsilon}}(t_{i})-\widetilde{\mathbf{W}}_{z}(t_{i})\right|^{2}\right]
=\mathcal{O}\left(dn^{-G(p,\chi)}\left(\log\left(n\right)\right)^{4}\right),
\end{align}
which means that the approximation error converges to zero not only in the case of a fixed number $d$ of scales,
but also when
$d$ grows with
$n$ according to $d=n^{G(p,\chi)-\vartheta}$ for some $\vartheta>0$.
When $p,\chi\to \infty$ so that $G(p,\chi)\to 1/2$, the number of scales we can handle is $d=n^{1/2-\vartheta}$ for some $\vartheta>0$. Note that this rate is faster than that of discretized short-time Fourier transform shown in \cite[Theorem~3.2]{wu2025uncertainty}.
See Remark \ref{remark:approximation} at the end of Section \ref{sec:proof:prop:approximation} for details.

We shall mention that Proposition \ref{prop:approximation} provides a link between the discretized AWT of non-Gaussian noise and the continuous-time AWT of Gaussian noise, in a manner analogous to the relationship between stochastic differential equations and their discretized stochastic difference equations driven by non-Gaussian innovations \cite{Unser2014I,Unser2014II}.
Future work could focus on achieving a better convergence rate and extending the results above to more general noise models. Since nonstationary and non-Gaussian noise is everywhere, we also need to study what is the suitable continuous model whose discretization satisfies the causal condition and Assumption \ref{assumption:weak}, and establish the relevant results.

\section{Conclusion}\label{sec:conclusion}
This paper presented a statistical framework for analyzing the magnitude and phase of the AWT of noisy signals. Modeling the observation as a deterministic signal corrupted by stationary Gaussian noise, we derived joint and marginal probability density functions for the magnitude and phase in the general and noise-only settings. The analysis yielded SNR-dependent concentration bounds, explicit covariance expressions, and independence results, offering mathematical foundations for their distributional and dependence structures. We further connected these results to practical tasks by characterizing ridge detection error probabilities, the regularity of scalogram contours and the relationship between the magnitude and phase.
 In addition, we extended the framework to non-Gaussian noise by establishing a Gaussian approximation for the discretized AWT of weakly dependent noise sequences.
The results lay the groundwork for developing rigorous AWT-based algorithms and wavelet-domain statistical inference.

\section{Proofs}\label{sec:proof}

\subsection{Proof of Lemma \ref{prop:circular}}\label{sec:proof:prop:circular}
Because $\psi$ is a wavelet and the mean $\mathbb{E}[\Phi(t)]$ is constant in time,
the zero-mean property $\int_{\mathbb{R}}\psi(u)du=0$ implies that
\begin{align}\label{WPhi_minus_mean}
W_{\Phi}(t,s)=\int_{\mathbb{R}}\left(\Phi(\tau)-\mathbb{E}[\Phi(\tau)]\right)\frac{1}{\sqrt{s}}\overline{\psi\left(\frac{\tau-t}{s}\right)}d\tau.
\end{align}
By substituting (\ref{spectral_rep:samplepath}) into (\ref{WPhi_minus_mean}),
the AWT of $\Phi$ has the representation
\begin{align}\label{spect:WPhi}
W_{\Phi}(t,s) = \sqrt{s}\int_{\mathbb{R}}e^{it\lambda} \overline{\widehat{\psi}(s\lambda)}Z_{F}(d\lambda),\ t\in \mathbb{R},s>0.
\end{align}
For any $\ell,\ell'\in\{1,2,...,n\}$,
using (\ref{spect:WPhi}), we have
\begin{align}\notag
(s_{\ell}s_{\ell'})^{-\frac{1}{2}}\Gamma(\ell,\ell')
=& (s_{\ell}s_{\ell'})^{-\frac{1}{2}}\mathbb{E}\left[W_{\Phi}(t_{\ell},s_{\ell})\overline{W_{\Phi}(t_{\ell'},s_{\ell'})}\right]
\\\notag=&\mathbb{E}\left[\int_{\mathbb{R}} e^{it_{\ell}\lambda} \overline{\widehat{\psi}(s_{\ell}\lambda)} Z_{F}(d\lambda)\int_{\mathbb{R}} e^{-it_{\ell'}\lambda'} \widehat{\psi}(s_{\ell'}\lambda') \overline{Z_{F}(d\lambda')}\right]
\\\label{entry_Gamma_pre}=&
\int_{\mathbb{R}}\int_{\mathbb{R}} e^{it_{\ell}\lambda} \overline{\widehat{\psi}(s_{\ell}\lambda)}  e^{-it_{\ell'}\lambda'} \widehat{\psi}(s_{\ell'}\lambda') \mathbb{E}\left[Z_{F}(d\lambda)\overline{Z_{F}(d\lambda')}\right].
\end{align}
The equation (\ref{entry_Gamma}) follows by applying (\ref{ortho}) to (\ref{entry_Gamma_pre}).
On the other hand,
\begin{align}\notag
(s_{\ell}s_{\ell'})^{-\frac{1}{2}}C_{\ell\ell'}
=&(s_{\ell}s_{\ell'})^{-\frac{1}{2}}\mathbb{E}\left[W_{\Phi}(t_{\ell},s_{\ell})W_{\Phi}(t_{\ell'},s_{\ell'})\right]
\\\notag=&\mathbb{E}\left[\int_{\mathbb{R}} e^{it_{\ell}\lambda} \overline{\widehat{\psi}(s_{\ell}\lambda)} Z_{F}(d\lambda)\hspace{-0.15cm}\int_{\mathbb{R}} e^{it_{\ell'}\lambda'} \overline{\widehat{\psi}(s_{\ell'}\lambda')} Z_{F}(d\lambda')\right]
\\\label{entry_pseudoC_pre}=&
\int_{\mathbb{R}}\int_{\mathbb{R}} e^{it_{\ell}\lambda} \overline{\widehat{\psi}(s_{\ell}\lambda)}  e^{it_{\ell'}\lambda'} \overline{\widehat{\psi}(s_{\ell'}\lambda')} \mathbb{E}\left[Z_{F}(d\lambda)Z_{F}(d\lambda')\right].
\end{align}
By (\ref{ortho}), for any measurable sets $\Delta_{1},\Delta_{2}\subset \mathbb{R}$,
\begin{align}\notag
\mathbb{E}\left[Z_{F}(\Delta_{1})
Z_{F}(\Delta_{2})\right]=&
\mathbb{E}\left[Z_{F}(\Delta_{1})
\overline{Z_{F}(-\Delta_{2})}\right]
\\\label{EZFZF}=&
F(\Delta_{1}\cap(-\Delta_{2})).
\end{align}
Applying (\ref{EZFZF}) to (\ref{entry_pseudoC_pre}) yields (\ref{entry_pseudoC}).
Under Assumption \ref{assump:analytic}, we have $C_{\ell\ell'}=0$ because
the supports of $\widehat{\psi}(s_{\ell}\cdot)$ and $\widehat{\psi}(-s_{\ell'}\cdot)$
are disjoint.
\qed

\subsection{Proof of Proposition \ref{prop:jpdf:Sn}}\label{sec:proof:prop:jpdf:Sn}
(a) Equation~\eqref{jpdf:|WY|:n:nonnull} is obtained by expressing the variable $\mathbf{w}$
in the probability density function~(\ref{pdf_WY}) of $W_{Y}$ in polar coordinates and integrating with respect to the angular components.

(b) In the case $n=1$,
\begin{align}\notag
p\left(r;|W_{Y}(t,s)|\right)
= &\frac{r}{\pi \mathbb{E}[|W_{\Phi}(t,s)|^2]}
\int_{-\pi}^{\pi}\textup{exp}\left(-\frac{\left|re^{i\theta}-W_{f}(t,s)\right|^2}{\mathbb{E}[|W_{\Phi}(t,s)|^2]}\right)d\theta
\\\notag=&\frac{r}{\pi \mathbb{E}[|W_{\Phi}(t,s)|^2]}
\textup{exp}\left(-\frac{r^2+|W_{f}(t,s)|^2}{\mathbb{E}[|W_{\Phi}(t,s)|^2]}\right)
\\\notag&\times\int_{-\pi}^{\pi}\textup{exp}\left(\frac{2r\Re(e^{i\theta}\overline{W_{f}(t,s)})}{\mathbb{E}[|W_{\Phi}(t,s)|^2]}\right)
d\theta.
\end{align}
By expressing $W_{f}(t,s)=|W_{f}(t,s)|e^{i\varphi}$, where $\varphi\in [0,2\pi)$, we obtain
\begin{align}\notag
p\left(r;|W_{Y}(t,s)|\right)
=&
\frac{r}{\pi \mathbb{E}[|W_{\Phi}(t,s)|^2]}
\textup{exp}\left(-\frac{r^2+|W_{f}(t,s)|^2}{\mathbb{E}[|W_{\Phi}(t,s)|^2]}\right)
\\\notag&\times\int_{-\pi}^{\pi}\textup{exp}\left(\frac{2|W_{f}(t,s)|r}{\mathbb{E}[|W_{\Phi}(t,s)|^2]}\cos(\theta-\varphi)\right)d\theta
\\\notag=&\frac{2r}{\mathbb{E}[|W_{\Phi}(t,s)|^2]}
\textup{exp}\left(-\frac{r^2+|W_{f}(t,s)|^2}{\mathbb{E}[|W_{\Phi}(t,s)|^2]}\right)
I_{0}\left(\frac{2|W_{f}(t,s)|r}{\mathbb{E}[|W_{\Phi}(t,s)|^2]}
\right).
\end{align}

(c) To prove (\ref{jpdf:|WY|:n=2:nonnull}), we consider the decomposition
\begin{align*}
W_{Y}(t_{2},s_{2}) = V+\frac{\Gamma_{12}}{\Gamma_{11}}W_{\Phi}(t_{1},s_{1})+W_{f}(t_{2},s_{2}),
\end{align*}
where
\begin{align*}
V = W_{\Phi}(t_{2},s_{2})-\frac{\Gamma_{12}}{\Gamma_{11}} W_{\Phi}(t_{1},s_{1}).
\end{align*}
Note that $V$ is a complex Gaussian random variable satisfying
$$\mathbb{E}[V]=0,\
\mathbb{E}[V^{2}]=0,\ \textup{and}\
\mathbb{E}[|V|^{2}]=\Gamma_{22}-\frac{|\Gamma_{12}|^{2}}{\Gamma_{11}}.$$
Moreover, $V$ is independent of $W_{\Phi}(t_{1},s_{1})$.
Conditioned on $W_{\Phi}(t_{1},s_{1})$,
$|W_{Y}(t_{2},s_{2})|$ has the Rician distribution with
noncentrality parameter $m(W_{\Phi}(t_{1},s_{1}))$, which is defined as
$$
m\left(W_{\Phi}(t_{1},s_{1})\right)=\left|\frac{\Gamma_{12}}{\Gamma_{11}}W_{\Phi}(t_{1},s_{1})+W_{f}(t_{2},s_{2})\right|,
$$
and spread parameter $\mathbb{E}[|V|^{2}]/2$.
That is, the conditional probability density function is given by
\begin{align*}
p\left(r_{2};|W_{Y}(t_{2},s_{2})|\mid W_{\Phi}(t_{1},s_{1})\right)
=
p_{\textup{Rice}}\left(r_{2};m\left(W_{\Phi}(t_{1},s_{1})\right),
\frac{1}{2}\mathbb{E}[|V|^{2}]\right)
\end{align*}
for $r_{2}\geq0$.
On the other hand, conditioned on $W_{\Phi}(t_{1},s_{1})$, $|W_{Y}(t_{1},s_{1})|$ is deterministically equal to $|W_{f}(t_{1},s_{1})+W_{\Phi}(t_{1},s_{1})|$.
This implies that, conditioned on $W_{\Phi}(t_{1},s_{1})$,
the random variables $|W_{Y}(t_{1},s_{1})|$ and $|W_{Y}(t_{2},s_{2})|$ are independent.
Hence, their joint conditional probability density function is given by
\begin{align}\notag
&p\left(r_{1},r_{2}; |W_{Y}(t_{1},s_{1})|,|W_{Y}(t_{2},s_{2})|\mid W_{\Phi}(t_{1},s_{1})\right)
\\\notag=&\delta\left(r_{1}-\left|W_{f}(t_{1},s_{1})+W_{\Phi}(t_{1},s_{1})\right|\right)
 p_{\textup{Rice}}\left(r_{2};m\left(W_{\Phi}(t_{1},s_{1})\right),
\frac{1}{2}\mathbb{E}[|V|^{2}]\right),
\end{align}
where $\delta$ denotes the Dirac delta distribution on $\mathbb{R}$.
From the laws of probability,
\begin{align}\notag
&p\left(r_{1},r_{2}; |W_{Y}(t_{1},s_{1})|,|W_{Y}(t_{2},s_{2})|\right)
\\\notag=&\int_{\mathbb{C}}\delta\left(r_{1}-\left|W_{f}(t_{1},s_{1})+z\right|\right)
 p_{\textup{Rice}}\left(r_{2};m\left(z\right),
\frac{1}{2}\mathbb{E}[|V|^{2}]\right)
 \frac{1}{\pi\Gamma_{11}}
\textup{exp}\left(-\frac{|z|^{2}}{\Gamma_{11}}\right)dz.
\end{align}
By a change of variables,
\begin{align}\notag
&p\left(r_{1},r_{2}; |W_{Y}(t_{1},s_{1})|,|W_{Y}(t_{2},s_{2})|\right)
\\\notag=&\frac{r_{1}}{\pi\Gamma_{11}}\int_{0}^{2\pi}
 \hspace{-0.2cm}p_{\textup{Rice}}\left(r_{2};m\left(r_{1}e^{i\theta}-W_{f}(t_{1},s_{1})\right),
\frac{1}{2}\mathbb{E}[|V|^{2}]\right)
\textup{exp}\left(-\frac{|r_{1}e^{i\theta}-W_{f}(t_{1},s_{1})|^{2}}{\Gamma_{11}}\right)d\theta.
\end{align}
\qed

\subsection{Proof of Corollary \ref{prop:jpdf:S2}}\label{sec:proof:prop:jpdf:S2}
(a) Consider the polar representation of
$W_{\Phi}(t_{\ell},s_{\ell})$, where $\ell\in\{1,2,...,n\}$:
\begin{align*}
W_{\Phi}(t_{\ell},s_{\ell})= |W_{\Phi}(t_{\ell},s_{\ell})|e^{i\Theta_{\ell}},
\end{align*}
where $\Theta_{\ell}=\arg \left(W_{\Phi}(t_{\ell},s_{\ell})\right)\in [0,2\pi)$.
Denote
$$\mathbf{W}_{\Phi}:=[W_{\Phi}(t_{1},s_{1})\ W_{\Phi}(t_{2},s_{2})\ \cdots\ W_{\Phi}(t_{n},s_{n})]^{\top},$$
$|\mathbf{W}_{\Phi}|:=[|W_{\Phi}(t_{1},s_{1})|\ |W_{\Phi}(t_{2},s_{2})|\ \cdots\ |W_{\Phi}(t_{n},s_{n})|]^{\top},$
and $\Theta=[\Theta_{1}\ \Theta_{2}\ \cdots\ \Theta_{n}]^{\top}$.
Because the probability density function of $\mathbf{W}_{\Phi}$ is given by
\begin{align*}
p\left(\mathbf{w};\mathbf{W}_{\Phi}\right)
=
\frac{1}{\pi^{n}\textup{det}(\Gamma)}
\textup{exp}\left(-\mathbf{w}^{*}\Gamma^{-1}\mathbf{w}\right), \ \mathbf{w}\in \mathbb{C}^{n},
\end{align*}
where $\Gamma$ is defined in (\ref{def:Gamma}) and (\ref{entry_Gamma}),
\begin{align}\label{jpdf:R&Theta}
p\left(r_{1:n},\theta_{1:n};|\mathbf{W}_{\Phi}|,\Theta\right)
=
\frac{r_{1}r_{2}\cdots r_{n}}{\pi^{n}\textup{det}(\Gamma)}
\textup{exp}\left(-\begin{bmatrix}r_{1}e^{i\theta_{1}} \\ r_{2}e^{i\theta_{2}}\\ \vdots \\ r_{n}e^{i\theta_{n}} \end{bmatrix}^{*}\Gamma^{-1}\left[\begin{array}{cc}r_{1}e^{i\theta_{1}}\\ r_{2}e^{i\theta_{2}} \\\vdots\\ r_{n}e^{i\theta_{n}}\end{array}\right]\right),
\end{align}
where  $r_{1:n}=(r_{1},r_{2},...,r_{n})\in \mathbb{R}_{\geq 0}^{n},\ \theta_{1:n}=(\theta_{1},\theta_{2},...,\theta_{n})\in[0,2\pi)^{n}.$
By  integrating  (\ref{jpdf:R&Theta}) with respect to the angular components, we obtain
\begin{align}\notag
p\left(r_{1:n};|\mathbf{W}_{\Phi}|\right)
=
\frac{r_{1}r_{2}\cdots r_{n}}{\pi^{n}\textup{det}(\Gamma)}
J\left(\textup{diag}(r_{1:n})\Gamma^{-1}\textup{diag}(r_{1:n})\right),
\end{align}
where the function $J$ is  defined in (\ref{def:functionB}).

(b)
For any Hermitian matrix $\mathbf{M}=[M_{\ell\ell'}]_{1\leq\ell,\ell'\leq 2}$,
\begin{align}\notag
\begin{bmatrix} e^{-i\theta_{1}} & e^{-i\theta_{2}}\end{bmatrix}\mathbf{M}\begin{bmatrix} e^{i\theta_{1}} \\ e^{i\theta_{2}}\end{bmatrix}
=&
M_{11} +  M_{22} +  \left(M_{12} e^{i(\theta_2 - \theta_1)} + M_{21} e^{-i(\theta_2 - \theta_1)}\right)
\\\notag=&M_{11} +  M_{22} +2 \Re \left(M_{12} e^{i(\theta_2 - \theta_1)}\right).
\end{align}
Hence, the function $J$, defined in (\ref{def:functionB}) with $n=2$, simplifies to
\begin{align}\label{B(M):n=2}
J(\mathbf{M})=&e^{- (M_{11} +  M_{22})}
\int_0^{2\pi} \int_0^{2\pi} \exp \left(- 2 \Re \left(M_{12} e^{i(\theta_2 - \theta_1)}\right) \right) d\theta_1 d\theta_2.
\end{align}
Let $M_{12} = |M_{12}| e^{i\phi}$. Because
$$
\Re\left(M_{12} e^{i(\theta_2 - \theta_1)}\right)= |M_{12}| \cos(\theta_2 - \theta_1 + \phi),
$$
(\ref{B(M):n=2}) implies that
\begin{align}\notag
J(\mathbf{M}) =& e^{- (M_{11} +  M_{22})}
\int_0^{2\pi}\hspace{-0.2cm} \int_0^{2\pi}\hspace{-0.2cm} \exp \left( -2 |M_{12}| \cos(\theta_2 - \theta_1 + \phi) \right) d\theta_1 d\theta_2.
\end{align}
By the formula
\[
\int_0^{2\pi} e^{-2 |M_{12}| \cos \theta} d\theta = 2\pi I_0\left(2 |M_{12}|\right),
\]
\begin{align}\label{B(M):n=2:simplify1}
J(\mathbf{M}) =
4\pi^2 e^{- (M_{11} +  M_{22})} I_0\left(2 |M_{12}|\right).
\end{align}
By substituting (\ref{B(M):n=2:simplify1}) into (\ref{jpdf:S:n}) with
\begin{align}\notag
\mathbf{M}= \textup{diag}(r_{1:2})\Gamma^{-1}\textup{diag}(r_{1:2})
=\frac{1}{\textup{det}(\Gamma)}
\begin{bmatrix}
r_{1}^{2}\Gamma_{22} & -r_{1}r_{2}\Gamma_{12} \\
-r_{1}r_{2}\overline{\Gamma_{12}} & r_{2}^{2}\Gamma_{11}
\end{bmatrix},
\end{align}
we obtain
\begin{align}\notag
p\left(r_{1:2};|\mathbf{W}_{\Phi}|\right)
=
\frac{4r_{1}r_{2}}{\textup{det}(\Gamma)}\exp\left(- \frac{r_{1}^{2}\Gamma_{22} +  r_{2}^{2}\Gamma_{11}}{\textup{det}(\Gamma)}\right) I_0\left(\frac{2|\Gamma_{12}|}{\textup{det}(\Gamma)}r_{1}r_{2}\right).
\end{align}
\qed

\subsection{Proof of Proposition \ref{lemma:convprob}}\label{sec:proof:lemma:convprob}
By denoting $\varphi=\textup{arg}(W_{f}(t,s))$,
\begin{align}\label{expan:WX}
|W_{Y}(t,s)|^2
=&|W_{f}(t,s)+ W_{\Phi}(t,s)|^{2}
\\\notag=&\left||W_{f}(t,s)|e^{i\varphi}+ W_{\Phi}(t,s)\right|^{2}
\\\notag=&\left||W_{f}(t,s)|+ N_{R}+iN_{I}\right|^{2}
\\\notag=&|W_{f}(t,s)|^{2}+2|W_{f}(t,s)|N_{R} +N^{2}_{R}+N^{2}_{I},
\end{align}
where
\begin{align}\label{def:NR}
N_{R} = \Re\left(e^{-i\varphi}W_{\Phi}(t,s)\right)
=\cos(\varphi)\Re(W_{\Phi}(t,s))+\sin(\varphi)\Im(W_{\Phi}(t,s))
\end{align}
and
\begin{align}\label{def:NI}
N_{I} =& \Im\left(e^{-i\varphi}W_{\Phi}(t,s)\right)
=\cos(\varphi)\Im(W_{\Phi}(t,s))-\sin(\varphi)\Re(W_{\Phi}(t,s)).
\end{align}
Using (\ref{expan:WX}), for any $\varepsilon>0$, we have
\begin{align}\notag
\mathbb{P}\left(\left|\frac{|W_{Y}(t,s)|^2}{|W_{f}(t,s)|^2}-1\right|>\varepsilon\right)
=&\mathbb{P}\left(\left|\frac{2N_{R}}{|W_{f}(t,s)|} +\frac{N^{2}_{R}+N^{2}_{I}}{|W_{f}(t,s)|^{2}}\right|>\varepsilon\right)
\\\label{markov0}\leq&\mathbb{P}\left(\frac{|N_{R}|}{|W_{f}(t,s)|}>\frac{\varepsilon}{4}\right)
+
\mathbb{P}\left(\frac{N^{2}_{R}+N^{2}_{I}}{|W_{f}(t,s)|^{2}}>\frac{\varepsilon}{2}\right).
\end{align}
As demonstrated in the proof of Lemma \ref{prop:circular}, $\Re(W_{\Phi}(t,s))$ and $\Im(W_{\Phi}(t,s))$ are independent normal random variables. Consequently, (\ref{def:NR}) and (\ref{def:NI}) imply that
the random vector $(N_{R},N_{I})$ follows a normal distribution with a mean of zero and a covariance matrix given by
\begin{align*}
\mathbb{E}\left[\begin{pmatrix}N_{R}N_{R} & N_{R}N_{I}\\ N_{I}N_{R} & N_{I}N_{I}\end{pmatrix}\right] =
\begin{pmatrix}\sigma_{s}^{2} & 0\\ 0 & \sigma_{s}^{2}\end{pmatrix},
\end{align*}
where
\begin{align*}
\sigma_{s}^{2} = \frac{1}{2}\mathbb{E}\left[|W_{\Phi}(t,s)|^2\right]=\frac{s}{2}\int_{0}^{\infty}|\widehat{\psi}(s\lambda)|^{2}F(d\lambda)
\end{align*}
and $F$ is the spectral measure of the covariance function of $\Phi$.
Hence,  $N_{R}/\sigma_{s}$ follows a standard normal distribution,
and $(N^{2}_{R}+N^{2}_{I})/\sigma_{s}^{2}$ follows a chi-square distribution with two degrees of freedom, i.e., an exponential distribution with a mean of two.
These observations enable us to estimate the right-hand side of (\ref{markov0}) as follows
\begin{align}\notag
 \mathbb{P}\left(\frac{|N_{R}|}{|W_{f}(t,s)|}>\frac{\varepsilon}{4}\right)
=&\mathbb{P}\left(\frac{|N_{R}|}{\sigma_{s}}>\frac{\varepsilon}{4}\frac{|W_{f}(t,s)|}{\sigma_{s}}\right)
\\\label{markov1}\leq& \textup{exp}\left(-\frac{\varepsilon^{2}}{32}\frac{|W_{f}(t,s)|^{2}}{\sigma_{s}^{2}}\right)
\end{align}
and
\begin{align}\notag
\mathbb{P}\left(\frac{N^{2}_{R}+N^{2}_{I}}{|W_{f}(t,s)|^{2}}>\frac{\varepsilon}{2}\right)
=&\mathbb{P}\left(\frac{N^{2}_{R}+N^{2}_{I}}{\sigma_{s}^{2}}>\frac{\varepsilon}{2}\frac{|W_{f}(t,s)|^{2}}{\sigma_{s}^{2}}\right)
\\\label{markov2}=&\textup{exp}\left(-\frac{\varepsilon}{4}\frac{|W_{f}(t,s)|^{2}}{\sigma_{s}^{2}}\right).
\end{align}
By substituting (\ref{markov1}) and (\ref{markov2}) into (\ref{markov0}),
\begin{align}\label{markov3}
\mathbb{P}\left(\left|\frac{|W_{Y}(t,s)|^2}{|W_{f}(t,s)|^2}-1\right|>\varepsilon\right)
\leq\textup{exp}\left(-\frac{\varepsilon^{2}}{32}\frac{|W_{f}(t,s)|^{2}}{\sigma_{s}^{2}}\right)
+
\textup{exp}\left(-\frac{\varepsilon}{4}\frac{|W_{f}(t,s)|^{2}}{\sigma_{s}^{2}}\right).
\end{align}
Proposition \ref{lemma:convprob} follows from (\ref{markov3})
and the inequality
\begin{align}\notag
\mathbb{P}\left(\left|\frac{|W_{Y}(t,s)|^2}{|W_{f}(t,s)|^2}-1\right|>\varepsilon\right)
=&\mathbb{P}\left(\left|\frac{|W_{Y}(t,s)|}{|W_{f}(t,s)|}-1\right|>
\left(\frac{|W_{Y}(t,s)|}{|W_{f}(t,s)|}+1\right)^{-1}\varepsilon\right)
\\\notag\geq&\mathbb{P}\left(\left|\frac{|W_{Y}(t,s)|}{|W_{f}(t,s)|}-1\right|>
\varepsilon\right).
\end{align}
\qed

\subsection{Proof of Proposition \ref{prop:ridge}}\label{sec:proof:prop:ridge}

First of all, we define the following events:
\begin{align}\label{def:E}
E=\left\{|W_{Y}(t,s_{f}(t))|<(1-\delta) |W_{Y}(t,s)|\right\},
\end{align}
\begin{align}\label{def:Es}
E_{s}=\left\{\left|\frac{|W_{Y}(t,s)|}{|W_{f}(t,s)|}-1\right|<\varepsilon\right\},
\end{align}
and
\begin{align}\label{def:Esf}
E_{s_{f}(t)}=\left\{\left|\frac{|W_{Y}(t,s_{f}(t))|}{|W_{f}(t,s_{f}(t))|}-1\right|<\varepsilon\right\}.
\end{align}
The definition of $\varepsilon$ in (\ref{threshold}) implies that $\mathbb{P}\left(E\cap E_{s}\cap E_{s_{f}(t)}\right)=0$.
Denote the complement of $E_{s}$ (resp. $E_{s_{f}(t)}$) by  $E_{s}^{c}$ (resp. $E^{c}_{s_{f}(t)}$).
From (\ref{def:E}), (\ref{def:Es}), and (\ref{def:Esf}), we obtain
\begin{align}\notag
\mathbb{P}\left(E\right)
=& \mathbb{P}\left(E\cap E_{s}\cap E_{s_{f}(t)}\right)+\mathbb{P}\left(E\cap (E_{s}\cap E_{s_{f}(t)})^{c}\right)
\\\label{P(E)}\leq&\mathbb{P}\left(E_{s}^{c}\right)
+\mathbb{P}(E_{s_{f}(t)}^{c}).
\end{align}
By applying Proposition \ref{lemma:convprob} to $\mathbb{P}(E_{s}^{c})$
and $\mathbb{P}(E_{s_{f}(t)}^{c})$
in (\ref{P(E)}),
\begin{align}\notag
\mathbb{P}\left(E\right)\leq&
\textup{exp}\left(-\frac{\varepsilon^{2}}{16}q(t,s_{f}(t))\right)
+\textup{exp}\left(-\frac{\varepsilon}{2}q(t,s_{f}(t))\right)
\\\notag&+\textup{exp}\left(-\frac{\varepsilon^{2}}{16}q(t,s)\right)
+\textup{exp}\left(-\frac{\varepsilon}{2}q(t,s)\right).
\end{align}
\qed

\subsection{Proof of Proposition \ref{thm:contour:c1}}\label{sec:proof:thm:contour:c1}
We first recall two theorems which will be used in the proof.

\noindent{\bf Sard's Theorem} \cite[Chapter 6]{lee2013introduction}
{\it Let \( g : \mathbb{R}^n \to \mathbb{R}^m \) be a function of class \( C^k \), where \( k \ge \max(n - m + 1, 1) \). Then the set of \emph{critical values} of \(g \), that is,
\[
\{ g(x) \in \mathbb{R}^m : Dg(x) \text{ is not of full rank} \},
\]
has Lebesgue measure zero in \( \mathbb{R}^m \). Here, $Dg(x)$ represents the Jacobian matrix of $g$ at $x$.
}
\\
\noindent{\bf Regular Level Set Theorem} \cite[Corollary A.26]{lee2018introduction}  {\it Let
$M$ and $N$ be smooth manifolds, and let $\Phi:M\rightarrow N$ be a smooth map. Every regular level
set of $\Phi$   is a properly embedded submanifold of $M$ whose codimension is equal to
$dim(N)$.}
\\

Because the level curves of $|W_{Y}|$ and $|W_{Y}|^2$ are identical, it suffices to analyze those of the latter.
Since $W_{Y}=W_{f}+W_{\Phi}$ and $W_{\Phi}$
is a random field, for each random element $\omega$ in the probability space $\Omega$,
we denote the corresponding sample path of $W_{\Phi}$ by $W_{\Phi}(\cdot,\cdot;\omega)$.
According to \cite[Proposition 1]{liu2024analyzing}, the assumptions on the Fourier transform of the wavelet $\psi$
ensure the existence of an event $\mathcal{E}\subset \Omega$ with $\mathbb{P}(\mathcal{E})=1$ such that
 $W_{\Phi}(\cdot,\cdot;\omega)$ is twice continuously differentiable on $\mathbb{R}_{+}\times \mathbb{R}$ for all $\omega\in \mathcal{E}$.

On the other hand, since both $f$ and $\psi$ are square-integrable, the Plancherel theorem yields
\begin{align}\label{W_f_Plancherel}
W_{f}(t,s) = \frac{\sqrt{s}}{2\pi}\int_{\mathbb{R}}e^{it\lambda}\widehat{f}(\lambda)\overline{\widehat{\psi}(s\lambda)}d\lambda.
\end{align}
The assumptions on the Fourier transform of $\psi$, together with (\ref{W_f_Plancherel}),
imply that $W_{f}$ is also twice continuously differentiable on $\mathbb{R}\times \mathbb{R}_{+}$.
Hence, $W_{Y}(\cdot,\cdot;\omega)$
is twice continuously differentiable on $\mathbb{R}\times \mathbb{R}_{+}$ for $\omega\in\mathcal{E}$.

In what follows, we fix a $\omega\in\mathcal{E}$.
Because the conditions of Sard's theorem are satisfied by the function $|W_{Y}(\cdot,\cdot;\omega)|^{2}$,
the set of its critical values has Lebesgue measure zero in $\mathbb{R}_{+}$.
That is, there exists a subset $I(\omega)\subset\mathbb{R}_{+}$ such that
$\mathbb{R}_{+}\setminus I(\omega)$ has Lebesgue measure zero and the gradient of
$|W_{Y}(t,s;\omega)|^2$ does not vanish at all points $(t,s)\in K$, where
\begin{align*}
K = \left\{(t',s')\in \mathbb{R}\times \mathbb{R}_{+}\mid |W_{Y}(t',s';\omega)|^2 \in  I(\omega)\right\}.
\end{align*}

For any $(t_{0},s_{0})\in K$, because the gradient of $|W_{Y}(t_0,s_0;\omega)|^2$ does not vanish
and
it consists of two components,
we consider the following two cases:

{\bf Case 1:} $\frac{\partial |W_{Y}|^2}{\partial t}(t_{0},s_{0})\neq 0.$
By the implicit function theorem,  in a neighborhood $U$ of $t_{0}$, there exists a unique differentiable function $\varphi_{1}$
with $\varphi_{1}(t_{0}) = s_{0}$
such that
\begin{align*}
|W_{Y}|^2\left(t,\varphi_{1}(t)\right)  = |W_{Y}|^2\left(t_{0},s_{0}\right)
\end{align*}
for all $t\in U$.

{\bf Case 2:} $\frac{\partial |W_{Y}|^2}{\partial s}(t_{0},s_{0})\neq 0.$
By the implicit function theorem,  in a neighborhood $V$ of $t_{0}$, there exists a unique differentiable function $\varphi_{2}$
with $\varphi_{2}(s_{0}) = t_{0}$
such that
\begin{align*}
|W_{Y}|^2\left(\varphi_{2}(s),s\right)  = |W_{Y}|^2\left(t_{0},s_{0}\right)
\end{align*}
for all $t\in V$.

Regarding the last statement in Proposition \ref{thm:contour:c1}, because $|W_{Y}|^2\left(\cdot,\cdot;\omega\right)$ is a real-valued $C^{2}$ function on $\mathbb{R}\times \mathbb{R}_{+}$, and each $c\in I(\omega)$ is a regular value (i.e., the gradient of $|W_{Y}|^2\left(\cdot,\cdot;\omega\right)$ is nonzero on $L_{c}$),
it follows from the regular level set theorem
that $L_{c}$ is an embedded submanifold of dimension one.
\qed

\subsection{Proof of Proposition \ref{prop:jpdf:Thetan}}\label{sec:proof:prop:jpdf:Thetan}
Equation (\ref{jpdf:ThetaY:n:nonnull}) is obtained by expressing the variable $\mathbf{w}$
in the probability density function~(\ref{pdf_WY}) of $W_{Y}$ in polar coordinates and integrating with respect to the radial components.

For the case $n=1$, from (\ref{jpdf:ThetaY:n:nonnull}), we have
\begin{align}\label{proof:phase:pdf1}
&p(\theta;\Theta_{Y}(t,s))
\\\notag=&\frac{1}{\pi \sigma_{s}^{2}}\int_{0}^{\infty}\textup{exp}
\left\{-\frac{\left(r\cos\theta-\mathfrak{R}(W_{f}(t,s))\right)^{2}
}{\sigma_{s}^{2}}
-\frac{\left(r\sin\theta-\mathfrak{I}(W_{f}(t,s))\right)^{2}}{\sigma_{s}^{2}}
\right\}rdr
\\\notag=&\frac{1}{\pi \sigma_{s}^{2}}\textup{exp}\left\{-\frac{1}{\sigma_{s}^{2}}\left|W_{f}(t,s)\right|^{2}\right\}
\\\notag&\times\int_{0}^{\infty}\textup{exp}\left\{\frac{-r^2+2r\mathfrak{R}(W_{f}(t,s))\cos\theta}{\sigma_{s}^{2}}
+\frac{2r\mathfrak{I}(W_{f}(t,s))\sin\theta}{\sigma_{s}^{2}}\right\}rdr
\end{align}
for $\theta\in [0,2\pi)$, where $\sigma^{2}_{s} = \mathbb{E}[|W_{\Phi}(t,s)|^{2}]$.
By the equation
\begin{equation*}
\frac{\mathfrak{R}(W_{f}(t,s))}{|W_{f}(t,s)|}\cos\theta
 +\frac{\mathfrak{I}(W_{f}(t,s))}{|W_{f}(t,s)|}\sin\theta
 =\cos\left(\theta-\theta_{f}\right),
\end{equation*}
the probability density function of $\Theta_{Y}(t,s)$ in (\ref{proof:phase:pdf1})  can be rewritten as
\begin{align}\label{proof:phase:pdf2}
p(\theta;\Theta_{Y}(t,s))
=&\frac{1}{\pi \sigma_{s}^{2}}\textup{exp}\left\{-\frac{1}{\sigma_{s}^{2}}\left|W_{f}(t,s)\right|^{2}\right\}
\\\notag&\times\int_{0}^{\infty}\textup{exp}\left\{-\frac{1}{\sigma_{s}^{2}}
\left[r^2-2|W_{f}(t,s)|\cos(\theta-\theta_{f})r\right]\right\}rdr.
\end{align}
By completing the square,
\begin{align}\label{integral_formula}
&\int_{0}^{\infty}\hspace{-0.3cm}\textup{exp}\left\{\hspace{-0.1cm}-\frac{1}{\sigma_{s}^{2}}\left[r^2-2|W_{f}(t,s)|\cos(\theta-\theta_{f})r\right]\right\}rdr
 \\\notag=&\frac{\sigma_{s}^{2}}{2}\hspace{-0.05cm}+\hspace{-0.05cm}\frac{\sqrt{q\pi}\sigma_{s}^{2}}{2}
\cos(\theta-\theta_{f}) \textup{exp}\left(q\cos^{2}(\theta-\theta_{f})\right)
\hspace{-0.1cm}B_{q}(\theta-\theta_{f}),
\end{align}
where $q=q(t,s)$ is defined in (\ref{def:SNR}), and $B_{q}$ is defined in (\ref{def:B_q}).
The equation (\ref{margianl_pdf_phase}) follows by substituting (\ref{integral_formula}) into (\ref{proof:phase:pdf2}).
\qed

\subsection{Proof of Corollary \ref{corollary:phase:jpdf:null}}\label{sec:proof:corollary:phase:jpdf:null}

(a) In the null case $Y=\Phi$, (\ref{jpdf:ThetaY:n:nonnull}) simplifies to
\begin{align}\label{jpdf:phase:n}
p\left(\theta_{1:n};\Theta_{\Phi}\right)
=&
\frac{1}{\pi^{n}\textup{det}(\Gamma)}
\\\notag& \times\int_{0}^{\infty}\hspace{-0.1cm}\int_{0}^{\infty}\cdots\int_{0}^{\infty}\hspace{-0.1cm}\left(\overset{n}{\underset{\ell=1}{\prod}}r_{\ell}\right)
 \exp
\left(-\mathbf{r}^{\top}\mathbf{M}\mathbf{r}\right)dr_{1}dr_{2}\cdots dr_{n},
\end{align}
where $\mathbf{r} = [ r_{1}\  r_{2} \ \cdots\  r_{n}]^{\top}$
and
\begin{align*}
\mathbf{M} = \textup{diag}\left(e^{-i\theta_1}, \dots, e^{-i\theta_n} \right)\Gamma^{-1} \textup{diag}\left(e^{i\theta_1}, \dots, e^{i\theta_n}\right).
\end{align*}
Let's express $\mathbf{r}$ in spherical coordinates by making the change of variables
$\mathbf{r}= \rho \omega$, where $\rho=|\mathbf{r}|$ is the radial part and
$\omega = [\omega_{1},\ldots,\omega_{n}]^{\top}$ is a point
on the positive orthant of the unit sphere.
Under this transformation, each component becomes $r_{\ell}=\rho \omega_{\ell}$ for $\ell\in\{1,2,\ldots,n\}$.
The integral in (\ref{jpdf:phase:n}) then becomes
\begin{align}\notag
p\left(\theta_{1:n};\Theta_{\Phi}\right)
=
\frac{1}{\pi^{n}\textup{det}(\Gamma)}
 \int_{0}^{\infty}\int_{S_{+}^{n-1}}
 \rho^{n}\left(\overset{n}{\underset{\ell=1}{\prod}}\omega_{\ell}\right)
 \exp
\left(-\rho^{2}\omega^{\top}\mathbf{M}\omega\right)\rho^{n-1}d\omega d\rho,
\end{align}
where $d\omega$ is the surface measure on $S_{+}^{n-1}$.
Because $\mathbf{M}$ is positive definite, we have $\omega^{\top}\mathbf{M}\omega>0$ and
\begin{align}\label{jpdf:phase:n:v2}
p\left(\theta_{1:n};\Theta_{\Phi}\right)
=
\frac{(n-1)!}{2\pi^{n}\textup{det}(\Gamma) }
\int_{S_{+}^{n-1}}
\left(\overset{n}{\underset{\ell=1}{\prod}}\omega_{\ell}\right)
\left(\omega^{\top}\mathbf{M}\omega\right)^{-n}d\omega.
\end{align}

(b) For the case $n=2$, applying the change of variables
$\omega_{1}=\cos u$ and $\omega_{2}=\sin u$, where $u\in[0,\pi/2)$,
(\ref{jpdf:phase:n:v2}) can be rewritten as
\begin{align}\notag
2\pi^{2}\textup{det}(\Gamma) p\left(\theta_{1:2};\Theta_{\Phi}\right)
=&
\int_{0}^{\pi/2}\hspace{-0.2cm}
\frac{\sin u\cos u }{\left(M_{22}\sin^{2}u+2 \Re(M_{12})\cos u \sin u+M_{11}\cos^{2}u \right)^{2}}
du
\\\label{jpdf:phase:n:v3}=&
\int_{0}^{\pi/2}
\frac{\sin u\cos^{-3}u }{\left(M_{22}\tan^{2}u+ 2 \Re(M_{12}) \tan u+M_{11} \right)^{2}}
du,
\end{align}
where $\mathbf{M}=[M_{\ell\ell'}]_{1\leq \ell,\ell'\leq 2}$.
By substituting $v = \tan u$ into (\ref{jpdf:phase:n:v3}) and noting that
$du =(1+v^2)^{-1}dv$ and $\sin u \cos^{-3}u=v(1+v^2)$, we have
\begin{align}\notag
2\pi^{2}\textup{det}(\Gamma) p\left(\theta_{1:2};\Theta_{\Phi}\right)
=&
\int_{0}^{\infty}
\frac{v}{\left( M_{22}v^2 +2 \Re(M_{12}) v+M_{11}\right)^{2}}
dv
\\\notag=&\frac{1}{2(M_{11}M_{22}-|\Re(M_{12})|^2)}
-\frac{\Re(M_{12})}{2(M_{11}M_{22}-|\Re(M_{12})|^2)^{3/2}}
\\\label{jpdf:phase:n:v4}&\ \ \ \times\Big[\frac{\pi}{2}-
\tan^{-1}\Big(\frac{\Re(M_{12})}{\sqrt{M_{11}M_{22}-|\Re(M_{12})|^2}}
\Big)\Big].
\end{align}
By substituting
\begin{align*}
\mathbf{M}=\frac{1}{\textup{det}(\Gamma)}\begin{bmatrix} \Gamma_{22} & -\Gamma_{12}e^{i(\theta_{2}-\theta_{1})} \\
-\overline{\Gamma_{12}}e^{i(\theta_{1}-\theta_{2})} & \Gamma_{11}\end{bmatrix}
\end{align*}
into (\ref{jpdf:phase:n:v4}),
\begin{align*}
p\left(\theta_{1:2};\Theta_{\Phi}\right)
=&\frac{1 }{4\pi^{2}}\frac{\textup{det}(\Gamma)}{\Gamma_{11}\Gamma_{22}-|\Re(\Gamma_{12}e^{i(\theta_{2}-\theta_{1})})|^2}
\\\notag&+\frac{1 }{4\pi^{2}}\frac{\textup{det}(\Gamma)\Re(\Gamma_{12}e^{i(\theta_{2}-\theta_{1})})}
{(\Gamma_{11}\Gamma_{22}-|\Re(\Gamma_{12}e^{i(\theta_{2}-\theta_{1})})|^2)^{3/2}}
\\&\ \ \ \times\Big[\frac{\pi}{2}+
\tan^{-1}\Big(\frac{\Re(\Gamma_{12}e^{i(\theta_{2}-\theta_{1})})}{\sqrt{\Gamma_{11}\Gamma_{22}-|\Re(\Gamma_{12}e^{i(\theta_{2}-\theta_{1})})|^2}}
\Big)\Big].
\end{align*}
\qed

\subsection{Proof of Proposition \ref{cov:mag:generalpsi}}\label{sec:proof:cov:mag:generalpsi}
First of all, for any $t\in \mathbb{R}$ and $s>0$,
\begin{align}\notag
|W_{Y}(t,s)|^2
=& \left(\mathfrak{R}(W_{f}(t,s))\right)^{2}+ 2\mathfrak{R}(W_{f}(t,s)) \mathfrak{R}(W_{\Phi}(t,s))
\\\notag&
+\left( \mathfrak{R}(W_{\Phi}(t,s))\right)^{2}
+\left(\mathfrak{I}(W_{f}(t,s))\right)^{2}
\\\notag&+ 2\mathfrak{I}(W_{f}(t,s)) \mathfrak{I}(W_{\Phi}(t,s))+\left( \mathfrak{I}(W_{\Phi}(t,s))\right)^{2}.
\end{align}
By (\ref{spectral_rep:samplepath}) and It$\hat{\textup{o}}$'s formula,
\begin{align*}
\Box(W_{\Phi}(t,s))
= \sqrt{s}\int_{\mathbb{R}}e^{it\lambda} \overline{\widehat{\Box(\psi)}}(s\lambda)Z_{F}(d\lambda)
\end{align*}
and
\begin{align}\notag
\left(\Box(W_{\Phi}(t,s))\right)^2 = s\int_{\mathbb{R}}|\widehat{\Box(\psi)}(su)|^2
+s\int_{\mathbb{R}^{2}}'e^{it(u+v)} \overline{\widehat{\Box(\psi)}}(su)
\overline{\widehat{\Box(\psi)}}(sv)Z_{F}(du)Z_{F}(dv)
F(du),
\end{align}
where $\Box\in \{\Re,\Im\}.$
Hence,
\begin{align}\notag
|W_{Y}(t,s)|^2-\mathbb{E}\left[|W_{Y}(t,s)|^2\right]
=&
 2\Re(W_{f}(t,s))\sqrt{s}\int_{\mathbb{R}}e^{it\lambda} \overline{\widehat{\Re(\psi)}}(s\lambda)Z_{F}(d\lambda)
\\\notag&+ 2\Im(W_{f}(t,s))\sqrt{s}\int_{\mathbb{R}}e^{it\lambda} \overline{\widehat{\Im(\psi)}}(s\lambda)Z_{F}(d\lambda)
\\\notag&+s\int_{\mathbb{R}^{2}}'e^{it(u+v)} \overline{\widehat{\Re(\psi)}}(su)
\overline{\widehat{\Re(\psi)}}(sv)Z_{F}(du)Z_{F}(dv)
\\\notag&+s\int_{\mathbb{R}^{2}}'e^{it(u+v)} \overline{\widehat{\Im(\psi)}}(su)
\overline{\widehat{\Im(\psi)}}(sv)Z_{F}(du)Z_{F}(dv).
\end{align}
By the product formula \cite{major1981lecture,nourdin2012normal}, for any $t_{1},t_{2}\in \mathbb{R}$ and $s_{1},s_{2}>0$,
\begin{align}\notag
&\textup{Cov}\left(|W_{Y}(t_{1},s_{1})|^2, |W_{Y}(t_{2},s_{2})|^2\right)
\\\notag=&4\underset{\Box,\triangle \in \{\Re,\Im\}}{\sum} \Box(W_{f}(t_{1},s_{1}))\triangle(W_{f}(t_{2},s_{2}))
\sqrt{s_{1}s_{2}}\int_{\mathbb{R}}e^{i(t_{1}-t_{2})\lambda} \overline{\widehat{\Box(\psi)}}(s_{1}\lambda)\widehat{\triangle(\psi)}(s_{2}\lambda)F(d\lambda)
\\\notag&+2s_{1}s_{2}\left[\int_{\mathbb{R}}e^{i(t_{1}-t_{2})\lambda} \overline{\widehat{\Box(\psi)}}(s_{1}u)
\widehat{\triangle(\psi)}(s_{2}u)F(du)\right]^{2}.
\end{align}
\qed

\subsection{Proof of Proposition \ref{prop:Cov_general}}\label{sec:proof:prop:Cov_general}
The equation (\ref{format:cov:mag}) follows directly from the definition of covariance.
Hence, it suffices to derive (\ref{prop:pdf:product}).
Let $g: \mathbb{C}\setminus\{0\}\times \mathbb{C}\mapsto \mathbb{C}\setminus\{0\}\times \mathbb{C}$ be defined by $g(w_{1},w_{2}) = (w_{1},w_{1}w_{2})$, which has the inverse
$$g^{-1}(w_{1},w_{2}) =\left(g^{-1}_{1}(w_{1},w_{2}),g^{-1}_{2}(w_{1},w_{2})\right) =  (w_{1},w_{2}/w_{1}).$$
The complex Jacobian matrix of $g$ is a $4\times4$ matrix given by
\begin{align*}
\mathbf{J}_{g}:=&\left[\begin{array}{cccc}
\partial g^{-1}_{1}/\partial w_{1} & \partial g^{-1}_{1}/\partial w_{2} & \partial g^{-1}_{1}/\partial \overline{w_{1}} & \partial g^{-1}_{1}/\partial \overline{w_{2}}
\\
\partial g^{-1}_{2}/\partial w_{1} & \partial g^{-1}_{2}/\partial w_{2} & \partial g^{-1}_{2}/\partial \overline{w_{1}} & \partial g^{-1}_{2}/\partial \overline{w_{2}}
\\
\partial \overline{g^{-1}_{1}}/\partial w_{1} & \partial \overline{g^{-1}_{1}}/\partial w_{2} & \partial \overline{g^{-1}_{1}}/\partial \overline{w_{1}} & \partial \overline{g^{-1}_{1}}/\partial \overline{w_{2}}
\\
\partial \overline{g^{-1}_{2}}/\partial w_{1} & \partial \overline{g^{-1}_{2}}/\partial w_{2} & \partial \overline{g^{-1}_{2}}/\partial \overline{w_{1}} & \partial \overline{g^{-1}_{2}}/\partial \overline{w_{2}}
\end{array}\right]
\\=&\left[\begin{array}{cccc}
1 & 0 & 0 & 0
\\
-w_{2}/w_{1}^{2} & 1/w_{1} & 0 & 0
\\
0 & 0 & 1 & 0
\\
0 & 0 & -\overline{w_{2}}/\overline{w_{1}}^{2} & 1/\overline{w_{1}}
\end{array}\right].
\end{align*}
It has the determinant  $1/|w_{1}|^{2}$.
By (\ref{pdf_WY}) with $\mathbf{W}_{f} = [W_{f}(t_{1},s_{1})\ W_{f}(t_{2},s_{2})]^{\top}$
and $\mathbf{W}_{Y} = [W_{Y}(t_{1},s_{1})\ W_{Y}(t_{2},s_{2})]^{\top}$,
\begin{align}\notag
&p\left(w_{1},z;W_{Y}(t_{1},s_{1}),W_{Y}(t_{1},s_{1})W_{Y}(t_{2},s_{2})\right)
\\\notag=& |w_{1}|^{-2} p\left(g^{-1}(w_{1},z);\mathbf{W}_{Y}\right)
\\\notag=&
\frac{1}{\pi^{2}\textup{det}(\Gamma)|w_{1}|^{2}}
\textup{exp}\left(-\left(
\left[\begin{array}{cc}w_{1}\\ z/w_{1}\end{array}\right]-\mathbf{W}_{f}\right)^{*}
\Gamma^{-1}
\left(\left[\begin{array}{cc}w_{1}\\ z/w_{1}\end{array}\right]-\mathbf{W}_{f}\right)\right)
\end{align}
for any $w_{1},z\in \mathbb{C}$, 
which implies that the marginal probability density function of $W_{Y}(t_{1},s_{1})W_{Y}(t_{2},s_{2})$ has the integral form
\begin{align}\notag
&p\left(z;W_{Y}(t_{1},s_{1})W_{Y}(t_{2},s_{2})\right)
\\\notag=&\int_{0}^{\infty}\int_{0}^{2\pi}
p\left(re^{i\theta},z;W_{Y}(t_{1},s_{1}),W_{Y}(t_{1},s_{1})W_{Y}(t_{2},s_{2})\right)
r d\theta dr
\\\notag=&
\int_{0}^{\infty}\int_{0}^{2\pi} \frac{1}{\pi^{2}\textup{det}(\Gamma)r}
\textup{exp}\left(-\left(
\left[\begin{array}{cc}re^{i\theta}\\ zr^{-1}e^{-i\theta}\end{array}\right]-\mathbf{W}_{f}\right)^{*}
\Gamma^{-1}
\left(\left[\begin{array}{cc}re^{i\theta}\\ zr^{-1}e^{-i\theta}\end{array}\right]-\mathbf{W}_{f}\right)\right)d\theta dr.
\end{align}
\qed

\subsection{Proof of Corollary \ref{prop:Cov_null}}\label{sec:proof:prop:Cov_null}
For the null case, i.e., $f=0$, (\ref{prop:pdf:product}) can be simplified as follows
\begin{align}\label{null_product}
&p\left(z;W_{\Phi}(t_{1},s_{1})W_{\Phi}(t_{2},s_{2})\right)
\\\notag=&
\int_{0}^{\infty}\frac{1}{\pi^{2}\textup{det}(\Gamma)r}\left\{\int_{0}^{2\pi}
\textup{exp}\left(-
\left[\begin{array}{cc}re^{i\theta}\\ zr^{-1}e^{-i\theta}\end{array}\right]^{*}
\Gamma^{-1}
\left[\begin{array}{cc}re^{i\theta}\\ zr^{-1}e^{-i\theta}\end{array}\right]\right)d\theta\right\} dr,
\end{align}
where $z\in \mathbb{C}$,
\begin{align*}
\Gamma^{-1}=
\frac{1}{\det(\Gamma)}\left[\begin{array}{rr}\Gamma_{22} & -\Gamma_{12}\\ -\overline{\Gamma_{12}} & \Gamma_{11}\end{array}\right],
\end{align*}
and
$\det(\Gamma)=\Gamma_{11}\Gamma_{22}-|\Gamma_{12}|^{2}>0$.
Note that
\begin{align}\label{complex_exponent}
\left[\begin{array}{cc}re^{i\theta}\\ zr^{-1}e^{-i\theta}\end{array}\right]^{*}
\Gamma^{-1}
\left[\begin{array}{cc}re^{i\theta}\\ zr^{-1}e^{-i\theta}\end{array}\right]
=\frac{\Gamma_{22} r^{2}+\Gamma_{11}|z|^{2}r^{-2}-2\mathfrak{R}\left(z\Gamma_{12}e^{-i2\theta}\right)}
{\det(\Gamma)}.
\end{align}
By substituting (\ref{complex_exponent}) into (\ref{null_product}),
\begin{align}\label{null_product2}
p\left(z;W_{\Phi}(t_{1},s_{1})W_{\Phi}(t_{2},s_{2})\right)
=&
\left\{\int_{0}^{2\pi}
\textup{exp}
\left(
\frac{2\mathfrak{R}\left(z\Gamma_{12} e^{-i2\theta}
\right)}{\det(\Gamma)}\right)
d\theta\right\}
\\\notag&\times\left\{\int_{0}^{\infty}\frac{1}{\pi^{2}\textup{det}(\Gamma)r}
\textup{exp}\left(-
\frac{\Gamma_{22} r^{2}+\Gamma_{11}|z|^{2}r^{-2}}
{\det(\Gamma)}
\right)dr\right\}.
\end{align}
Let $z\Gamma_{12}= |z\Gamma_{12}|e^{i\nu}$ for some $\nu\in [0,2\pi)$.
For the first integral in (\ref{null_product2}), we have
\begin{align}\notag
\int_{0}^{2\pi}
\textup{exp}
\left(
\frac{2\mathfrak{R}\left(z\Gamma_{12}e^{-i2\theta}
\right)}{\det(\Gamma)}\right)
d\theta
=&\int_{0}^{2\pi}
\textup{exp}
\left(
\frac{2|z\Gamma_{12}|}{\det(\Gamma)}\mathfrak{R}\left(e^{i(\nu-2\theta)}
\right)\right)
d\theta
\\\notag=&
\int_{0}^{2\pi}
\textup{exp}
\left(
-\frac{2|z\Gamma_{12}|}{\det(\Gamma)}\mathfrak{R}\left(e^{i(2\theta-\nu)}
\right)\right)
d\theta
\\\notag=&\frac{1}{2}
\int_{-\nu}^{4\pi-\nu}
\textup{exp}
\left(
-\frac{2|z\Gamma_{12}|}{\det(\Gamma)}\mathfrak{R}\left(e^{iu}
\right)\right)
du
\\\notag=&\int_{0}^{2\pi}
\textup{exp}
\left(
-\frac{2|z\Gamma_{12}|}{\det(\Gamma)}\mathfrak{R}\left(e^{iu}
\right)\right)
du
\\\label{I0_app}=&2\pi I_{0}\left(\frac{2|\Gamma_{12}|}{\det(\Gamma)}|z|\right),
\end{align}
where
\begin{align*}
I_{0}(x)= \frac{1}{2\pi}\int_{0}^{2\pi}e^{x\cos(u)}du,\ x\in \mathbb{R},
\end{align*}
is the modified Bessel function of the first kind of order zero.
For the second integral in (\ref{null_product2}),
by the formula
\begin{align*}
\int_{0}^{\infty}r^{-1}e^{-ar^{2}-br^{-2}}dr=K_{0}\left(2\sqrt{ab}\right),\ a,b>0,
\end{align*}
where $K_{0}$ is the modified Bessel function of the second kind of order zero,
we have
\begin{align}\label{K0_app}
\int_{0}^{\infty}\frac{1}{\pi^{2}\textup{det}(\Gamma)r}
\textup{exp}\left(-
\frac{\Gamma_{22} r^{2}+\Gamma_{11}|z|^{2}r^{-2}}
{\det(\Gamma)}
\right)dr
=\frac{1}{\pi^{2}\textup{det}(\Gamma)}K_{0}\left(\frac{2\sqrt{\Gamma_{11}\Gamma_{22}}}{\det(\Gamma)}|z|\right).
\end{align}
By substituting (\ref{I0_app}) and  (\ref{K0_app}) into (\ref{null_product2}),
\begin{align}\label{null_product3}
p\left(z; W_{\Phi}(t_{1},s_{1})W_{\Phi}(t_{2},s_{2})\right)
=
\frac{2}{\pi\textup{det}(\Gamma)}I_{0}\left(\frac{2|\Gamma_{12}|}{\det(\Gamma)}|z|\right)
K_{0}\left(\frac{2\sqrt{\Gamma_{11}\Gamma_{22}}}{\det(\Gamma)}|z|\right)
\end{align}
for any $z\in \mathbb{C}$.
The equation (\ref{null_product3}) shows that the probability density function of $W_{\Phi}(t_{1},s_{1})W_{\Phi}(t_{2},s_{2})$
at any $z\in \mathbb{C}$
depends on $|z|$.
By (\ref{null_product3}) and the formula \cite[p. 684]{gradshteyn2014table}
\begin{align*}
\int_{0}^{\infty}r^{2}I_{0}(ar)K_{0}(br)dr = \frac{\pi}{2}\frac{1}{b^{3}}{}_2F_1\left(\frac{3}{2},\frac{3}{2};1;\frac{a^2}{b^2}\right)
\end{align*}
for any $a\geq 0$ and $b>a$,we have
\begin{align}\notag
\mathbb{E}\left[\left|W_{\Phi}(t_{1},s_{1})\right|\left|W_{\Phi}(t_{2},s_{2})\right|\right]
=&\int_{0}^{2\pi}\int_{0}^{\infty}r^{2}p\left(r; W_{\Phi}(t_{1},s_{1})W_{\Phi}(t_{2},s_{2})\right)drd\theta
\\\notag=&\frac{4\pi}{\pi\textup{det}(\Gamma)}\left[\frac{\pi}{2}\left(\frac{\textup{det}(\Gamma)}
{2\sqrt{\Gamma_{11}\Gamma_{22}}}\right)^{3}
{}_2F_1\left(\frac{3}{2},\frac{3}{2};1;\frac{|\Gamma_{12}|^{2}}{\Gamma_{11}\Gamma_{22}}\right)
\right]
\\\notag=&
\frac{\pi}{4}\left(\textup{det}(\Gamma)\right)^{2}
(\Gamma_{11}\Gamma_{22})^{-\frac{3}{2}}
{}_2F_1\left(\frac{3}{2},\frac{3}{2};1;\frac{|\Gamma_{12}|^{2}}{\Gamma_{11}\Gamma_{22}}\right)
\\\notag=&
\frac{\pi}{4}\left(\Gamma_{11}\Gamma_{22}\right)^{\frac{1}{2}}
\left(1-\frac{|\Gamma_{12}|^{2}}{\Gamma_{11}\Gamma_{22}}\right)^{2}
{}_2F_1\left(\frac{3}{2},\frac{3}{2};1;\frac{|\Gamma_{12}|^{2}}{\Gamma_{11}\Gamma_{22}}\right).
\end{align}
On the other hand,
\begin{align}\label{EabsW}
\mathbb{E}\left[\left|W_{\Phi}(t_{\ell},s_{\ell})\right|\right]=\frac{1}{2}(\pi\Gamma_{\ell\ell})^{\frac{1}{2}}
\end{align}
for  $\ell\in\{1,2\}$. Therefore,
\begin{align}\label{Cov_abs_W}
&\textup{Cov}\left(\left|W_{\Phi}(t_{1},s_{1})\right|,\left|W_{\Phi}(t_{2},s_{2})\right|\right)
\\\notag=&\frac{\pi}{4}\sqrt{\Gamma_{11}\Gamma_{22}}
\left[\left(1-\frac{|\Gamma_{12}|^{2}}{\Gamma_{11}\Gamma_{22}}\right)^2
{}_2F_1\left(\frac{3}{2},\frac{3}{2};1;\frac{|\Gamma_{12}|^{2}}{\Gamma_{11}\Gamma_{22}}\right)-1\right].
\end{align}
On the other hand, by the definition $\Gamma_{\ell\ell}=[\left|W_{\Phi}(t_{\ell},s_{\ell})\right|^{2}]$ and
 (\ref{EabsW}), we have
\begin{align}\label{Var_abs_W}
\textup{Var}\left(\left|W_{\Phi}(t_{\ell},s_{\ell})\right|\right)
=\left(1-\frac{\pi}{4}\right)\Gamma_{\ell\ell},\
\ell=1,2.
\end{align}
The equation (\ref{eq:correlation_mag}) follows by combining (\ref{Cov_abs_W}) and (\ref{Var_abs_W}).
The asymptotic behavior of the correlation coefficient given in (\ref{implication:corollary:mag}) follows from the property that
for any $a,b,c>0$,
\begin{align*}
\underset{x\rightarrow0+}{\lim} \frac{c}{abx}  \left({}_2F_1\left(a,b;c;x\right)-1\right)=1.
\end{align*}
\qed

\subsection{Proof of Corollary \ref{prop:Cov2_phase_null}}\label{sec:proof:prop:Cov2_phase_null}
By Proposition \ref{prop:jpdf:Thetan}(a),
\begin{align}\notag
&\textup{Cov}\left(\Theta_{\Phi}(t_{1},s_{1}),\Theta_{\Phi}(t_{2},s_{2})\right)
\\\notag=&\int_{0}^{2\pi}\int_{0}^{2\pi}\int_{0}^{\infty}\int_{0}^{\infty} \left(\theta_{1}-\pi\right)\left(\theta_{2}-\pi\right)\frac{r_{1}r_{2}}{\pi^{2}\textup{det}(\Gamma)}
\\\notag&\times\textup{exp}
\left(-\begin{bmatrix}r_{1}e^{i\theta_{1}} \\ r_{2}e^{i\theta_{2}}\end{bmatrix}^{*}\Gamma^{-1}\begin{bmatrix}r_{1}e^{i\theta_{1}} \\ r_{2}e^{i\theta_{2}}\end{bmatrix}\right)
dr_{1}dr_{2}d\theta_{1}d\theta_{2}
\\\notag=&
\int_{0}^{2\pi}\int_{0}^{2\pi}\int_{0}^{\infty}\int_{0}^{\infty} \left(\theta_{1}-\pi\right)\left(\theta_{2}-\pi\right) \frac{r_{1}r_{2}}{\pi^{2}\textup{det}(\Gamma)}
\\\notag&\times\textup{exp}\left(-\frac{\Gamma_{11}r_{1}^{2}+\Gamma_{22}r_{2}^{2}-2\Re(e^{i(\theta_{1}-\theta_{2})}
\overline{\Gamma_{12}})r_{1}r_{2}}{\textup{det}(\Gamma)}\right)
 dr_{1}dr_{2}d\theta_{1}d\theta_{2}
\\\notag=&
\frac{\textup{det}(\Gamma)}{\Gamma_{11}\Gamma_{22}}\int_{0}^{2\pi}\int_{0}^{2\pi}
\int_{0}^{\infty}\int_{0}^{\infty}\left(\theta_{1}-\pi\right)\left(\theta_{2}-\pi\right) \frac{r_{1}r_{2}}{\pi^{2}}
\\\label{cov:theta_theta1_pre}&\textup{exp}\left(-r_{1}^{2}-r_{2}^{2}+\frac{2|\Gamma_{12}|}{\sqrt{\Gamma_{11}\Gamma_{22}}}
\cos\Big(\theta_{1}-\theta_{2}-\arg(\Gamma_{12})\Big)r_{1}r_{2}
\right)
 dr_{1}dr_{2}d\theta_{1}d\theta_{2}.
\end{align}
From (\ref{cov:theta_theta1_pre}),
\begin{align}\label{cov:theta_theta1}
\textup{Cov}\left(\Theta_{\Phi}(t_{1},s_{1}),\Theta_{\Phi}(t_{2},s_{2})\right)
=h(\frac{|\Gamma_{12}|}{\sqrt{\Gamma_{11}\Gamma_{22}}},\arg(\Gamma_{12})),
\end{align}
where $h(u,v)=\left(1-u^2\right)\widetilde{h}(u,v)$ and
\begin{align}\notag
\widetilde{h}(u,v) =&
\int_{0}^{2\pi}\int_{0}^{2\pi}\int_{0}^{\infty}\int_{0}^{\infty} \left(\theta_{1}-\pi\right)\left(\theta_{2}-\pi\right)\frac{r_{1}r_{2}}{\pi^{2}}
\\\notag&\textup{exp}\left(-r_{1}^{2}-r_{2}^{2}
\right)\textup{exp}\Big(
2u
\cos\big(\theta_{1}-\theta_{2}-v\big)r_{1}r_{2}
\Big)
dr_{1}dr_{2}d\theta_{1}d\theta_{2}
\end{align}
for $u\geq 0$ and $v\in[0,2\pi)$.
By the fact $\lim_{u\rightarrow0}h(u,v)=0$ and $\widetilde{h}(0,v)=0$ for any $v\in[0,2\pi)$,
\begin{align}\notag
\underset{u\rightarrow0+}{\lim} \frac{h(u,v)}{u}
=&
\underset{u\rightarrow0+}{\lim}
-2u\widetilde{h}(u,v)
\\\notag&+(1-u^2)\int_{0}^{2\pi}\int_{0}^{2\pi}\int_{0}^{\infty}\int_{0}^{\infty} \left(\theta_{1}-\pi\right)\left(\theta_{2}-\pi\right)\frac{r_{1}r_{2}}{\pi^{2}}
\\\notag&\times\textup{exp}\left(-r_{1}^{2}-r_{2}^{2}
\right)\textup{exp}\Big(
2u
\cos\big(\theta_{1}-\theta_{2}-v\big)r_{1}r_{2}
\Big)
\\\notag&\times 2
\cos\big(\theta_{1}-\theta_{2}-v\big)r_{1}r_{2}dr_{1}dr_{2}d\theta_{1}d\theta_{2}
\\\notag
=&\frac{2}{\pi^{2}}\left(\int_{0}^{\infty}r^2 e^{-r^{2}}
dr
\right)^2
\int_{0}^{2\pi}\int_{0}^{2\pi}\left(\theta_{1}-\pi\right)\left(\theta_{2}-\pi\right)
\cos\big(\theta_{1}-\theta_{2}-v\big) d\theta_{1}d\theta_{2}
\\\label{cov:theta_theta2a}=&
\frac{\pi}{2}\cos(v).
\end{align}
By applying (\ref{cov:theta_theta2a}) to (\ref{cov:theta_theta1}),
\begin{align}\notag
\textup{Cov}\left(\Theta_{\Phi}(t_{1},s_{1}),\Theta_{\Phi}(t_{2},s_{2})\right)\sim \frac{\pi}{2}\cos\big(\arg(\Gamma_{12})\big) \frac{|\Gamma_{12}|}{\sqrt{\Gamma_{11}\Gamma_{22}}}
\end{align}
as the ratio  $|\Gamma_{12}|/\sqrt{\Gamma_{11}\Gamma_{22}}$ tends to zero.
\qed


\subsection{Proof of Corollary \ref{prop:Cov_phase_null}}\label{sec:proof:prop:Cov_phase_null}
Because
\begin{align}\label{start:proof:circov}
\mathbb{E}\left[e^{i(\Theta_{\Phi}(t_{1},s_{1})-\Theta_{\Phi}(t_{2},s_{2}))}\right]
=\int_{0}^{\infty}\int_{0}^{2\pi}\hspace{-0.25cm}e^{i\theta}p\left(re^{i\theta}; W_{\Phi}(t_{1},s_{1})\overline{W_{\Phi}(t_{2},s_{2})}\right)
rd\theta dr,
\end{align}
we begin by deriving the probability density function of
$W_{\Phi}(t_{1},s_{1})\overline{W_{\Phi}(t_{2},s_{2})}$.
Let $h: \mathbb{C}\setminus\{0\}\times \mathbb{C}\mapsto \mathbb{C}\setminus\{0\}\times \mathbb{C}$ be defined by $h(w_{1},w_{2}) = (w_{1},w_{1}\overline{w_{2}})$, which has the inverse
$$h^{-1}(w_{1},w_{2}) :=\left(h^{-1}_{1}(w_{1},w_{2}),h^{-1}_{2}(w_{1},w_{2})\right) =  (w_{1},\overline{w_{2}}/\overline{w_{1}}).$$
The complex Jacobian matrix of $h$ is given by
\begin{align*}
\mathbf{J}_{h}:=&\left[\begin{array}{cccc}
\partial h^{-1}_{1}/\partial w_{1} & \partial h^{-1}_{1}/\partial w_{2} & \partial h^{-1}_{1}/\partial \overline{w_{1}} & \partial h^{-1}_{1}/\partial \overline{w_{2}}
\\
\partial h^{-1}_{2}/\partial w_{1} & \partial h^{-1}_{2}/\partial w_{2} & \partial h^{-1}_{2}/\partial \overline{w_{1}} & \partial h^{-1}_{2}/\partial \overline{w_{2}}
\\
\partial \overline{h^{-1}_{1}}/\partial w_{1} & \partial \overline{h^{-1}_{1}}/\partial w_{2} & \partial \overline{h^{-1}_{1}}/\partial \overline{w_{1}} & \partial \overline{h^{-1}_{1}}/\partial \overline{w_{2}}
\\
\partial \overline{h^{-1}_{2}}/\partial w_{1} & \partial \overline{h^{-1}_{2}}/\partial w_{2} & \partial \overline{h^{-1}_{2}}/\partial \overline{w_{1}} & \partial \overline{h^{-1}_{2}}/\partial \overline{w_{2}}
\end{array}\right]
\\=&\left[\begin{array}{cccc}
1 & 0 & 0 & 0
\\
0 & 0 & -\overline{w_{2}}/\overline{w_{1}}^{2} & 1/\overline{w_{1}}
\\
0 & 0 & 1 & 0
\\
-w_{2}/w_{1}^{2} & 1/w_{1} & 0  & 0
\end{array}\right],
\end{align*}
whose determinant is $1/|w_{1}|^{2}$.
By (\ref{pdf_WY}) with $f=0$
and $\mathbf{W}_{\Phi} = [W_{\Phi}(t_{1},s_{1})\ W_{\Phi}(t_{2},s_{2})]^{\top}$,
\begin{align}\notag
&p\left(w_{1},z;W_{\Phi}(t_{1},s_{1}),W_{\Phi}(t_{1},s_{1})\overline{W_{\Phi}(t_{2},s_{2})}\right)
\\\notag=& |w_{1}|^{-2} p\left(h^{-1}(w_{1},z);\mathbf{W}_{\Phi}\right)
\\\notag=&
\frac{1}{\pi^{2}\textup{det}(\Gamma)|w_{1}|^{2}}
\textup{exp}\left(-
\left[\begin{array}{cc}w_{1}\\ \overline{z}/\overline{w_{1}}\end{array}\right]^{*}
\Gamma^{-1}
\left[\begin{array}{cc}w_{1}\\ \overline{z}/\overline{w_{1}}\end{array}\right]\right),\ w_{1},q\in \mathbb{C},
\end{align}
which implies that the marginal density function of $W_{\Phi}(t_{1},s_{1})\overline{W_{\Phi}(t_{2},s_{2})}$ has the integral form
\begin{align}\label{null_product:phase}
&p\left(z;W_{\Phi}(t_{1},s_{1})\overline{W_{\Phi}(t_{2},s_{2})}\right)
\\\notag=&\int_{0}^{\infty}\int_{0}^{2\pi}
p\left(re^{i\theta},z; W_{\Phi}(t_{1},s_{1}),W_{\Phi}(t_{1},s_{1})\overline{W_{\Phi}(t_{2},s_{2})}\right)r d\theta dr
\\\notag=&
\int_{0}^{\infty}\hspace{-0.2cm}\int_{0}^{2\pi} \frac{1}{\pi^{2}\textup{det}(\Gamma)r}
\textup{exp}\left(-
\left[\hspace{-0.15cm}\begin{array}{cc}re^{i\theta}\\ \overline{z}r^{-1}e^{i\theta}\end{array}\hspace{-0.15cm}\right]^{*}
\Gamma^{-1}
\left[\hspace{-0.15cm}\begin{array}{cc}re^{i\theta}\\ \overline{z}r^{-1}e^{i\theta}\end{array}\hspace{-0.15cm}\right]\right)d\theta dr
\\\notag=&
\int_{0}^{\infty}\frac{2}{\pi\textup{det}(\Gamma)r}
\textup{exp}\left(-
\left[\begin{array}{cc}r\\ \overline{z}r^{-1}\end{array}\right]^{*}
\Gamma^{-1}
\left[\begin{array}{cc}r\\ \overline{z}r^{-1}\end{array}\right]\right)dr.
\end{align}
Because
\begin{align*}
\Gamma=
\left[\begin{array}{cc}\Gamma_{11} & \Gamma_{12}\\ \overline{\Gamma_{12}} & \Gamma_{22}\end{array}\right]
\ \textup{and}\
\Gamma^{-1}=
\frac{1}{\det(\Gamma)}\left[\begin{array}{rr}\Gamma_{22} & -\Gamma_{12}\\ -\overline{\Gamma_{12}} & \Gamma_{11}\end{array}\right],
\end{align*}
where $\det(\Gamma_{Y})=\Gamma_{11}  \Gamma_{22}-|\Gamma_{12}|^{2}>0$, we have
\begin{align}\label{complex_exponent:phase}
\left[\begin{array}{cc}r\\ \overline{z}r^{-1}\end{array}\right]^{*}
\Gamma^{-1}
\left[\begin{array}{cc}r\\ \overline{z}r^{-1}\end{array}\right]
=\frac{\Gamma_{22}r^{2}+\Gamma_{11}|z|^{2}r^{-2}-2\mathfrak{R}\left(\overline{z}\Gamma_{12}\right)}
{\det(\Gamma)}.
\end{align}
By substituting (\ref{complex_exponent:phase}) into (\ref{null_product:phase}),
\begin{align}\label{null_product2s}
&p\left(z;W_{\Phi}(t_{1},s_{1})\overline{W_{\Phi}(t_{2},s_{2})}\right)
\\\notag=&
\textup{exp}
\left(
\frac{2\mathfrak{R}\left(\overline{z}\Gamma_{12}
\right)}{\det(\Gamma)}\right)
\int_{0}^{\infty}\frac{2}{\pi\textup{det}(\Gamma)r}
\textup{exp}\left(-
\frac{\Gamma_{22}r^{2}+\Gamma_{11}|z|^{2}r^{-2}}
{\det(\Gamma)}
\right)dr.
\end{align}
For the integral in (\ref{null_product2s}),
by the formula
\begin{align*}
\int_{0}^{\infty}r^{-1}e^{-ar^{2}-br^{-2}}dr=K_{0}\left(2\sqrt{ab}\right),\ a,b>0,
\end{align*}
where $K_{0}$ is the modified Bessel function of the second kind of order zero,
we have
\begin{align}\label{K0_apps}
&\int_{0}^{\infty}\frac{2}{\pi\textup{det}(\Gamma)r}
\textup{exp}\left(-
\frac{\Gamma_{22}r^{2}+\Gamma_{11}|z|^{2}r^{-2}}
{\det(\Gamma)}
\right)dr
=\frac{2}{\pi\textup{det}(\Gamma)}K_{0}\left(\frac{2\sqrt{\Gamma_{11}\Gamma_{22}}}{\det(\Gamma)}|z|\right).
\end{align}
By substituting  (\ref{K0_apps}) into (\ref{null_product2s}),
\begin{align}\label{null_product3s}
p\left(z;W_{\Phi}(t_{1},s_{1})\overline{W_{\Phi}(t_{2},s_{2})}\right)
=
\frac{2}{\pi\textup{det}(\Gamma)}\textup{exp}
\left(
\frac{2\mathfrak{R}\left(\overline{z}\Gamma_{12}
\right)}{\det(\Gamma)}\right)
K_{0}\left(\frac{2\sqrt{\Gamma_{11}\Gamma_{22}}}{\det(\Gamma)}|z|\right)
\end{align}
for any $z\in \mathbb{C}$.
By substituting (\ref{null_product3s}) into (\ref{start:proof:circov}),
\begin{align}\notag
&\mathbb{E}\left[e^{i(\Theta_{\Phi}(t_{1},s_{1})-\Theta_{\Phi}(t_{2},s_{2}))}\right]
\\\notag=&\int_{0}^{\infty}\int_{0}^{2\pi}e^{i\theta}p\left(re^{i\theta}; W_{\Phi}(t_{1},s_{1})\overline{W_{\Phi}(t_{2},s_{2})}\right)
rd\theta dr
\\\label{null_product4}=&
\int_{0}^{\infty}\int_{0}^{2\pi}e^{i\theta}
\left[
\frac{2}{\pi\textup{det}(\Gamma)}\textup{exp}
\left(
\frac{2\mathfrak{R}\left(re^{-i\theta}\Gamma_{12}
\right)}{\det(\Gamma)}\right)
K_{0}\left(\frac{2\sqrt{\Gamma_{11}\Gamma_{22}}}{\det(\Gamma)}r\right)
\right]
rd\theta dr.
\end{align}
For any $r>0$,
\begin{align}\notag
&\int_{0}^{2\pi}e^{i\theta}
\textup{exp}
\left(
\frac{2r|\Gamma_{12}|\mathfrak{R}\left(e^{-i(\theta-\arg\left(\Gamma_{12}\right))}
\right)}{\det(\Gamma)}\right)
d\theta
\\\notag=&\int_{0}^{2\pi}e^{i\theta}
\textup{exp}
\left(
\frac{2r|\Gamma_{12}|}
{\det(\Gamma)}\cos\left(\theta-\arg\left(\Gamma_{12}\right)\right)\right)
d\theta
\\\notag=&e^{i\arg\left(\Gamma_{12}\right)}\int_{0}^{2\pi}e^{iu}
\textup{exp}
\left(
\frac{2r|\Gamma_{12}|}
{\det(\Gamma)}\cos\left(u\right)\right)
du
\\\notag=& 2\pi e^{i\arg\left(\Gamma_{12}\right)}I_{1}\left(\frac{2r|\Gamma_{12}|}
{\det(\Gamma)}\right),
\end{align}
where $I_{1}$ is the modified Bessel function of the first kind of order one.
Hence, (\ref{null_product4}) becomes
\begin{align}\label{null_product4s}
&\mathbb{E}\left[e^{i(\Theta_{\Phi}(t_{1},s_{1})-\Theta_{\Phi}(t_{2},s_{2}))}\right]
\\\notag=&
\frac{2}{\pi\textup{det}(\Gamma)}
\left(2\pi e^{i\arg\left(\Gamma_{12}\right)}\right)
\int_{0}^{\infty}
I_{1}\left(\frac{2|\Gamma_{12}|}
{\det(\Gamma)}r\right)
K_{0}\left(\frac{2\sqrt{\Gamma_{11}\Gamma_{22}}}{\det(\Gamma)}r\right)
rdr
\\\notag=&
\frac{4e^{i\arg\left(\Gamma_{12}\right)}}{\textup{det}(\Gamma)}
\left(\frac{\det(\Gamma)}{2\sqrt{\Gamma_{11}\Gamma_{22}}}\right)^{2}
\int_{0}^{\infty}
I_{1}\left(\frac{|\Gamma_{12}|}{\sqrt{\Gamma_{11}\Gamma_{22}}}u
\right)
K_{0}\left(u\right)
udu
\\\notag=&
e^{i\arg\left(\Gamma_{12}\right)}
\frac{\det(\Gamma)}{\Gamma_{11}\Gamma_{22}}
\int_{0}^{\infty}
I_{1}\left(\frac{|\Gamma_{12}|}{\sqrt{\Gamma_{11}\Gamma_{22}}}u
\right)
K_{0}\left(u\right)
udu.
\end{align}
By the series representation (\ref{def:;Bessel}) of the modified Bessel function $I_{1}$, for any $c\in[0,1)$,
\begin{align}\notag
\int_{0}^{\infty}
I_{1}\left(cu
\right)
K_{0}\left(u\right)
udu
=&\int_{0}^{\infty}\left[\frac{1}{2}cu
\overset{\infty}{\underset{k=0}{\sum}}\frac{\left(\frac{1}{4}c^{2}u^{2}\right)^{k}}{k! \Gamma(k+2)}\right]
K_{0}\left(u\right)
udu
\\\notag=&
\frac{1}{2}c\overset{\infty}{\underset{k=0}{\sum}}\left(\frac{c^{2}}{4}\right)^{k}\frac{1}{k! \Gamma(k+2)}
\int_{0}^{\infty}u^{2k+2}K_{0}\left(u\right)
du
\\\notag=&
\frac{1}{2}c\overset{\infty}{\underset{k=0}{\sum}}\left(\frac{c^{2}}{4}\right)^{k}\frac{1}{k! \Gamma(k+2)}
\left[2^{2k+1}\Gamma\left(k+\frac{3}{2}\right)\Gamma\left(k+\frac{3}{2}\right)\right]
\\\notag=&
c\overset{\infty}{\underset{k=0}{\sum}}\left(c^{2}\right)^{k}\frac{1}{k! \Gamma(k+2)}
\left[\Gamma\left(k+\frac{3}{2}\right)\Gamma\left(k+\frac{3}{2}\right)\right]
\\\label{formula:F2}=&c\frac{\pi}{4}{}_2F_1\left(\frac{3}{2},\frac{3}{2};2;c^{2}\right).
\end{align}
By (\ref{null_product4s}) and (\ref{formula:F2}), the circular covariance between $\Theta_{\Phi}(t_{1},s_{1})$ and $\Theta_{\Phi}(t_{2},s_{2})$ has the form
\begin{align}\notag
\mathbb{E}\left[e^{i(\Theta_{\Phi}(t_{1},s_{1})-\Theta_{\Phi}(t_{2},s_{2}))}\right]
=&
\frac{\pi}{4}e^{i\arg\left(\Gamma_{12}\right)}
\frac{\det(\Gamma)}{\Gamma_{11}\Gamma_{22}}\frac{|\Gamma_{12}|}{\sqrt{\Gamma_{11}\Gamma_{22}}}
{}_2F_1\left(\frac{3}{2},\frac{3}{2};2;\frac{|\Gamma_{12}|^{2}}{\Gamma_{11}\Gamma_{22}}\right)
\\\notag=&\frac{\pi}{4}e^{i\arg\left(\Gamma_{12}\right)}
\left(\hspace{-0.13cm}1-\frac{|\Gamma_{12}|^{2}}{\Gamma_{11}\Gamma_{22}}\hspace{-0.05cm}\right)\hspace{-0.1cm}\frac{|\Gamma_{12}|}{\sqrt{\Gamma_{11}\Gamma_{22}}}
{}_2F_1\left(\hspace{-0.05cm}\frac{3}{2},\frac{3}{2};2;\frac{|\Gamma_{12}|^{2}}{\Gamma_{11}\Gamma_{22}}\hspace{-0.03cm}\right).
\end{align}
Finally, because the circular symmetry of $[W_{\Phi}(t_{1},s_{1})\ W_{\Phi}(t_{2},s_{2})]$ implies that the random vector
$$\left[e^{i\arg(W_{\Phi}(t_{1},s_{1}))}\ e^{i\arg(W_{\Phi}(t_{2},s_{2}))}\right]$$
is also circularly symmetric.
Hence, we have $\mathbb{E}\left[e^{i(\Theta_{\Phi}(t_{1},s_{1})+\Theta_{\Phi}(t_{2},s_{2}))}\right]=0$.
\qed

\subsection{Proof of Corollary \ref{prop:Cov_phaseAM_null}}\label{sec:proof:prop:Cov_phaseAM_null}
For any $r_{1}>0$ and $\theta_{2}\in[0,2\pi)$,
\begin{align}\label{phaseAM_pdf1}
&p\left(r_{1},\theta_{2};
|W_{\Phi}(t_{1},s_{1})|,\Theta_{\Phi}(t_{2},s_{2})\right)
\\\notag=&\int_{0}^{\infty}\int_{0}^{2\pi}p\left(r_{1},r_{2},\theta_{1},\theta_{2};
|W_{\Phi}(t_{1},s_{1})|,|W_{\Phi}(t_{2},s_{2})|,
\Theta_{\Phi}(t_{1},s_{1}),\Theta_{\Phi}(t_{2},s_{2})\right)
d\phi d\rho
\\\notag
=&\int_{0}^{\infty}\int_{0}^{2\pi}\hspace{-0.2cm}
\frac{r_{1}r_{2}}{\pi^2\textup{det}(\Gamma)}
\textup{exp}\left(-\begin{bmatrix}r_{1}e^{i\theta_{1}} \\ r_{2} e^{i\theta_{2}} \end{bmatrix}^{*}\hspace{-0.1cm}\Gamma^{-1}\hspace{-0.1cm}\left[\begin{array}{cc}r_{1}e^{i\theta_{1}}\\ r_{2}e^{i\theta_{2}}\end{array}\right]\right)
d\theta_1  dr_2
\\\notag=&
\int_{0}^{\infty}\int_{0}^{2\pi}\frac{r_{1}r_{2}}{\pi^2\textup{det}(\Gamma)}
\textup{exp}\left(-\frac{\Gamma_{22}r_{1}^2+\Gamma_{11}r_{2}^{2}}
{\textup{det}(\Gamma)}\right)
\textup{exp}\left(\frac{2r_{1}r_{2}\Re(e^{i(\theta_{1}-\theta_{2})}
\overline{\Gamma_{12}})}{\textup{det}(\Gamma)}\right)
 d\theta_1  dr_2.
\end{align}
Note that
\begin{align*}
\int_{0}^{2\pi}\hspace{-0.2cm}\textup{exp}\left(\frac{2r_{1}r_{2}\Re(e^{i(\theta_{1}-\theta_{2})}
\overline{\Gamma_{12}})}{\textup{det}(\Gamma)}\right)  d\theta_1
=2\pi I_{0}\left(\frac{2|\Gamma_{12}|r_{1}r_{2}}{\textup{det}(\Gamma)}\right),
\end{align*}
which does not change along with $\theta_{2}$.
Hence, (\ref{phaseAM_pdf1}) can be rewritten as
\begin{align}\notag
&p\left(r_{1},\theta_{2};
|W_{\Phi}(t_{1},s_{1})|,\Theta_{\Phi}(t_{2},s_{2})\right)
\\\notag=&\frac{2r_{1}}{\pi\textup{det}(\Gamma)}
\int_{0}^{\infty}\hspace{-0.25cm}r_{2}
\textup{exp}\left(\hspace{-0.1cm}-\frac{\Gamma_{22}r_{1}^2+\Gamma_{11}r_{2}^{2}}{\textup{det}(\Gamma)}\right)
I_{0}\hspace{-0.1cm}\left(\frac{2|\Gamma_{12}|r_{1}r_{2}}{\textup{det}(\Gamma)}\right)dr_2.
\end{align}
By using the formula \cite[p. 713]{gradshteyn2014table}
\begin{align*}
\int_{0}^{\infty} r e^{-ar^2}I_{0}(br)dr = \frac{1}{2a}\textup{exp}\left(\frac{b^2}{4a}\right),\ a,b>0,
\end{align*}
with $a =\Gamma_{11}/\textup{det}(\Gamma)$ and $b=2|\Gamma_{12}|r_{1}/\textup{det}(\Gamma)$,
\begin{align}\notag
&p\left(r_{1},\theta_{2};
|W_{\Phi}(t_{1},s_{1})|,\Theta_{\Phi}(t_{2},s_{2})\right)
\\\notag=&\frac{2r_{1}}{\pi\textup{det}(\Gamma)}\textup{exp}\left(-\frac{\Gamma_{22}r_{1}^2}{\textup{det}(\Gamma)}\right)
\left[\frac{\textup{det}(\Gamma)}{2\Gamma_{11}}\textup{exp}\left(\frac{|\Gamma_{12}|^2r_{1}^2}{\Gamma_{11}\textup{det}(\Gamma)}\right)
\right]
\\\notag=&\frac{r_{1}}{\pi\Gamma_{11}}
\textup{exp}\left(-\frac{\Gamma_{11}\Gamma_{22}-|\Gamma_{12}|^2}
{\Gamma_{11}\textup{det}(\Gamma)}r_{1}^{2}\right)
\\\notag=&
\frac{r_{1}}{\pi\Gamma_{11}}
\textup{exp}\left(-\frac{1}
{\Gamma_{11}}r_{1}^{2}\right).
\end{align}
From Corollaries \ref{prop:jpdf:S2} and \ref{corollary:phase:jpdf:null}, we know that
\begin{align}\notag
p\left(r_{1};
|W_{\Phi}(t_{1},s_{1})|\right)=\frac{2r_{1}}{\Gamma_{11}}\textup{exp}\left(-\frac{r_{1}^{2}}{\Gamma_{11}}\right),\ r_{1}\geq 0,
\end{align}
and $\Theta_{\Phi}(t_{2},s_{2})$ is uniformly distributed over $[0,2\pi).$
Therefore,
\begin{align}\notag
p\left(r_{1},\theta_{2};
|W_{\Phi}(t_{1},s_{1})|,\Theta_{\Phi}(t_{2},s_{2})\right)
=p\left(r_{1};
|W_{\Phi}(t_{1},s_{1})|\right)p\left(\theta_{2};
\Theta_{\Phi}(t_{2},s_{2})\right),
\end{align}
which implies that $|W_{\Phi}(t_{1},s_{1})|$ and $\Theta_{\Phi}(t_{2},s_{2})$ are independent.
\qed

\subsection{Proof of Proposition \ref{prop:approximation}}\label{sec:proof:prop:approximation}
By the Skorohod representation theorem \cite{wu2007strong}
(see also \cite[Theorem 3.1]{mies2023sequential} and \cite[Theorem 3.1]{wu2025uncertainty}),
on a sufficiently rich probability space, there exist random variables $\{\widehat{\epsilon}_{j}\}_{j\in \mathbb{Z}}$ and mean-zero independent Gaussian random
variables $\{z_{j}\}_{j\in \mathbb{Z}}$, defined on the same probability space,
such that $\{\widehat{\epsilon}_{j}\}_{j\in \mathbb{Z}}$ has the same distribution as $\{\epsilon_{j}\}_{j\in \mathbb{Z}}$ and $\textup{Var}(z_{j}) = \textup{Var}(\epsilon_{j})$ for all $j\in \mathbb{Z}$.
Let $\widetilde{W}_{\widehat{\epsilon}}$ denote the discretized AWT of $\{\widehat{\epsilon}_{j}\}_{j\in \mathbb{Z}}$,
and
denote the probability measure of
$\{\widetilde{W}_{\widehat{\epsilon}}(t_{i},s_{\ell})\mid 1\leq i\leq n,\ 1\leq \ell \leq d\}$
by $\nu_{\widehat{\epsilon}}^{(n)}$.
We have
\begin{align}\label{epsilon_D_proxy}
\nu_{\epsilon}^{(n)}=\nu_{\widehat{\epsilon}}^{(n)}.
\end{align}

Define
$$I_{n} = [-\lceil \sqrt{n}As_{d}\rceil+1, \lceil \sqrt{n}As_{d}\rceil+n]\cap \mathbb{Z}.$$
For $h\in I_{n}$,
define the partial sum
\begin{align*}
S_{\widehat{\epsilon}}(h) = \overset{h}{\underset{k=-\lceil \sqrt{n}As_{d}\rceil+1}{\sum}}\widehat{\epsilon}_{k}
\end{align*}
with $S_{\widehat{\epsilon}}(h)=0$ when $h\leq -\lceil \sqrt{n}As_{d}\rceil$.
By (\ref{def:discreteCWT2}) and the summation by parts,
for any $1\leq i\leq n$ and $1\leq \ell \leq d$,
\begin{align}\notag
\widetilde{W}_{\widehat{\epsilon}}(t_{i},s_{\ell})
=&\frac{1}{n^{1/4}}\underset{j\in D^{(\ell)}}{\sum}\left(S_{\widehat{\epsilon}}(i+j)-S_{\widehat{\epsilon}}(i+j-1)\right)\overline{\psi_{j}^{(\ell)}}
\\\notag=&\frac{1}{n^{1/4}}S_{\widehat{\epsilon}}\left(i+\lceil \sqrt{n}As_{\ell}\rceil\right) \overline{\psi_{\lceil \sqrt{n}As_{\ell}\rceil+1}^{(\ell)}}
\\\notag&-\frac{1}{n^{1/4}}S_{\widehat{\epsilon}}\left(i-\lceil \sqrt{n}As_{\ell}\rceil-1\right)\overline{\psi_{-\lceil \sqrt{n}As_{\ell}\rceil}^{(\ell)}}
\\\label{def:discreteCWT3}&-\frac{1}{n^{1/4}}\underset{j\in D^{(\ell)}}{\sum}S_{\widehat{\epsilon}}\left(i+j\right)\left(\ \overline{\psi^{(\ell)}_{j+1}}-\overline{\psi^{(\ell)}_{j}}\ \right).
\end{align}
Similarly, the discretized AWT of the Gaussian sequence $\{z_{j}\}_{j\in \mathbb{Z}}$ can be expressed as
\begin{align}\notag
\hspace{-0.3cm}\widetilde{W}_{z}(t_{i},s_{\ell})
=& \frac{1}{n^{1/4}}S_{z}\left(i+\lceil \sqrt{n}As_{\ell}\rceil\right)\overline{\psi_{\lceil \sqrt{n}As_{\ell}\rceil+1}^{(\ell)}}
\\\notag&-\frac{1}{n^{1/4}}S_{z}\left(i-\lceil \sqrt{n}As_{\ell}\rceil-1\right)\overline{\psi_{-\lceil \sqrt{n}As_{\ell}\rceil}^{(\ell)}}
\\\label{def:discreteCWT_z}&-\frac{1}{n^{1/4}}\underset{j\in D^{(\ell)}}{\sum}S_{z}\left(i+j\right)\left(\ \overline{\psi^{(\ell)}_{j+1}}-\overline{\psi^{(\ell)}_{j}}\ \right),
\end{align}
where
\begin{align*}
S_{z}(h) = \overset{h}{\underset{k=-\lceil \sqrt{n}As_{d}\rceil+1 }{\sum}}z_{k},\ h\in I_{n},
\end{align*}
and $S_{z}(h)=0$ for $h\leq-\lceil \sqrt{n}As_{d}\rceil$.
By (\ref{def:discreteCWT3}) and (\ref{def:discreteCWT_z}),
\begin{align}\notag
&n^{1/4}\left(\widetilde{W}_{\widehat{\epsilon}}(t_{i},s_{\ell}) -\widetilde{W}_{z}(t_{i},s_{\ell})\right)
\\\notag=&
\left[S_{\widehat{\epsilon}}\left(i+\lceil \sqrt{n}As_{\ell}\rceil\right)- S_{z}\left(i+\lceil \sqrt{n}As_{\ell}\rceil\right)\right]\overline{\psi_{\lceil \sqrt{n}As_{\ell}\rceil+1}^{(\ell)}}
\\\notag&-
\left[S_{\widehat{\epsilon}}\left(i-\hspace{-0.1cm}\lceil \sqrt{n}As_{\ell}\rceil\hspace{-0.1cm}-\hspace{-0.05cm}1\right)-S_{z}\left(i-\hspace{-0.1cm}\lceil \sqrt{n}As_{\ell}\rceil\hspace{-0.1cm}-\hspace{-0.05cm}1\right)\right]\hspace{-0.05cm}
\overline{\psi_{-\lceil \sqrt{n}As_{\ell}\rceil}^{(\ell)}}
\\\notag&-\underset{j\in D^{(\ell)}}{\sum}\left[ S_{\widehat{\epsilon}}\left(i+j\right) -S_{z}\left(i+j\right)\right]\left(\ \overline{\psi^{(\ell)}_{j+1}}-\overline{\psi^{(\ell)}_{j}}\ \right).
\end{align}
From (\ref{def:discrete_psi}), both $|\psi_{\lceil \sqrt{n}As_{\ell}\rceil+1}^{(\ell)}|$ and $|\psi_{-\lceil \sqrt{n}As_{\ell}\rceil}^{(\ell)}|$
are bounded by $n^{-1/4}s_{\ell}^{-1/2}\|\psi\|_{L^{\infty}}$. Hence,
\begin{align}\notag
\left|\widetilde{W}_{\widehat{\epsilon}}(t_{i},s_{\ell}) -\widetilde{W}_{z}(t_{i},s_{\ell})\right|
\leq&
\frac{\|\psi\|_{L^{\infty}}}{\sqrt{n s_{\ell}}}\left|S_{\widehat{\epsilon}}\left(i+\lceil \sqrt{n}As_{\ell}\rceil \right)- S_{z}\left(i+\lceil \sqrt{n}As_{\ell}\rceil\right)\right|
\\\notag&+\frac{\|\psi\|_{L^{\infty}}}{\sqrt{ ns_{\ell}}}\left|S_{\widehat{\epsilon}}\left(i-\lceil \sqrt{n}As_{\ell}\rceil-1\right)-S_{z}\left(i-\lceil \sqrt{n}As_{\ell}\rceil-1\right)\right|
\\\label{diff:WzWepsilon}&+\frac{1}{n^{1/4}}
(\underset{j\in D^{(\ell)}}{\max}\left|S_{\widehat{\epsilon}}\left(i+j\right) -S_{z}\left(i+j\right)\right|)
\underset{j\in D^{(\ell)}}{\sum}\left|\psi^{(\ell)}_{j+1}-\psi^{(\ell)}_{j}\right|.
\end{align}
For the third term on the right-hand side of (\ref{diff:WzWepsilon}), if $\psi$ is absolutely continuous, then its derivative $\psi'$ exists almost everywhere, and
\begin{align}\notag
\underset{j\in D^{(\ell)}}{\sum}\left|\psi^{(\ell)}_{j+1}-\psi^{(\ell)}_{j}\right|
=&\overset{\lceil \sqrt{n}As_{\ell}\rceil}{\underset{j=-\lceil \sqrt{n}As_{\ell}\rceil}{\sum}}
\left|\frac{1}{n^{1/4}s_{\ell}^{1/2}}\psi\left(\frac{j+1}{\sqrt{n}s_{\ell}}\right)-\frac{1}{n^{1/4}s_{\ell}^{1/2}}
\psi\left(\frac{j}{\sqrt{n}s_{\ell}}\right)\right|
\\\notag=&
\frac{1}{n^{1/4}s_{\ell}^{1/2}}\overset{\lceil \sqrt{n}As_{\ell}\rceil}{\underset{j=-\lceil \sqrt{n}As_{\ell}\rceil}{\sum}}
\left|\int_{jn^{-1/2}s_{\ell}^{-1}}^{(j+1)n^{-1/2}s_{\ell}^{-1}}
\psi'(\tau)d\tau\right|
\\\label{upperbound_sum_psi}\leq&\frac{1}{n^{1/4}s_{\ell}^{1/2}}
\|\psi'\|_{L^{1}}.
\end{align}
By substituting (\ref{upperbound_sum_psi}) into (\ref{diff:WzWepsilon}),
\begin{align}\notag
\sqrt{s_{\ell}}\left|\widetilde{W}_{\widehat{\epsilon}}(t_{i},s_{\ell}) -\widetilde{W}_{z}(t_{i},s_{\ell})\right|
\leq& \frac{1}{\sqrt{n}}\left(\|\psi'\|_{L^{1}}+2\|\psi\|_{L^{\infty}}\right)
\\\label{diff:WzWepsilon2}&\left(\underset{-\lceil \sqrt{n}As_{\ell}\rceil-1 \leq j\leq \lceil \sqrt{n}As_{\ell}\rceil}{\max}\left|S_{\widehat{\epsilon}}\left(i+j\right) -S_{z}\left(i+j\right)\right|\right).
\end{align}
Since $D^{(1)}\subset D^{(2)}\subset\cdots\subset D^{(d)}$,
\begin{align}\notag
&\underset{1\leq i \leq n}{\max}\ \underset{1\leq \ell \leq d}{\max}\ \underset{-\lceil \sqrt{n}As_{\ell}\rceil-1 \leq j\leq \lceil \sqrt{n}As_{\ell}\rceil}{\max}\left|S_{\widehat{\epsilon}}\left(i+j\right) -S_{z}\left(i+j\right)\right|
\\\notag=&\underset{-\lceil \sqrt{n}As_{d}\rceil \leq h\leq \lceil \sqrt{n}As_{d}\rceil+n}{\max}\left|S_{\widehat{\epsilon}}\left(h\right) -S_{z}\left(h\right)\right|
\\\label{max_diff_S}=&\underset{h\in I_{n}}{\max}\left|S_{\widehat{\epsilon}}\left(h\right) -S_{z}\left(h\right)\right|,
\end{align}
where the last equality follows from that
$S_{\widehat{\epsilon}}(-\lceil \sqrt{n}As_{d}\rceil)=S_{z}(-\lceil \sqrt{n}As_{d}\rceil)=0$.
By (\ref{diff:WzWepsilon2}) and (\ref{max_diff_S}), for any $\varepsilon> 0$,
\begin{align}\notag
&\mathbb{P}\left(\underset{1\leq \ell \leq d}{\max}\ \underset{1\leq i\leq n}{\max}\sqrt{s_{\ell}}\left|\widetilde{W}_{\widehat{\epsilon}}(t_{i},s_{\ell})
-\widetilde{W}_{z}(t_{i},s_{\ell})\right|>\varepsilon\right)
\\\notag\leq&
\mathbb{P}\left( \frac{1}{\sqrt{n}}\left(\|\psi'\|_{L^{1}}+2\|\psi\|_{L^{\infty}}\right)\underset{h\in I_{n}}{\max}\left|S_{\widehat{\epsilon}}\left(h\right) -S_{z}\left(h\right)\right|>\varepsilon\right)
\\\label{control:1}=&
\mathbb{P}\left( \underset{h\in I_{n}}{\max}\left|S_{\widehat{\epsilon}}\left(h\right) -S_{z}\left(h\right)\right|>\frac{\varepsilon \sqrt{n}}{\|\psi'\|_{L^{1}}+2\|\psi\|_{L^{\infty}}}\right).
\end{align}
Under Assumptions \ref{assumption:causal} and \ref{assumption:weak}, Theorem 3.1 in \cite{wu2025uncertainty}
shows that
\begin{align}\label{cited_result:Wu}
\mathbb{E}\left[\underset{h\in I_{n}}{\max}\left|
\left(S_{\widehat{\epsilon}}\left(h\right) -S_{z}\left(h\right)\right)\right|^{2}\right]
=\mathcal{O}\left(|I_{n}|^{1-G(p,\chi)}\left(\log|I_{n}|\right)^{4}\right),
\end{align}
where the cardinality of $I_{n}$, denoted by $|I_{n}|$, is equal to $2\lceil \sqrt{n}As_{d}\rceil+n$.
By applying the Markov inequality to (\ref{control:1}),
\begin{align}\notag
&\mathbb{P}\left(\underset{1\leq \ell \leq d}{\max}\ \underset{1\leq i\leq n}{\max}\sqrt{s_{\ell}}\left|\widetilde{W}_{\widehat{\epsilon}}(t_{i},s_{\ell})
-\widetilde{W}_{z}(t_{i},s_{\ell})\right|>\varepsilon\right)
\\\notag=&
\left(\|\psi'\|_{L^{1}}+2\|\psi\|_{L^{\infty}}\right)^{2}\varepsilon^{-2}n^{-1}\mathcal{O}\left(|I_{n}|^{1-G(p,\chi)}\left(\log|I_{n}|\right)^{4}\right)
\\\label{control_Markov}=&
\varepsilon^{-2}\mathcal{O}\left(n^{-G(p,\chi)}\left(\log n\right)^{4}\right).
\end{align}
By choosing
\begin{align*}
\varepsilon =\varepsilon_{n}:= n^{-G(p,\chi)/3}(\log n)^{4/3+\delta},
\end{align*}
where $\delta$ is an arbitrary positive number, (\ref{control_Markov}) becomes
\begin{align}\label{control:3}
\mathbb{P}\left(\underset{\begin{subarray}{c}
1\leq i\leq n \\ 1\leq \ell \leq d\end{subarray}}{\max}\sqrt{s_{\ell}}|\widetilde{W}_{\widehat{\epsilon}}(t_{i},s_{\ell})
-\widetilde{W}_{z}(t_{i},s_{\ell})|>\varepsilon_{n}\hspace{-0.15cm}\right)
=o\left(\varepsilon_{n}\right)
\end{align}
when $n\rightarrow\infty$.
Denote
$$
\widetilde{\mathbf{W}}_{\square}=\{\sqrt{s_{\ell}}\widetilde{W}_{\square}(t_{i},s_{\ell})\mid 1\leq i\leq n,\ 1\leq \ell \leq d\},
$$
where $\square\in \{\epsilon,\widehat{\epsilon},z\}$.
For any $B\in \mathcal{B}(\mathbb{C}^{n\times d})$, by (\ref{epsilon_D_proxy}),
\begin{align}\label{control:4}
\nu_{\epsilon}^{(n)}\left(B\right) =& \mathbb{P}\left(\widetilde{\mathbf{W}}_{\epsilon}\in B\right)=\mathbb{P}\left(\widetilde{\mathbf{W}}_{\widehat{\epsilon}}\in B\right)
\\\notag=&
\mathbb{P}\left(\widetilde{\mathbf{W}}_{\widehat{\epsilon}}\in B, |\widetilde{\mathbf{W}}_{\widehat{\epsilon}}-\widetilde{\mathbf{W}}_{z}|_{\infty}\leq \varepsilon_{n}\right)
+
\mathbb{P}\left(\widetilde{\mathbf{W}}_{\widehat{\epsilon}}\in B, |\widetilde{\mathbf{W}}_{\widehat{\epsilon}}-\widetilde{\mathbf{W}}_{z}|_{\infty}> \varepsilon_{n}\right)
\\\notag\leq&
\mathbb{P}\left(\widetilde{\mathbf{W}}_{z}\in B^{\varepsilon_{n}}\right)
+\mathbb{P}\left( |\widetilde{\mathbf{W}}_{\widehat{\epsilon}}-\widetilde{\mathbf{W}}_{z}|_{\infty}> \varepsilon_{n}\right).
\end{align}
By (\ref{control:3}) and (\ref{control:4}),
\begin{align}\notag
\nu_{\epsilon}^{(n)}\left(B\right) <&
\nu_{z}^{(n)}\left(B^{\varepsilon_{n}}\right)+ o\left(\varepsilon_{n}\right)
\end{align}
when $n$ is sufficiently large. Similarly,
$
\nu_{z}^{(n)}\left(B\right) <
\nu_{\epsilon}^{(n)}\left(B^{\varepsilon_{n}}\right)+ o\left(\varepsilon_{n}\right)
$
for all sufficiently large $n$.
Therefore,
the L$\acute{\textup{e}}$vy-Prokhorov distance between $\nu_{z}^{(n)}$ and $\nu_{\epsilon}^{(n)}$
is bounded by $o\left(\varepsilon_{n}\right)$ when $n$ is sufficiently large.
\qed

\begin{Remark}\label{remark:approximation}
By (\ref{diff:WzWepsilon2}),
the Euclidean norm in $\mathbb{R}^{d}$ of $\widetilde{\mathbf{W}}_{\widehat{\epsilon}}(t_{i})-\widetilde{\mathbf{W}}_{z}(t_{i})$
is bounded as follows.
\begin{align}\notag
&\left|\widetilde{\mathbf{W}}_{\widehat{\epsilon}}(t_{i})-\widetilde{\mathbf{W}}_{z}(t_{i})\right|^{2}
=
\overset{d}{\underset{\ell=1}{\sum}}
\left|\widetilde{W}_{\widehat{\epsilon}}(t_{i},s_{\ell}) -\widetilde{W}_{z}(t_{i},s_{\ell})\right|^{2}
\\\notag\leq& \left(\|\psi'\|_{L^{1}}+2\|\psi\|_{L^{\infty}}\right)^{2}n^{-1}
\overset{d}{\underset{\ell=1}{\sum}}s_{\ell}^{-1}\hspace{-0.1cm}\underset{-\lceil \sqrt{n}As_{\ell}\rceil-1 \leq j\leq \lceil \sqrt{n}As_{\ell}\rceil}{\max}\left|S_{\widehat{\epsilon}}\left(i+j\right) -S_{z}\left(i+j\right)\right|^{2},
\end{align}
which implies that
\begin{align}\notag
&\left(\|\psi'\|_{L^{1}}+2\|\psi\|_{L^{\infty}}\right)^{-2}n\underset{1\leq i\leq n}{\max}\ \left|\widetilde{\mathbf{W}}_{\widehat{\epsilon}}(t_{i})-\widetilde{\mathbf{W}}_{z}(t_{i})\right|^{2}
\\\notag
\leq&\overset{d}{\underset{\ell=1}{\sum}}s_{\ell}^{-1}\ \underset{-\lceil \sqrt{n}As_{\ell}\rceil \leq h\leq \lceil \sqrt{n}As_{\ell}\rceil+n}{\max}\left|S_{\widehat{\epsilon}}\left(h\right) -S_{z}\left(h\right)\right|^{2}.
\end{align}
Hence,
\begin{align}\notag
&\left(\|\psi'\|_{L^{1}}+2\|\psi\|_{L^{\infty}}\right)^{-2}n\mathbb{E}\left[\underset{1\leq i\leq n}{\max}\ \left|\widetilde{\mathbf{W}}_{\widehat{\epsilon}}(t_{i})-\widetilde{\mathbf{W}}_{z}(t_{i})\right|^{2}\right]
\\\notag\leq&\overset{d}{\underset{\ell=1}{\sum}}s_{\ell}^{-1}\mathbb{E}\left[\underset{-\lceil \sqrt{n}As_{\ell}\rceil \leq h\leq \lceil \sqrt{n}As_{\ell}\rceil+n}{\max}\left|S_{\widehat{\epsilon}}\left(h\right) -S_{z}\left(h\right)\right|^{2}\right].
\end{align}
By applying (\ref{cited_result:Wu})
to the inequality above,
\begin{align}\notag
&n\mathbb{E}\left[\underset{1\leq i\leq n}{\max}\left|\widetilde{\mathbf{W}}_{\widehat{\epsilon}}(t_{i})-\widetilde{\mathbf{W}}_{z}(t_{i})\right|^{2}\right]
\\\notag=&\overset{d}{\underset{\ell=1}{\sum}}\frac{1}{s_{\ell}}\
\mathcal{O}\left(\left(2\lceil \sqrt{n}As_{\ell}\rceil+n\right)^{1-G(p,\chi)}\left(\log\left(2\lceil \sqrt{n}As_{\ell}\rceil+n\right)\right)^{4}\right),
\end{align}
which implies that
\begin{align}\notag
\mathbb{E}\left[\underset{1\leq i\leq n}{\max}\left|\widetilde{\mathbf{W}}_{\widehat{\epsilon}}(t_{i})-\widetilde{\mathbf{W}}_{z}(t_{i})\right|^{2}\right]
=d\mathcal{O}\left(n^{-G(p,\chi)}\left(\log\left(n\right)\right)^{4}\right)
\end{align}
if the scale values $\{s_{\ell}\}_{1\leq \ell\leq d}$ are bounded.
\end{Remark}

\bibliographystyle{IEEEtran}
\bibliography{ref_nonnull2027}

\end{document}